\RequirePackage{pdf14} 
\documentclass[
a4paper,
12pt,
twoside,
abstracton,
headsepline,
chapterprefix,
titlepage=false,
cleardoublepage=empty
]{scrreprt} 

\usepackage{texFiles/style/thesisStyle}
\usepackage{datetime}
\usepackage{gitinfo2}

\newcommand{\mthauthor}{Nils-Arne Dreier}%
\newcommand{\mthhometown}{Bünde}%
\newcommand{\mthdate}{2020}%
\newcommand{\mthtitle}{Hardware-Oriented Krylov Methods for High-Performance Computing}%
\newcommand{\mthsubtitle}{}%
\newcommand{\mthlogo}{images/logos/logowwu}
\newcommand{\mththesistitle}{Inaugural Dissertation}
\newcommand{\mththesisgraduation}{zur Erlangung des Doktorgrades der Naturwissenschaften}
\newcommand{\mththesisdegree}{Dr. rer. nat.}
\newcommand{\mthdepartmentI}{im Fachbereich Mathematik und Informatik}
\newcommand{\mthdepartmentII}{der Mathematisch-Naturwissenschaftlichen Fakultät}
\newcommand{\mthuni}{der Westfälischen Wilhelms-Universität Münster}

\newcommand{\mthdean}{Prof.~Dr.~Xiaoyi Jiang}%
\newcommand{\mthdeanaff}{Westfälische Wilhelms-Universität Münster}%
\newcommand{\mthdeancity}{Münster, DE}%
\newcommand{\mthreviewerA}{Prof.~Dr.~Christian Engwer}%
\newcommand{\mthreviewerAaff}{Westfälische Wilhelms-Universität Münster}%
\newcommand{\mthreviewerAcity}{Münster, DE}%
\newcommand{\mthreviewerB}{Laura Grigori, PhD}%
\newcommand{\mthreviewerBaff}{INRIA Paris}%
\newcommand{\mthreviewerBcity}{Paris, FR}%
\newcommand{\bb}[1]{\mathbb{#1}}

\newcommand{\bbR}{\bb{R}}
\newcommand{\bbP}{\bb{P}}
\newcommand{\krylov}[4][]{\mathcal{K}_{#1}^{#4}\left(#2,#3\right)}
\newcommand{\SubA}[1][]{\mathbb{S}_{#1}}

\newcommand{\Identity}[1][]{\mathbb{I}_{#1}}
\newcommand{\Prec}{M}

\newcommand{\dimA}{n}
\newcommand{\solution}{x^*}
\newcommand{\blocksolution}{X^*}

\newcommand{\cond}[1][]{\kappa\def\temp{#1}\ifx\temp\empty\else\left(#1\right)\fi}
\newcommand{\polyX}{\mathrm{x}}

\newcommand{\vspan}[1][]{\operatorname{span}\def\temp{#1}\ifx\temp\empty\else\left(#1\right)\fi}

\newcommand{\diag}[1]{\operatorname{diag}\left({#1}\right)}

\newcommand{\normalizer}[2][]{\operatorname{Norm}_{#2}\def\temp{#1}\ifx\temp\empty\else\left(#1\right)\fi}

\newcommand{\transpose}[1]{{#1}^{\mathsf{T}}}

\newcommand{\inverse}[1]{{#1}^{-1}}
\newcommand{\transinv}[1]{{#1}^{\mathsf{-T}}}

\renewcommand{\tilde}[1]{\widetilde{#1}}

\newcommand{\machineeps}{\varepsilon_{\text{mach}}}
\newcommand{\tol}{\varepsilon_{\text{tol}}}

\newcommand{\mcJ}{\mathcal{J}}
\newcommand{\sproduct}[3][]{\langle #2, #3 \rangle_{#1}}
\newcommand{\norm}[2][]{\|#2\|_{#1}}
\newcommand{\trace}[1]{\operatorname{tr}\left(#1\right)}
\newcommand{\Cheby}[2][]{T^{#2}\def\temp{#1}\ifx\temp\empty\else\left(#1\right)\fi}
\newcommand{\SCheby}[2][]{\widetilde{T}_{\lambda_1,\lambda_\dimA}^{#2}\def\temp{#1}\ifx\temp\empty\else\left(#1\right)\fi}
\newcommand{\dune}[1][]{\textsc{Dune}\ifthenelse{\equal{#1}{}}{}{\textsc{-{#1}}}\xspace}
\newcommand{\cc}{C\texttt{++}}
\DeclareMathOperator*{\argmin}{arg\,min}

\makeatletter
\newsavebox{\@brx}
\newcommand{\llangle}[1][]{\savebox{\@brx}{\(\m@th{#1\langle}\)}%
  \mathopen{\copy\@brx\kern-0.5\wd\@brx\usebox{\@brx}}}
\newcommand{\rrangle}[1][]{\savebox{\@brx}{\(\m@th{#1\rangle}\)}%
  \mathclose{\copy\@brx\kern-0.5\wd\@brx\usebox{\@brx}}}
\makeatother
\begin{document}%
\sloppy
\pagenumbering{roman}

\hypersetup{pageanchor=false}
\includepdf[pages={1}]{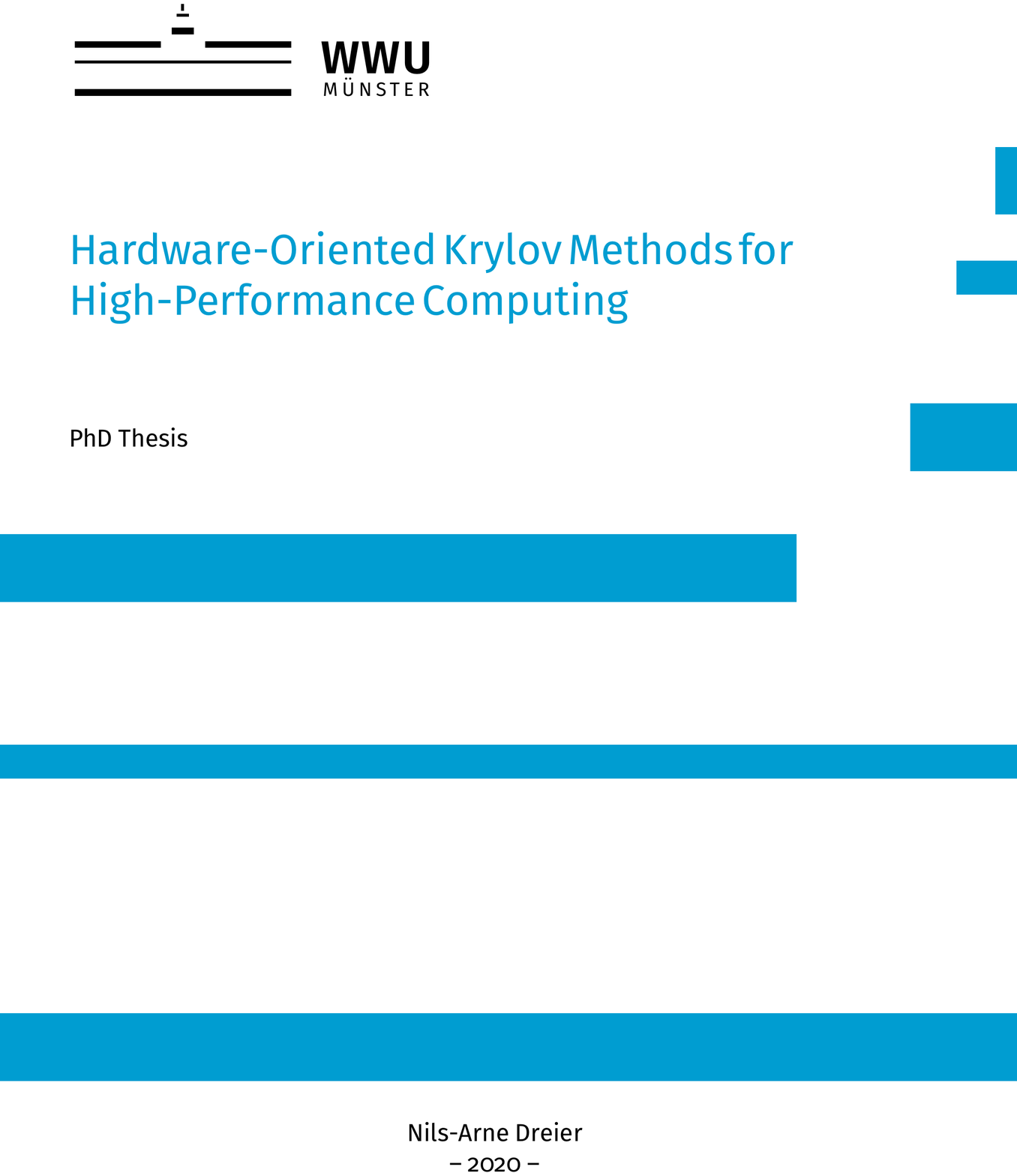}
\graphicspath{{..}}
\subject{\includestandalone{\mthlogo}\vspace{1cm}\\{\LARGE\textbf{Fach: Mathematik}}}%
\title{\color{mainColor}\vspace{0.3cm}\mthtitle}%
\subtitle{\color{mainColor}\mthsubtitle}%
\author{\vspace{0.55cm}\\{\LARGE\textbf{\color{mainColor}\mththesistitle}}\\\mththesisgraduation\\\ \\\textendash\ \mththesisdegree\ \textendash\\\ \\\mthdepartmentI\\\mthdepartmentII\\\mthuni}%
\date{}%
\publishers{{\small eingereicht von\\}{\color{mainColor}\textbf\mthauthor}\\{\small aus\\}{\color{mainColor}\textbf\mthhometown}\\\textendash\ \mthdate\ \textendash}%
{\singlespacing\maketitle[3]}%
\graphicspath{{images/}}
\thispagestyle{empty}\vspace*{\fill}
\begin{tabularx}{\textwidth}{l @{\extracolsep{\fill}} r}\hline
	\ &\ \\
	Dekan: & \mthdean\\
	\multicolumn{2}{r}{\small \mthdeanaff}\\
	\multicolumn{2}{r}{\small \mthdeancity}\\
	\ &\ \\
	Erster Gutachter: & \mthreviewerA\\
	\multicolumn{2}{r}{\small \mthreviewerAaff}\\
	\multicolumn{2}{r}{\small \mthreviewerAcity}\\
	\ &\ \\
	Zweiter Gutachter: & \mthreviewerB\\
	\multicolumn{2}{r}{\small \mthreviewerBaff}\\
	\multicolumn{2}{r}{\small \mthreviewerBcity}\\
	\ &\ \\
	Tag der mündlichen Prüfung: & 08.03.2021\\
	\ &\ \\
	Tag der Promotion: & 08.03.2021\\
	\ &\ \\\hline
\end{tabularx}
%
\chapter*{Abstract}

Krylov subspace methods are an essential building block in numerical simulation
software.
The efficient utilization of modern hardware is a challenging problem in the
development of these methods.
In this work, we develop Krylov subspace methods to solve linear systems with
multiple right-hand sides, tailored to modern hardware in high-performance
computing.

To this end, we analyze an innovative block Krylov subspace framework that
allows to balance the computational and data-transfer costs to the hardware.
Based on the framework, we formulate commonly used Krylov methods.
For the CG and BiCGStab methods, we introduce a novel stabilization approach as
an alternative to a deflation strategy.
This helps us to retain the block size, thus leading to a simpler and more
efficient implementation.

In addition, we optimize the methods further for distributed memory systems and
the communication overhead.
For the CG method, we analyze approaches to overlap the communication and
computation and present multiple variants of the CG method, which
differ in their communication properties.
Furthermore, we present optimizations of the orthogonalization procedure in the
GMRes method.
Beside introducing a pipelined Gram-Schmidt variant that overlaps the global
communication with the computation of inner products, we present a novel
orthonormalization method based on the TSQR algorithm, which is
communication-optimal and stable.
For all optimized method, we present tests that show their superiority in a
distributed setting.

\begin{otherlanguage}{ngerman}
\chapter*{Zusammenfassung}

Krylovraummethoden stellen einen essentiellen Bestandteil numerischer
Simulationssoftware dar.
Die effiziente Nutzung moderner Hardware ist ein herausforderndes
Problem bei der Entwicklung solcher Methoden.
Gegenstand dieser Dissertation ist die Formulierung von Krylovraumverfahren zur
Lösung von linearen Gleichungssystemen mit mehreren rechten Seiten, welche die
Eigenschaften moderner Hardware berücksichtigen.

Dazu untersuchen wir ein innovatives Blockkrylovraum-Framework, welches es
ermöglicht die Berechnungs- und Datentransferkosten der Blockkrylovraummethode
an die Hardware anzupassen.
Darauf aufbauend formulieren wir mehrere Krylovraummethoden.
Für die CG und BiCGStab Methoden führen wir eine neuartige
Stabilisierungstrategie ein, die es ermöglicht die Spaltenanzahl des Residuums
beizubehalten.
Diese ersetzt die bekannte Deflationstrategien und ermöglicht eine einfachere
und effizientere Implementierung der Methoden.

Des Weiteren optimieren wir die Methoden bezüglich der Kommunikation auf
Systemen mit verteiltem Speicher.
Für die CG Methode untersuchen wir Strategien, um die Kommunikation mit
Berechnungen zu überlappen.
Dazu stellen wir mehrere Varianten des Algorithmus vor, welche sich durch ihre
Kommunikationseigenschaften unterscheiden.
Außerdem werden für die GMRes Methode optimierte Varianten der
Orthonormalisierung entwickeln.
Neben einem Gram-Schmidt Verfahren, welches Berechnungen und Kommunikation
überlappt, präsentieren wir eine neue Methode, welche auf dem TSQR-Algorithmus
aufbaut und Stabilität sowie geringe Kommunikationskosten vereint.
Für alle optimierten Varianten zeigen wir numerische Tests, welche die
Verbesserungen auf Systemen mit verteiltem Speicher demonstrieren.

\end{otherlanguage}



%
\cleardoublepage%
\makeackn{

  I would like to express my deep gratitude to all people who have supported me
over the last years.
  First, I thank Prof.~Dr.~Christian Engwer for giving me the opportunity to work on
this topic, all the creative discussions, motivation, great guidance and for
being an excellent supervisor.
  I thank all my colleges in our workgroup for the pleasant atmosphere, in
particular I thank Liesel Sommer and Marcel Koch for proof-reading this thesis
and giving useful hints for improvements.
  Furthermore, I thank Prof.\ Dr.\ Robert Klöfkorn for giving me the opportunity to work
for a few weeks in Bergen, collecting valuable experience and enjoying the
Norwegian nature.
  All implementations of algorithms in the thesis are based on the \dune
software framework, thus I thank all the developers of \dune for making this
project happen.

\begin{otherlanguage}{ngerman}

  Diese Arbeit wäre ohne die bedingungslose Unterstützung meiner Eltern Kirsten
und Eckhard nicht möglich gewesen.
  Danke für all die finanzielle und moralische Unterstützung während meiner
gesamten Studienzeit.

Der größte Dank gilt meiner Frau Eileen.
Danke dafür, dass es dich in meinem Leben gibt und für all den Rückhalt, die
Unterstützung und Liebe, die es mir sehr erleichtert haben diese Arbeit zu
verfassen.

\end{otherlanguage}
}


%
\cleardoublepage
\hypersetup{pageanchor=true}
\phantomsection
\setcounter{tocdepth}{1}%
{\color{mainColor}\tableofcontents}%
\thispagestyle{empty}%

\clearpage%
\listoffigures%
%
\listoftables%
%
\listofalgorithms%
%
\lstlistoflistings%

\newpage
\thispagestyle{empty}
\hbox{}
\cleardoublepage%
\pagenumbering{arabic}%
\chapter{Introduction} \epigraph{\itshape We can only see a short distance
ahead, but we can see plenty there that needs to be done.}{\scshape Alan Turing}
\section{Motivation}
In the last decades, High-Performance Computing (HPC) became an
essential part of science, industry and our every day life.
Engineers use it to optimize the shape of cars and aircraft.
Meteorologists use it to create the daily weather forecast.
Physicists use it to simulate quantum mechanics which
helps to understand the elements of our universe.
It is widely used to simulate the global climate on large supercomputers.
Petroleum engineers design offshore platforms using
HPC to make them more efficient.
In medical research, scientists simulate an entire heart
or brain to investigate the sources of strokes and heart attacks.
Other fields of application include sociology, biology and astrology.
Even during the COVID-19 pandemic, HPC is used to investigate medicine
that is effective to treat COVID-19 patients.

In all these applications, HPC brings great improvements.
As a consequence, the need for more and more computation power has grown
extremely.
\begin{figure}[ht]
  \centering
  \includegraphics[width=0.75\textwidth]{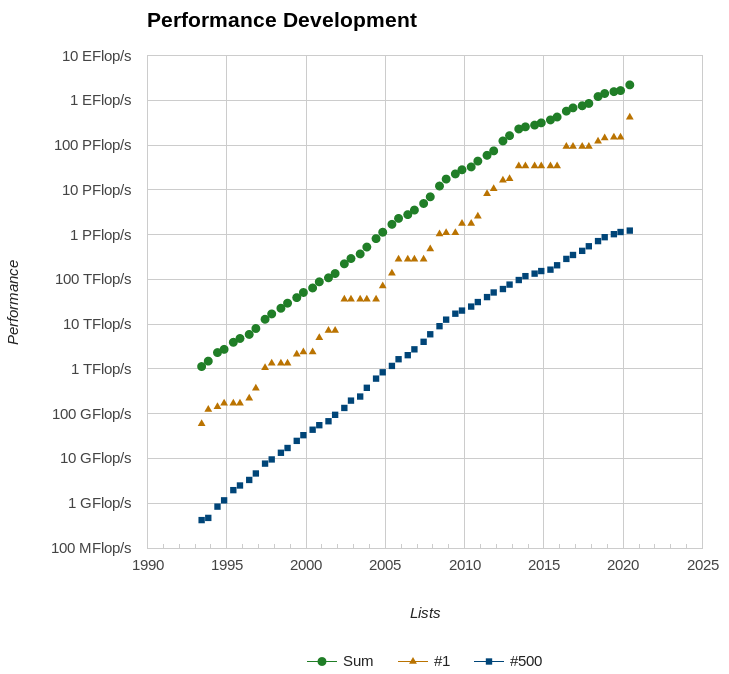}
  \caption[Performance development of the top 500
    supercomputers.]{Performance development of the top 500
    supercomputers.\\Taken from \url{https://top500.org}~\cite{top500}.}
\label{fig:top500}
\end{figure}
Figure~\ref{fig:top500} shows the development of the performance of
the fastest 500 supercomputers in the world.
It shows that the available computation performance has increased by a
factor of one million over the last 17 years.
Until the early 2000s, the increase of performance was due to an
increase of the frequency of the processors.
Since then, the frequency stagnates at approximately
$\SI{2}{\giga\hertz}$.
Due to the higher power consumption and heat production at higher
frequencies, it is not efficient to increase the frequency further.
Hence, an increase of performance is only possible by an
increase of parallelism.
Another challenge in HPC is the power consumption of the over-all
system.
Modern supercomputers consume power in the scale of megawatts.
That is comparable to a whole offshore wind turbine.
Furthermore, the fault-tolerance of large computers is a problem as
well.
The more components are involved, the higher is the probability that
components fail during the computation.
Due to all these challenges, it is very important to develop software
that uses the hardware efficiently.

The problem of solving large sparse linear systems is a building block
in many HPC codes that consume large parts of the computation time.
As direct solvers scale badly for large linear systems and consume far too
much memory, iterative solvers are used on supercomputers to solve
this kind of problems.
Especially Krylov solvers have been approved to solve this problem.
For several reasons, Krylov solvers only utilize a fraction of the
peak-performance of supercomputers.
This is shown in the HPCG benchmark list~\cite{top500}.
For example, Fugaku, currently the fastest supercomputer in the world, only
performed $\SI{13.366}{\peta\flop}$ in the HPCG benchmark where it
reaches $\SI{415.53}{\peta\flop}$ in the LINPACK benchmark.
This shows the potential for improvements.

In this thesis, we consider three aspects of this issue.
First, we consider the increasing parallelism of larger
machines.
This parallelism appears on three levels:
\begin{enumerate}
\item Instruction level: The instruction sets of modern CPU contain
  instructions that perform multiple floating-point operations.
  For example, Fused-Multiply-Add (FMA) instructions, where a
  multiplication is carried out together with an addition.
  Other examples are Single-Instruction-Multiple-Data (SIMD)
  instructions, where the same operation is applied on multiple
  data.
\item Shared memory level: Modern CPUs consist of multiple cores that
  work in parallel but operate on the same memory.
\item Distributed memory level: Supercomputers are build from multiple nodes
that communicate over a network.
  Modern supercomputers have hundreds to many hundred thousand nodes.
  This number is expected to grow even further in the future.
\end{enumerate}
All this parallelism must be exploited to use the supercomputer
efficiently.

Another aspect is the so called memory-wall.
The bandwidth between the memory and the CPU is limited, which hinders the CPU
to exploit its full performance.
This effect is often mitigated by a hierarchical cache.
However, this does only work if the loaded data is reused enough.
A quantity to measure the reuse of the data is the arithmetic intensity (flop
per byte, $\si{\flop\per\byte}$), which is a property of the used algorithm.

The last aspect is a consequence of the distributed memory
parallelism.
The communication costs grow if more nodes are involved in the
computation.
A typical communication pattern that is used in Krylov solvers is a
collective communication, e.g.\ a global sum.
This type of communication scales as $\mathcal{O}\left(\log(P)\right)$,
where $P$ is the number of processors.
This makes it essential to organize the communication well and
overlap the communication phase with other meaningful computations.

In the literature, the data transfer between the different cache levels
as well as the data movement between nodes are referred to as communication.
In the present thesis, we want to strictly separate between the data movement
between cache levels which could be seen as intra-node communication
and the data movement between nodes which could be seen as
inter-node communication.

We consider large sparse linear systems that need to be
solved for multiple right-hand sides.
This is a very common problem that appears in applications like
inverse problems or optimization.
We will see that this type of problem is quite well posed to solve or
mitigate all the mentioned issues.

\section{Related Work}

The first part of this thesis is strongly inspired by the work of
\citeauthor[]{frommer2019block}~\cite[]{frommer2017block,frommer2019block}
and the PhD thesis by \citeauthor[]{lund2018block}~\cite[]{lund2018block}.
They recently presented the block Krylov framework on which this thesis is
built on.

The second part of this thesis is related to the work of
\citeauthor[]{cools2017numerically}~\cite[]{cools2017numerically,cools2018analyzing,cools2019numerically}.
They presented pipelined Krylov methods that overlap the
collective communication of the inner products with computation.
In the field of communication-avoiding methods, \citeauthor[]{demmel2008communication}~\cite[]{demmel2008communication,demmel2012communication}
as well as
\citeauthor{carson2015communication}~\cite[]{carson2015communication}
and
\citeauthor[]{hoemmen2010communication}~\cite[]{hoemmen2010communication}
presented several methods and ideas to avoid communication in Krylov
methods.
These methods fuse the communication of multiple iterations into
one communication, to reduce the number of messages.
Therefore, they are known as $s$-step Krylov methods.
Also combinations of $s$-step Krylov methods and block Krylov methods
have been proposed~\cite[]{chronopoulos2010block}.

Another approach is to use multiple search directions in a Krylov
method.
This could be found for example in the multi preconditioning methods
of \citeauthor[]{spillane2016adaptive}~\cite[]{spillane2016adaptive}.
Another class of algorithms, that fall into this category, are the
enlarged Krylov methods.
The idea is to transfer the advantageous convergence properties from
the block Krylov methods for multiple right-hand sides to systems with
a single right-hand side.
They were presented by \citeauthor[]{grigori2017reducing}~\cite[]{grigori2017reducing}.

In practice, several linear algebra software frameworks provide
optimized Krylov methods.
For example, \textsc{PETSc}~\cite{petsc-efficient,petsc-web-page} provides
communication-avoiding and pipelined Krylov methods, but lacks block Krylov
methods.
\textsc{Trilinos}~\cite[]{trilinos-website} contains a linear algebra module
that contains block Krylov methods and corresponding communication-avoiding
methods.
Other packages focus more on scalable preconditioners,
e.g.\ \textsc{hypre}~\cite[]{falgout2002hypre}.

Software packages that focus on the solution of PDEs often only
provide textbook Krylov methods, that are not explicitly optimized for
high-performance computing.
For example, \dune~\cite{bastian2020dune,dune24:16,dunepaperI:08,dunepaperII:08},
\textsc{deal.II}~\cite{dealII92} and \textsc{NGSolve}~\cite{schoberl1997netgen}.


\section{Contributions and Outline}

As already mentioned, we distinguish intra- and inter-node communication.
We pick up this difference to structure this thesis into two parts.
Part~\ref{part:blockkrylov} refers to the first two aspects mentioned in the
motivation, i.e.\ the vectorization and memory-wall.
In the second part, we optimize the methods further for inter-node communication,
which refers to the last aspect in the description above.

In Chapter~\ref{chap:blockkrylovframework}, we review the block Krylov
framework by \citeauthor[]{frommer2017block} and analyze its building blocks
with respect to their performance on modern CPU architectures.
We provide a novel view onto the set of possible *-subalgebras, based on three
elementary cases and introduce a new class of *-subalgebras.
Furthermore, based on the performance analysis we provide a guideline for
choosing an appropriate *-subalgebra.

Based on this framework, we formulate block versions of the CG, GMRes
and BiCGStab method in the Chapters~\ref{chap:blockcg},
\ref{chap:blockgmres} and \ref{chap:blockbicgstab}, respectively.
For the CG and BiCGStab methods, we introduce a novel stabilization
strategy that replaces the deflation process used in most
methods in the literature.
The new strategy is better suited in our context as we depend on a fixed
number of columns in the block vectors.

In the second part, we optimize the methods with respect to inter-node
communication.
For that, we adopt the approaches by \citeauthor[]{cools2017numerically} for our
block CG method in Chapter~\ref{chap:communicationcg}.
This yields a novel pipelined block Krylov method that combines the advantages
of both approaches.

In Chapter~\ref{chap:communicationgmres}, we consider the
orthogonalization procedure of the block GMRes method.
We introduce a pipelined Gram-Schmidt orthogonalization and an
innovative reduction-based orthogonalization and compare it with the
classical Gram-Schmidt method, which is the standard in up-to-date
methods.
The new methods prove to perform better and are more stable
than the classical Gram-Schmidt method.

All newly introduced methods are validated with numerical experiments,
carried out on a modern Intel compute server or on the
supercomputer PALMAII of the University of Münster.

%
\chapter{A Brief Introduction to Krylov Methods}
\label{chap:KrylovIntro}

Krylov methods came up in the 1950s.
\citeauthor[]{lanczos1950iteration}~\cite[]{lanczos1950iteration}
presented his method for solving eigenvalue problems in
\citeyear[][]{lanczos1950iteration}.
At the same time,
\citeauthor[]{hestenes1952methods}~\cite[]{hestenes1952methods}
presented the Conjugate Gradient (CG) method for solving linear
systems.
Back then, the CG method was considered a direct method.
Later, around 1975, with the development of vector computers and
massive memory computers the methods became more popular as iterative
methods.
The term Krylov method goes back to the Russian mathematician
\citeauthor[]{krylov1931numerical}, who presented related work in
\citeyear[][]{krylov1931numerical}~\cite{krylov1931numerical}.
Nowadays, lots of Krylov methods were developed and became an
essential part of modern scientific computing.
\citeauthor[]{golub1989some}~\cite[]{golub1989some} gave a good
overview over the early developments of Krylov methods.
Recommendable books about Krylov methods are written by
\citeauthor[]{greenbaum1997iterative}~\cite{greenbaum1997iterative},
\citeauthor[]{saad2003iterative}~\cite{saad2003iterative},
\citeauthor[]{hackbusch1994iterative}~\cite[]{hackbusch1994iterative}
and
\citeauthor[]{trefethen1997numerical}~\cite[]{trefethen1997numerical}.

We start with some basic definitions.
For the rest of this chapter, we consider a linear system
\begin{align}
  \label{eq:lgs}
  A\solution =b,
\end{align}
where $A \in \mathcal{L}(\bbR^{\dimA}, \bbR^{\dimA})$ is an invertible linear
operator, $b \in \bbR^{\dimA}$ is a given right-hand side and
$\solution\in\bbR^{\dimA}$ is the desired solution.

\begin{defn}[Krylov space]
  For $k\in\mathbb{N}$ and $r\in\bbR^{\dimA}$, the vector space
  \begin{align}
    \krylov{A}{r}{k} = \vspan[r, Ar, \ldots, A^{k-1}r]
  \end{align}
  is called the order-$k$ Krylov space generated by $A$ and $r$.
  The quantity
  \begin{align}
    \nu(r,A) = \max_{k\in\mathbb{N}} \left(\dim
    \krylov{A}{r}{k}\right)
  \end{align}
  is called the grade of $r$ with respect to $A$.
\end{defn}

The following lemma summarizes the most important properties of the Krylov space.
\begin{lem}[Properties of the Krylov space]
  \label{lem:krylov_properties}
  The following properties of the Krylov space hold
  \begin{itemize}
  \item $\krylov{A}{r}{k} \subseteq \krylov{A}{r}{k+1}$
  \item $\nu(A,r) \leq 1 + \operatorname{rank}\left(A\right)$
    and $\nu(A,r) \leq \dimA$.
  \item The vector space of polynomials $\bbP^{k-1}$ can be
    embedded into the Krylov space $\krylov{A}{r}{k}$ with the embedding
    \begin{align}
      \label{eq:polynomial_isomorthism}
  \iota: \bbP^{k-1} &\to \krylov{A}{r}{k}\\
  p &\mapsto p(A)r.
\end{align}
If $k\leq \nu(A,r)$, then $\iota$ is an isomorphism.
\end{itemize}
\end{lem}

The objective of a Krylov method is to find an approximation
$x^{k}\in\krylov{A}{r^{0}}{k}$ to the solution $A^{-1}b$, where $r^{0} = b -
Ax^{0}$ is the initial residual for an initial guess $x^{0}\in\bbR^{\dimA}$.
For example, the CG method~\cite[]{hestenes1952methods} computes the best
approximation with respect to the energy error $\|\solution - x^{k}\|_{A}$, and
the GMRes method~\cite[]{saad1986gmres} computes the best approximation with
respect to the residual norm $\|Ax^{k} -b \|$.
To compute this approximation, it is often helpful to use an orthonormal basis of
the Krylov space.
This orthonormal basis can be computed with the Arnoldi process~\cite[]{arnoldi1951the}, that is based
on the Gram-Schmidt orthogonalization process.
It computes an orthonormal basis $V\in \bbR^{\dimA \times k}$ that satisfies the
so-called Arnoldi relation
\begin{align}
  \label{eq:arnoldi}
  AV^{k} = V^{k}H^{k} + h_{k+1,k}v^{k+1}\transpose{e_{k}},
\end{align} where $H^{k}\in\bbR^{k \times k}$ is a Hessenberg matrix,
$v^{k+1}\in \bbR^{\dimA}$ is the subsequent basis vector and $e_{k}$ is the
$k$th unity vector.
For example, the GMRes method uses this relation to minimize the euclidean norm
of the residual.
If $h_{k+1,k}$ is small, $H$ is a good approximation for the operator $A$
restricted on the Krylov space.
In the case where $A$ is symmetric, it follows from~\eqref{eq:arnoldi} that $H$
is a tridiagonal matrix.
This fact is used in the CG and MINRES~\cite[]{paige1975solution} methods, such
that the basis $V$ does not need to be stored explicitly.
Instead the approximation is updated during the iteration.
This property is called \textit{short recurrence}.

From the definition, it is clear that the solution of the system~\eqref{eq:lgs}
is contained in the Krylov space $\krylov{A}{r^{0}}{\nu(A,r^{0})}$.
Therefore, Krylov space methods that compute a best approximation in the Krylov
space, terminate after at least $\nu(A,r^{0})$ steps.
However, Krylov methods are usually used to compute a good approximation for the
solution that is achieved before $\nu(A,r)$ iterations are performed.
In general, it is not possible to provide an error estimation that ensures
convergence with fewer than $\nu(A,r)$ iterations, as the following example
shows.
\begin{ex}
  \label{ex:nonconverging_gmres}
Consider the following system
\begin{align}
  A &= \begin{pmatrix}
    0 & 1 & 0 & \cdots & 0\\
    \vdots&\ddots & \ddots&\ddots&\vdots\\
    0&\cdots&0&1&0\\
    0&\cdots&&0&1\\
    1 & 0 & \cdots && 0
  \end{pmatrix}
                       & b&=e_{1}=\begin{pmatrix}
                         1\\0\\\vdots\\0
                       \end{pmatrix}.
\end{align}
With initial guess $x^{0} = 0$, the initial residual is $r^{0}=e_{1}$.
For all $k < \dimA$, all vectors in the Krylov space would have the
last coefficient $0$, as the operator $A$ pushes the coefficients one
place further.
As the solution of the system is $e_{\dimA}$, the best approximation
in the Krylov space is $0$.
Only if $k \geq \dimA$ the error norm can be decreased.
\end{ex}

Therefore, to show any results about convergence rates additional assumptions
are necessary.
For example, there are results if the symmetric part
$\frac12(A+\transpose{A})$ is positive definite which can be found in the
excellent books of
\citeauthor[]{greenbaum1997iterative}~\cite[]{greenbaum1997iterative} or
\citeauthor[]{saad2003iterative}~\cite[]{saad2003iterative}.
For the CG method there exists the following famous estimation of the energy
error.
\begin{thm}[Convergence of CG method]
  \label{thm:cg_convergence}
  Let $A$ be symmetric positive definite and $e^{k}=\solution - x^{k}$
  the error of the $k$th CG iteration. Then the energy error of
  $e^{k}$ can be estimated by
  \begin{align}
    \|e^{k}\|_{A} \leq 2
    \left(\frac{\sqrt{\kappa}-1}{\sqrt{\kappa}+1}\right)^{k}\|e^{0}\|_{A},
  \end{align}
  where $\kappa = \|A\|\|A^{-1}\|$ is the condition number of $A$.
\end{thm}
The proof is geared to the one presented by \citeauthor[]{trefethen1997numerical}~\cite[]{trefethen1997numerical}.
\begin{proof}
  By Lemma~\ref{lem:krylov_properties}, we can identify every element in
  the Krylov space $\krylov{A}{r^{0}}{k}$ by a polynomial of degree
  $k-1$.
  In particular, we write the $k$th error of the CG method as
  \begin{equation}
    \label{eq:error_representation}
    \begin{aligned}
      e^{k} &= \solution - x^{k}\\
      &= \solution - x^{0} - p_{k-1}(A)r^{0}\\
      &= \solution - x^{0} - p_{k-1}(A)A(\solution - x^{0})\\
      &= q_{k}(A)e^{0},
    \end{aligned}
  \end{equation}
  for a polynomial $p_{k-1}\in\bbP^{k-1}$ and $q_{k}(\polyX) \colon= p_{k-1}(\polyX)\polyX + 1$.
  As the CG methods finds the best approximation with respect to the
  energy norm, we conclude
  \begin{align}
    \|e^{k}\|_{A} \leq \inf_{q_{k}} \|q_{k}(A)e^{0}\|_{A}.
  \end{align} Here the infimum is taken over all polynomials of degree $k$ with
absolute coefficient $1$.
  For the smallest and largest eigenvalues $\lambda_{\min}$ and
$\lambda_{\max}$, the polynomials that realize this infimum are given by the
scaled Chebyshev polynomials
  \begin{align}
    \tilde{T}_{k}(\polyX) =
    \left(T_{k}\left(\frac{-\lambda_{\max}-\lambda_{\min}}{\lambda_{\max}-\lambda_{\min}}\right)\right)^{-1}T_{k}\left(\frac{2\polyX-\lambda_{\max}-\lambda_{\min}}{\lambda_{\max}
    -\lambda_{\min}}\right),
  \end{align}
  where the $k$th Chebyshev polynomial $T_{k}$ is defined by the recursion
  formula
  \begin{align}
    \label{eq:chebyshev}
    T_{0}(\polyX) &= 1,\qquad T_{1}(\polyX) = \polyX\\
    T_{k+1}(\polyX)& = 2\polyX T_{k}(\polyX)-T_{k-1}(\polyX)
  \end{align}
  or directly by
  \begin{align}
    T_{k}(\polyX) = \cos\left(k\arccos\left(\polyX\right)\right).
  \end{align}
  One can show that the scaled Chebyshev polynomials minimize the
  $C^{\infty}$-norm on the interval
  $\left[\lambda_{\min},\lambda_{\max}\right]$ in the space of
  polynomials with absolute coefficient $1$.
  The $C^{\infty}$-norm is bounded by
  \begin{align}
    \label{eq:chebyshev_norm}
    \|\tilde{T}_{k}\|_{C^{\infty}} \leq 2\left(\frac{\sqrt{\kappa}-1}{\sqrt{\kappa}+1}\right)^{k},
  \end{align}
  where $\kappa = \norm{A}\norm{\inverse{A}} =
  \frac{\lambda_{\max}}{\lambda_{\min}}$ denotes the condition number of the
  operator $A$.
  As $A$ is symmetric positive definite the eigenvectors
  $\left(u_{i}\right)_{i}$ to the eigenvalues
  $\left(\lambda_{i}\right)_{i}$ build an orthonormal basis of
  $\bbR^{\dimA}$.
  We write the error $e^{0}$ in this basis as
  \begin{align}
    e^{0} = \sum_{i=0}^{\dimA -1} a_{i}u_{i}.
  \end{align}
  Then the energy error of $e^{0}$ is given by
  \begin{align}
    \|e^{0}\|_{A}^{2} = \sum_{i=0}^{\dimA-1}a_{i}^{2}\lambda_{i}
  \end{align}
  and the energy error of $e^{k}$ is given by
  \begin{align}
    \|e^{k}\|_{A}^{2} &= \|q_{k}(A)e^{0}\|_{A}^{2}\\
                      &= \sum_{i=0}^{\dimA-1}
                        q_{k}(\lambda_{i})^{2}a_{i}^{2}\lambda_{i}\\
                      &\leq \max_{i=0}^{\dimA-1}|q_{k}(\lambda_{i})|^{2}\|e^{0}\|_{A}^{2}\\
                      &\leq
                        4\left(\frac{\sqrt{\kappa}-1}
                        {\sqrt{\kappa}+1}\right)^{2k}\|e^{0}\|_{A}^{2},
  \end{align}
  where we used Equation~\eqref{eq:error_representation} and~\eqref{eq:chebyshev_norm}.
\end{proof}

The error bound given by Theorem~\ref{thm:cg_convergence} is sharp,
i.e.\ there exists data $A$ and $b$ such that equality holds.
But for fixed data better convergence could occur.
For example, if the initial residual $r^{0}$ is an eigenvector of $A$,
then the method would converge within one iteration, as the grade of
$r^{0}$ with respect to $A$ is $1$.

For the GMRes method a representation for the residual similar to
Equation~\eqref{eq:error_representation} can be formulated as
\begin{align}
  r^{k} &= b - Ax^{k}\\
        &= b - Ax^{0} - Ap_{k-1}(A)r^{0}\\
        &= r^{0} - Ap_{k-1}(A)r^{0}\\
        &= q_{k}(A)r^{0},
\end{align} with $q_{k}(\polyX) = 1 - p_{k-1}(\polyX)\polyX$.
From this equality, an error estimation could be derived, if the operator $A$ is
normal, i.e.\ diagonalizable.
We review this prove for the block variant of the GMRes method in Section~\ref{sec:gmres_convergence}.

Theorem~\ref{thm:cg_convergence} and the theory of the proof show that
the convergence behavior of Krylov space methods depend on the
condition number of the operator.
Therefore, it is common practice to use preconditioning.
That means the Krylov method is applied on the system
\begin{align}
  \inverse{M_{L}}A\inverse{M_{R}}y = \inverse{M_{L}}b,
\end{align} for some matrices $M_{L}, M_{R}$ for which the inverse can be
applied cheaply.
Once $y$ has been found, the solution of $Ax=b$ can be found easily by computing
$x = \inverse{M_{R}}y$.
The operators $M_{L}, M_{R}$ are chosen to improve the condition number of the
operator $\inverse{M}_{L}A\inverse{M_{R}}$ and hence to improve the convergence
of the Krylov method.
Often one of $M_{L}$ and $M_{R}$ is chosen to be the identity, resulting in
so-called left or right preconditioning.

Simple preconditioners depend on iterative splitting methods like
Jacobi or Gauß-Seidel iteration.
They split the operator into a sum of matrices
\begin{align}
  A = M+N,
\end{align}
where $N$ can be easily inverted.
For example, the Jacobi method chooses $N$ as the diagonal of $A$ and
the Gauss-Seidel method chooses $N$ as the lower triangular part of $A$.
Then the preconditioner is given by some iterations of the fixpoint iteration
\begin{align}
  x^{k+1} &= \inverse{N}\left(b-Mx^{k}\right)\\
          &= x^{k} + \inverse{N}r^{k}.
\end{align}
Other popular preconditioners compute incomplete factorization of the
operator $A$.
These preconditioners often only affect the large eigenvalues
of $A$.
Especially on very large systems this does not reduce the
condition number sufficiently, as the small eigenvalues are not affected.
More sophisticated preconditioners are multi-grid methods that use
restrictions of the operator $A$ to coarser spaces and apply the
simple preconditioners on that level too.
Thus, all ranges of eigenvalues are affected.
An alternative is to compute coarse spaces that contain
the eigenvectors of the small eigenvalues.
The preconditioner is then chosen as the projection onto the orthogonal
complement of this coarse spaces
(e.g.\ GenEO~\cite[]{spillane2014abstract}).

\begin{table}
  \renewcommand{\arraystretch}{2.8}
  \small
  \centering
  \caption{Overview of commonly used Krylov methods their properties
    and references.}
  \begin{tabular}{lcccc}\toprule
    Name & Requirements &  \makecell{Short\\Recursion} & Minimization & Reference\\
    \midrule
    \makecell[l]{Conjugate\\Gradients (CG)} & \makecell{symmetric\\
    positive definite}. & yes & $\|e^{k}\|_{A}$ & \cite[]{hestenes1952methods}\\
    \makecell[l]{General Minimal\\
    Residual (GMRes)} & none & no & $\|r^{k}\|_{2}$ & \cite[]{saad1986gmres}\\
    \makecell[l]{Biconjugate Gradients\\
    Stabilized (BiCGStab)} & none & yes & none & \cite[]{van1992bi}\\
    \makecell[l]{Minimum Residual\\(MINRes)} & symmetric & yes & $\|r^{k}\|_{2}$ & \cite[]{paige1975solution}\\
    \makecell[l]{Conjugate\\Residual (CR)} & symmetric & yes& $\|r^{k}\|_{2}$ &\cite[]{stiefel1955relaxationsmethoden},\cite[]{eisenstat1983variational} \\
    \makecell[l]{Quasi Minimal\\
    Residual (QMR)} & none & yes & none & \cite[]{freund1991qmr}\\
             \bottomrule
  \end{tabular}
  \label{tab:krylovmethodsoverview}
\end{table}
Table~\ref{tab:krylovmethodsoverview} shows an overview of widely used
Krylov methods for solving linear systems.
It shows the requirements for the operator and preconditioner as well
as whether it uses a short recursion.
Furthermore, the norm in which the error is minimized is given and
the citation in which the method was presented.

%
\cleardoublepage
\epigraphhead[500]{\epigraph{\itshape La théorie est mère de la pratique.}
{\textsc{Louis Pasteur}}}
\part{Block Krylov Methods}
\label{part:blockkrylov}

\chapter{A General Block Krylov Framework}
\label{chap:blockkrylovframework}
Block Krylov methods have been developed in the 1970s
 and 1980s to solve linear
systems with multiple right-hand sides \cite{oleary1980block}
or compute multiple eigenvectors \cite{underwood1975iterative}.
Recently, they have been rediscovered in the context of
high-performance computing to reduce the communication overhead.

The term ``block'' is quite overloaded in the field of numerical linear algebra.
In the context of matrix structures it means that the matrix is subdivided into
smaller matrices.
In the context of preconditioning it often refers to the block Jacobi method
that only considers the diagonal blocks of the system matrix to parallelize the
preconditoning, and in the context of Krylov methods it refers to the already
mentioned methods that are based on the work of \citeauthor[]{oleary1980block}~\cite[]{oleary1980block}.

We consider block Krylov methods to solve linear systems
with multiple right-hand sides.
Let $A\in \bbR^{\dimA\times\dimA}$ be an invertible linear operator and $B \in
\bbR^{\dimA\times s}$ a block vector.
A linear system with multiple right-hand sides, called a \textit{block system},
for the solution $\blocksolution\in\bbR^{\dimA\times s}$ is given by
\begin{align}
  \label{eq:blocklinearsystem}
  A\blocksolution &= B,
\end{align}
which is equivalent to
\begin{align}
  \label{eq:blocklinearsystem2}
  A\solution_i &= b_i& \forall i=1,\ldots,s,
\end{align}
where $\solution_{i}$ and $b_{i}$ denote the $i$th column of $\blocksolution$ and $B$,
respectively.

The basic idea of block Krylov methods is to make use of the sum of
all Krylov spaces of the linear systems in Equation \eqref{eq:blocklinearsystem2} to find a better
approximation for the solution.
\citeauthor[]{oleary1980block}~\cite[]{oleary1980block} showed that the
convergence of the block CG method is faster than that of the CG method and
independent of the $s-1$ smallest eigenvalues.
We recall this result in Theorem \ref{thm:blockcg}.

In the context of high-performance computing block Krylov methods have
another advantage.
During one iteration the operator (and preconditioner) is applied to
block vectors in $\bbR^{\dimA\times s}$, which is beneficial if the
matrix is explicitly stored.
It leads to a higher arithmetical intensity, which is crucial on
modern CPUs to achieve good performance.
Furthermore, it is well suited for the use of SIMD instructions if the
block vectors are stored in row-major format.

Several approaches have been proposed to use the faster convergence of
block Krylov methods for linear systems with a single right-hand side.
\citeauthor[]{grigori2017reducing}~\cite[]{grigori2017reducing,al2018enlarged}
proposed a method where they decompose the right-hand side $b$, based on the
domain decomposition, to obtain multiple right-hand sides which can be
used to solve the original problem.
This approach is also used in the PhD theses of
\citeauthor[]{moufawad2014enlarged}~\cite{moufawad2014enlarged}, \citeauthor[]{aldaas2018thesis}~\cite{aldaas2018thesis}
 and \citeauthor[]{tissot2019iterative}~\cite[]{tissot2019iterative}.

Other approaches are to choose additional right-hand sides randomly
(BRRHS-CG)~\cite[]{nikishin1995variable} or to choose additional initial guesses
randomly and solve all for the same right-hand side
(CoopCG)~\cite[]{bhaya2012coop}.
In principle, these approaches are also applicable to the methods
presented in this work.

One iteration of a block Krylov method has costs in order of
$\mathcal{O}(s^2\dimA+s^3)$, which could become a problem, if a lot of right-hand
sides are used, i.e.\ $s$ is large.
To mitigate this effect, but still take advantage of the higher operational
intensity of the operator and preconditioner application, we introduce a general framework of block Krylov spaces based on the
work of
\citeauthor[]{frommer2017block} \cite[]{frommer2017block,frommer2019block,lund2018block}.
That allows us to balance the information exchange between the different
right-hand sides and the computational blocking overhead.
Then, we provide a performance analysis of the building blocks, provide details
about our implementation and present some numerical tests that approve our
theory and show the advantages of the block Krylov framework.

\section{Block Krylov Spaces}
Let us start with the review of the block Krylov framework presented by
\citeauthor{frommer2017block}~\cite{frommer2017block}. Originally this
framework was introduced to evaluate functions of matrices.
Further development of the framework was done in the
thesis by \citeauthor{lund2018block}~\cite{lund2018block} and the paper by
\citeauthor{frommer2019block}~\cite{frommer2019block}.

In the subsequent of this work all methods and algorithms are built upon this
framework.
The standard Krylov methods can be obtained by choosing $s=1$, this is referred
to as the non-block case.
We start with the central definition of the block Krylov space.

\pagebreak
\begin{defn}[Block Krylov subspace]
  \label{defn:blockkrylovspace}
  Let $\SubA$ be a *-subalgebra of $\bbR^{s\times s}$ and $R\in\bbR^{\dimA\times s}$.
  The $k$th block Krylov space with respect to $A, R$ and $\SubA$ is defined by
  \begin{align}
    \krylov[\SubA]{A}{R}{k} = \left\{ \sum_{i=0}^{k-1} A^iRc_i \,\big|\,
    c_0,\ldots,c_{k-1}\in\SubA \right\} \subset \bbR^{\dimA \times s}.
  \end{align}
\end{defn}
\begin{rem}\hfill
  \vspace{-\topsep}
  \begin{itemize}
  \item A *-algebra is a vector space $\SubA$ equipped with a
    product and a conjunction.  In particular, it means, that for all
    elements $s\in\SubA$ and polynomials $p\in\bbP$ the evaluation of
    the polynomial for that element $p(s)$ is contained in the
    *-algebra.
  \item The Cayley-Hamilton theorem yields that every *-subalgebra
    of $\bbR^{s\times s}$ contains an identity.
    See for example \cite[]{bosch2014lineare}.
    This identity does not necessarily coincide with the identity in
    $\bbR^{s\times s}$.
    However, for the *-subalgebras $\SubA$ we consider in this work,
    the identities coincide $\Identity[\SubA] =
    \Identity[\bbR^{s\times s}]$. Other *-subalgebras would be
    pointless, as we will see later.
  \item We choose $\bbR^{\dimA\times s}$ as a vector space here.
    In principle every vector space over some field $F$ could be
    chosen.
    $\SubA$ is then a *-subalgebra of $F^{s\times s}$.
  \item For the rest of this thesis $\SubA$ denotes a *-subalgebra of
    $\bbR^{s\times s}$.
    \item The classical block Krylov methods as described by
      \citeauthor[]{oleary1980block} use $\SubA = \bbR^{s\times s}$.
  \end{itemize}
\end{rem}

For the convergence theory of Krylov methods, polynomials play an
important role, as we already saw in Theorem \ref{thm:cg_convergence}.
For the convergence theory in this framework, we introduce the more
generic $\SubA$-valued polynomials.
\begin{defn}
  A polynomial of the form
  \begin{align}
    \mathcal{P}(\polyX) &= \sum_{i=0}^{k}\polyX^{i}\gamma_{i}& \gamma_{i}\in\SubA
  \end{align}
  is called a $\SubA$-valued polynomial of degree $k$.
  We write $\bb{P}_{\SubA}^{k}$ for the space of $\SubA$-valued
  polynomials of degree $k$.
Inspired by the paper of \citeauthor[]{el2003block}~\cite[]{el2003block} we denote the product
\begin{align}
  \mathcal{P}(A)\circ Y = \sum_{i=0}^{k}A^{i}Y\gamma_{i},
\end{align}
where $Y\in\bbR^{\dimA \times s}$.
With this operation, the operator $\mathcal{P}(A)$ could be considered
as a linear operator on the space $\bbR^{\dimA \times s}$.
Furthermore, we define the right-sided product $\mathcal{P}\sigma \in \bb{P}_{\SubA}^{k}$ of
a $\SubA$-valued polynomial $\mathcal{P}\in\mathbb{P}_{\SubA}^{k}$, with
$\mathcal{P}(\polyX) = \sum_{i=0}^{k}\polyX^{i}\gamma_{i}$,
and $\sigma\in\SubA$ as
\begin{align}
  \left(\mathcal{P}\sigma\right)(\polyX) = \sum_{i=0}^{k}\polyX^{i}\gamma_{i}\sigma.
\end{align}
\end{defn}

From Definition \ref{defn:blockkrylovspace} we find the following two lemmas
immediately. The first one is in analogy with \eqref{eq:polynomial_isomorthism}.
\begin{lem}
  \label{lem:blockkrylov_polynomial}
  Every element $X\in\krylov[\SubA]{A}{R}{k}$ in the block
  Krylov space can be represented by a
  $\SubA$-valued
  polynomial
  $\mathcal{P}\in\bb{P}^{k-1}_{\SubA}$ of degree $k-1$
  \begin{align} X = \mathcal{P}(A) \circ R &= \sum_{i=0}^{k-1} A^{i}Rc_{i}&
c_{0},\ldots,c_{k-1}\in\SubA.
  \end{align}
\end{lem}

\begin{lem}
  \label{lem:krylovembedding}
  If $\SubA[1]$ and $\SubA[2]$ are two *-subalgebras of $\bbR^{s
    \times s}$,
  with $\SubA[1] \subseteq \SubA[2]$.
  Then
  \begin{align}
    \krylov[{\SubA[1]}]{A}{R}{k} \subseteq \krylov[{\SubA[2]}]{A}{R}{k}
  \end{align}
  holds.
\end{lem}

Analogous to the non-block case ($s=1$) we define the block grade of a
block vector in a block Krylov space.
This definition is inspired by \citeauthor[]{gutknecht2009block}~\cite[]{gutknecht2009block,gutknecht2007block}.
\begin{defn}[Block grade]
  In the setting of Definition \ref{defn:blockkrylovspace} we define the
  block grade of $R$ with respect of $A$ as
  \begin{align}
    \nu_{\SubA}(A,R) := \min \left\{k\in\bb{N}\,\big| \dim
    \krylov[\SubA]{A}{R}{k} = \dim \krylov[\SubA]{A}{R}{k+1}\right\}.
  \end{align}
\end{defn}

As \citeauthor[]{gutknecht2009block} show, the block grade defines the
minimal $k$ for which the solution is contained in the Krylov space.
We review this result in our context.
\begin{lem}
  We have
  \begin{align}
    \blocksolution \in X^{0} +
    \krylov[\SubA]{A}{R^{0}}{\nu_{\SubA}(A,R^{0})},
  \end{align}
  where $A\blocksolution = B$ and
  $R^{0} = B-AX^{0}$ is the residual for some initial guess
  $X^{0}\in\bbR^{\dimA\times s}$.
\end{lem}
\begin{proof}
  By definition of the block Krylov space we have for any $k\in\bb{N}$
  \begin{align}
    \krylov[\SubA]{A}{R^{0}}{k} \subseteq \krylov[\SubA]{A}{R^{0}}{k+1}.
  \end{align}
  From the definition of the block grade it follows that
  \begin{align}
    \krylov[\SubA]{A}{R^{0}}{\nu_{\SubA}(A,R^{0})} = \krylov[\SubA]{A}{R^{0}}{k}
  \end{align}
  for all $k \geq \nu_{\SubA}(A,R^{0})$.
  As $A$ is invertible, $\krylov[\SubA]{A}{R^{0}}{\nu_{\SubA}(A,R^{0})}$ is a
  $A$-invariant subspace of $\bbR^{\dimA\times s}$
  \begin{align}
    A\krylov[\SubA]{A}{R^{0}}{\nu_{\SubA}(A,R^{0})} =
    \krylov[\SubA]{A}{R^{0}}{\nu_{\SubA}(A,R^{0})}.
  \end{align}
  By definition $\krylov[\SubA]{A}{R^{0}}{\nu_{\SubA}(A,R^{0})}$
  contains $R^{0}$. This yields
  \begin{align}
    R^{0} \in A\krylov[\SubA]{A}{R^{0}}{\nu_{\SubA}(A,R^{0})}.
  \end{align}
  As $A$ is invertible, we can apply $\inverse{A}$ to get
  \begin{align}
    \blocksolution-X^{0} \in \krylov[\SubA]{A}{R^{0}}{\nu_{\SubA}(A,R^{0})}.
  \end{align}
  Adding $X^{0}$ completes the proof.
\end{proof}

This result is more of theoretical interest as the block grade in real
world problems is usually quite high.
In practice, often much fewer iterations are needed to reduce the
residual norm sufficiently.
As we will see later, a more practical relevant quantity is
defined by
\begin{align}
  \label{eq:blockkrylov_xi}
  \xi_{\SubA}(A,R^{0}) = \min \left\{k\in\bb{N}\,\big| \dim{\krylov[\SubA]{A}{R^{0}}{k} <
  k\dim \SubA}\right\}.
\end{align}
It is the iteration number in which the Krylov space does not grow by $k\dim
\SubA$ dimensions in every iteration.
This leads to a situation that must be treated numerically.

Next, we define a inner product on the vector space of block vectors
$\bbR^{n\times s}$ that is a generalization of the scalar product.
\begin{defn}[block inner product]
  A mapping $\sproduct[\SubA]{\cdot}{\cdot}: \bbR^{\dimA \times s}\times
 \bbR^{\dimA \times s} \to \SubA$ is called a \textit{block inner
   product} if the following conditions hold for all
 $X,Y,Z\in\bbR^{\dimA \times s}$ and $\gamma\in\SubA$:
 \begin{itemize}
 \item $\SubA$-linearity:
   $\sproduct[\SubA]{X+Y}{Z\gamma} = \sproduct[\SubA]{X}{Z}\gamma +
   \sproduct[\SubA]{Y}{Z}\gamma$
 \item symmetry: $\sproduct[\SubA]{X}{Y} =
   \transpose{\sproduct[\SubA]{Y}{X}}$
 \item definiteness: $\sproduct[\SubA]{X}{X}$ is positive definite for all
   full rank $X$, and $\sproduct[\SubA]{X}{X}=0$, if and only if
   $X=0$.
 \item normality: $\trace{\sproduct[\SubA]{X}{Y}} = \langle X, Y \rangle_{F}$
 \end{itemize}
 Here $\langle X, Y \rangle_{F} = \trace{\transpose{X}Y}$ denotes the
 Frobenius scalar product.

\end{defn}

\begin{rem}
  Note that the definitness condition implies that a block inner
  product can only be defined on *-subalgebras that contain full rank
  matrices.
  In particular, this yields that the identity of $\SubA$ is the same
  as the identity in $\bbR^{s\times s}$.
\end{rem}

\begin{defn}[normalizer]
  We call a map
  \begin{align}
    \normalizer{\SubA}: \bbR^{\dimA\times s}\to \SubA
  \end{align}
  a normalizer or scaling quotient, if for all $X\in\bbR^{\dimA\times s}$
  there exists a $Y$ such that
  \begin{align}
    X &= Y\normalizer[X]{\SubA}& \text{and}&& \sproduct[\SubA]{Y}{Y} = \Identity.
  \end{align}
\end{defn}
\pagebreak[3]
\begin{rem}\leavevmode\vspace{-\topsep}
  \begin{itemize}
  \item A normalizer can be computed by a QR factorization.
  \item In the work of \citeauthor[]{frommer2019block}, the scaling quotient is
    only required to be defined for full rank $X$.
    We use the more restrictive definition for our
    stabilization strategies and to resolve breakdowns in the block Krylov
    methods.
  \item We use the Householder algorithm to compute a normalizer in
    our code.
    The Gram-Schmidt orthogonalization progress would fail for
    rank-deficient block vectors.
    However, there is work to mitigate this, for example by replacing linear
    dependent columns with random vectors \cite[]{soodhalter2015block}.
  \item In algorithms we write the normalizer in python style syntax
    \begin{align}
      Y,\sigma = \normalizer[X]{\SubA}.
    \end{align}
  \end{itemize}
\end{rem}

Next, we take a look at different choices for the *-subalgebra $\SubA$. We first
introduce three elementary cases, that where also considered by
\citeauthor{frommer2017block} \cite{frommer2017block}, while we use a different
naming scheme.
\begin{defn}[elementary *-subalgebras]
  \label{defn:elementary-algebras}
  We consider the following three elementary cases of how to treat a
  block system within the Krylov framework.
  \begin{enumerate}
  \item global: The block system is considered as one linear
    system.
    This linear system can be represented by the Kronecker system
    \begin{align} \left(\Identity[s] \otimes A\right)\operatorname{vec}(X) &=
                                                                              \operatorname{vec}(B)\\
      \intertext{or}
      \begin{pmatrix}
        A\\
        &\ddots\\
        &&A
      \end{pmatrix}
           \begin{pmatrix}
             X_{1}\\
             \vdots\\
             X_{s}
           \end{pmatrix}
      &=
      \begin{pmatrix}
        B_{1}\\
        \vdots\\
        B_{s}
      \end{pmatrix}.
    \end{align}
    This choice corresponds to the *-subalgebra of multiples of the
    identity $\SubA[G] = \bbR \cdot\Identity$ and the Frobenius inner product
    \begin{align}
      \sproduct[{\SubA[G]}]{X}{Y} =
      \left(\sum_{i=1}^{s} \transpose{X_{i}}Y_{i}\right) \Identity.
    \end{align}
    This method goes back to \citeauthor[]{jbilou1999global}
    \cite[]{jbilou1999global} for the FOM and GMRes method.
  \item parallel: All columns of the block system are considered
    separately,
    but iterations are carried out simultaneously.
    This corresponds to the *-subalgebra of diagonal matrices
    $\SubA[P] = \diag{\bbR^{s}}$ and the inner block product
    \begin{align}
      \sproduct[{\SubA[P]}]{X}{Y} = \diag{\transpose{X}Y}.
    \end{align}
  \item block: The classic block Krylov case as presented by
    \citeauthor[]{oleary1980block}.
    This corresponds to the *-subalgebra $\SubA[B] = \bbR^{s\times s}$ and the
    inner product
    \begin{align}
      \sproduct[{\SubA[B]}]{X}{Y} = \transpose{X}Y.
    \end{align}
  \end{enumerate}
\end{defn}

From these three cases we compose more complex cases as the following definition
shows.
\begin{defn}[Relevant *-subalgebras]
  \label{defn:subalgebras}
  Let $p\in \bb{N}$ be a divider of $s$, $q=\frac{s}{p}$ and $X,Y \in
  \bbR^{\dimA\times s}$.
  We subdivide $X,Y$ column-wise into $\bbR^{\dimA\times p}$ matrices
  \begin{align}
    X &= \left[X_{1},\ldots,X_{q}\right] & Y&= \left[Y_{1},\ldots,Y_{q}\right].
  \end{align}
  Then we define the following *-subalgebras and corresponding block
  inner products:
  \begin{align}
    \text{block-parallel: }
    && \SubA[BP]^p &:= \diag{\left(\bbR^{p\times p}\right)^q}
    & \sproduct[{\SubA[BP]^p}]{X}{Y} &:= \diag{\transpose{X_1}Y_1,\ldots,\transpose{X_q}Y_q},\\
    \text{block-global: }
    && \SubA[BG]^p &:= \Identity[q] \otimes \bbR^{p\times p}
    & \sproduct[{\SubA[BG]^p}]{X}{Y} &:=\frac{1}{q} \sum_{i=1}^q
                                       \Identity[q] \otimes \transpose{X_i}Y_i,
  \end{align} where $\Identity[q]$ denotes the $q$ dimensional identity matrix and
$\diag{\left(\bbR^{p\times p}\right)^q}$ denotes the set of $s\times s$ matrices
where only the $p\times p$ diagonal matrices have non-zero values.
\end{defn}

In principle also a global-parallel combination would be possible.
As we want to make use of the advantages of the block strategy, we do not
consider this in the present thesis.
It would also be possible to apply the elementary case in a different order and
construct parallel-block or global-block methods.
The resulting *-algebras would be isomorphic to the ones of the block-parallel
and block-global methods.
Therefore, we restrict ourselves to the two mentioned cases.

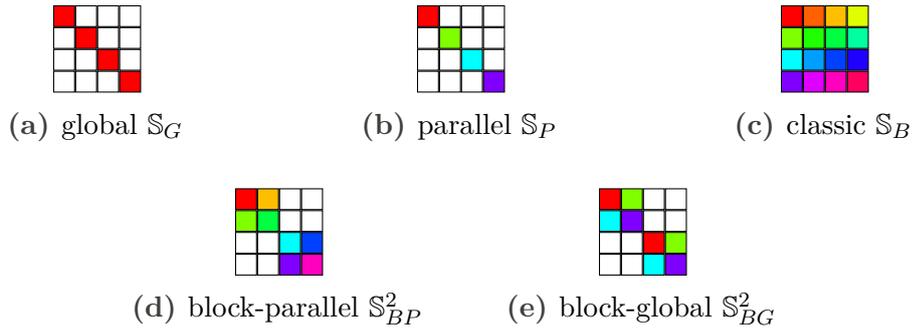
\begin{figure}
  \centering
  \begin{subfigure}[t]{0.3\textwidth}
    \centering
  \begin{tikzpicture}[Node/.style = {rectangle, draw,
        anchor=center,node distance=0mm}]
      \node[Node=5ex,fill=tikzcolor1] (1) {};
      \node[Node=5ex,fill=white,right=of 1] (2) {};
      \node[Node=5ex,fill=white,below=of 1] (3) {};
      \node[Node=5ex,fill=tikzcolor1,below=of 2] {};
      \node[Node=5ex,fill=white,right=of 2] (1b) {};
      \node[Node=5ex,fill=white,right=of 1b] (2b) {};
      \node[Node=5ex,fill=white,below=of 1b] {};
      \node[Node=5ex,fill=white,below=of 2b] {};
      \node[Node=5ex,fill=white, below=of 3] (1c) {};
      \node[Node=5ex,fill=white,right=of 1c] (2c) {};
      \node[Node=5ex,fill=white,below=of 1c] {};
      \node[Node=5ex,fill=white,below=of 2c] {};
      \node[Node=5ex,fill=tikzcolor1,right=of 2c] (1bc) {};
      \node[Node=5ex,fill=white,right=of 1bc] (2bc) {};
      \node[Node=5ex,fill=white,below=of 1bc] {};
      \node[Node=5ex,fill=tikzcolor1,below=of 2bc] {};
    \end{tikzpicture}
    \caption{global $\SubA[G]$}
  \end{subfigure}
  \begin{subfigure}[t]{0.3\textwidth}
    \centering
  \begin{tikzpicture}[Node/.style = {rectangle, draw,
        anchor=center,node distance=0mm}]
      \node[Node=5ex,fill=tikzcolor1] (1) {};
      \node[Node=5ex,fill=white,right=of 1] (2) {};
      \node[Node=5ex,fill=white,below=of 1] (3) {};
      \node[Node=5ex,fill=tikzcolor5,below=of 2] {};
      \node[Node=5ex,fill=white,right=of 2] (1b) {};
      \node[Node=5ex,fill=white,right=of 1b] (2b) {};
      \node[Node=5ex,fill=white,below=of 1b] {};
      \node[Node=5ex,fill=white,below=of 2b] {};
      \node[Node=5ex,fill=white, below=of 3] (1c) {};
      \node[Node=5ex,fill=white,right=of 1c] (2c) {};
      \node[Node=5ex,fill=white,below=of 1c] {};
      \node[Node=5ex,fill=white,below=of 2c] {};
      \node[Node=5ex,fill=tikzcolor9,right=of 2c] (1bc) {};
      \node[Node=5ex,fill=white,right=of 1bc] (2bc) {};
      \node[Node=5ex,fill=white,below=of 1bc] {};
      \node[Node=5ex,fill=tikzcolor13,below=of 2bc] {};
    \end{tikzpicture}
    \caption{parallel $\SubA[P]$}
  \end{subfigure}
  \begin{subfigure}[t]{0.3\textwidth}
    \centering
  \begin{tikzpicture}[Node/.style = {rectangle, draw,
        anchor=center,node distance=0mm}]
      \node[Node=5ex,fill=tikzcolor1] (1) {};
      \node[Node=5ex,fill=tikzcolor2,right=of 1] (2) {};
      \node[Node=5ex,fill=tikzcolor5,below=of 1] (3) {};
      \node[Node=5ex,fill=tikzcolor6,below=of 2] {};
      \node[Node=5ex,fill=tikzcolor3,right=of 2] (1b) {};
      \node[Node=5ex,fill=tikzcolor4,right=of 1b] (2b) {};
      \node[Node=5ex,fill=tikzcolor7,below=of 1b] {};
      \node[Node=5ex,fill=tikzcolor8,below=of 2b] {};
      \node[Node=5ex,fill=tikzcolor9, below=of 3] (1c) {};
      \node[Node=5ex,fill=tikzcolor10,right=of 1c] (2c) {};
      \node[Node=5ex,fill=tikzcolor13,below=of 1c] {};
      \node[Node=5ex,fill=tikzcolor14,below=of 2c] {};
      \node[Node=5ex,fill=tikzcolor11,right=of 2c] (1bc) {};
      \node[Node=5ex,fill=tikzcolor12,right=of 1bc] (2bc) {};
      \node[Node=5ex,fill=tikzcolor15,below=of 1bc] {};
      \node[Node=5ex,fill=tikzcolor16,below=of 2bc] {};
    \end{tikzpicture}
    \caption{classic $\SubA[B]$}
  \end{subfigure}\\[3ex]
    \begin{subfigure}[t]{0.3\textwidth}
    \centering
  \begin{tikzpicture}[Node/.style = {rectangle, draw,
        anchor=center,node distance=0mm}]
      \node[Node=5ex,fill=tikzcolor1] (1) {};
      \node[Node=5ex,fill=tikzcolor3,right=of 1] (2) {};
      \node[Node=5ex,fill=tikzcolor5,below=of 1] (3) {};
      \node[Node=5ex,fill=tikzcolor7,below=of 2] {};
      \node[Node=5ex,fill=white,right=of 2] (1b) {};
      \node[Node=5ex,fill=white,right=of 1b] (2b) {};
      \node[Node=5ex,fill=white,below=of 1b] {};
      \node[Node=5ex,fill=white,below=of 2b] {};
      \node[Node=5ex,fill=white, below=of 3] (1c) {};
      \node[Node=5ex,fill=white,right=of 1c] (2c) {};
      \node[Node=5ex,fill=white,below=of 1c] {};
      \node[Node=5ex,fill=white,below=of 2c] {};
      \node[Node=5ex,fill=tikzcolor9,right=of 2c] (1bc) {};
      \node[Node=5ex,fill=tikzcolor11,right=of 1bc] (2bc) {};
      \node[Node=5ex,fill=tikzcolor13,below=of 1bc] {};
      \node[Node=5ex,fill=tikzcolor15,below=of 2bc] {};
    \end{tikzpicture}
    \caption{block-parallel $\SubA[BP]^{2}$}
  \end{subfigure}
  \begin{subfigure}[t]{0.3\textwidth}
    \centering
  \begin{tikzpicture}[Node/.style = {rectangle, draw,
        anchor=center,node distance=0mm}]
      \node[Node=5ex,fill=tikzcolor1] (1) {};
      \node[Node=5ex,fill=tikzcolor5,right=of 1] (2) {};
      \node[Node=5ex,fill=tikzcolor9,below=of 1] (3) {};
      \node[Node=5ex,fill=tikzcolor13,below=of 2] {};
      \node[Node=5ex,fill=white,right=of 2] (1b) {};
      \node[Node=5ex,fill=white,right=of 1b] (2b) {};
      \node[Node=5ex,fill=white,below=of 1b] {};
      \node[Node=5ex,fill=white,below=of 2b] {};
      \node[Node=5ex,fill=white, below=of 3] (1c) {};
      \node[Node=5ex,fill=white,right=of 1c] (2c) {};
      \node[Node=5ex,fill=white,below=of 1c] {};
      \node[Node=5ex,fill=white,below=of 2c] {};
      \node[Node=5ex,fill=tikzcolor1,right=of 2c] (1bc) {};
      \node[Node=5ex,fill=tikzcolor5,right=of 1bc] (2bc) {};
      \node[Node=5ex,fill=tikzcolor9,below=of 1bc] {};
      \node[Node=5ex,fill=tikzcolor13,below=of 2bc] {};
    \end{tikzpicture}
    \caption{block-global $\SubA[BG]^{2}$}
  \end{subfigure}
  \caption[Schematic representation of different *-subalgebras.]{Schematic
representation of different *-subalgebras.
    Same colors mean the coefficients are coupled.
    White means restricted to zero.
    The top row shows the elementary cases. The bottom row shows the combined
*-subalgebras.
    Inspired by \cite[][Table~3.1]{lund2018block}.}
  \label{fig:subalgebras}
\end{figure}

Figure \ref{fig:subalgebras} shows a schematic representation of
elements in the different *-subalgebras of $\bbR^{4 \times 4}$.
Same colors mean a coupling of the coefficients.
White coefficients are restricted to be zero.

In the work of \citeauthor[]{frommer2017block} the block case is called
\textit{classic}, the parallel method is called \textit{loop-interchange} and
the block-parallel case is called \textit{hybrid}. The block-global method is
not considered in their work.
We chose the new naming scheme because it feels more natural, as we derive the
new cases from the three elementary ones.

To get a feeling for the different *-subalgebras we look at the following
example of computing the normalizer of a block vector $X$ in the
block-global case.
\begin{ex}[Normalizer in the block-global *-subalgebra]
  We can compute the normalizer of the block vector
  $[X_{1},\ldots,X_{q}]$ in the block-global case by computing
  a QR decomposition
  \begin{align}
    \begin{bmatrix}
      X_{1}\\
      \vdots\\
      X_{q}
    \end{bmatrix}
    =
    \begin{bmatrix}
      \tilde{Y}_{1}\\
      \vdots\\
      \tilde{Y}_{q}
    \end{bmatrix}
    \rho, \qquad \rho\in\bbR^{p\times p}.
  \end{align}
  To enforce the normalization we set the normalizer as
  \begin{align} \normalizer[X]{{{\SubA}_{BG}^{p}}}
    &= \Identity[q] \otimes \frac{1}{\sqrt{q}}\rho \in \SubA[BG].
  \end{align}
  Then, the block vector $Y = \sqrt{q}\tilde{Y}$ is normalized as
  \begin{align}
    \sproduct[{\SubA}_{BG}^{p}]{Y}{Y}
    &= \Identity[q] \otimes \sum_{i=0}^{q}\transpose{\tilde{Y}_{i}}\tilde{Y}_{i}\\
    &= \Identity[q]\otimes \Identity[p] = \Identity[s].
  \end{align}
\end{ex}
The following lemma derives directly from Definitions
\ref{defn:elementary-algebras}  and \ref{defn:subalgebras}.
\begin{lem}[Embeddings of *-subalgebras]
  \label{lem:algebraembedding}
  For $p_1, p_2\in\mathbb{N}$, where $p_1$ is a divisor of $p_2$ and $p_2$ is a
  divisor of $s$, we have the following embedding:
  \begin{align}
    \begin{array}{ccccccc}
      \SubA[P] & \subseteq & \SubA[BP]^{p_1} & \subseteq & \SubA[BP]^{p_2} & \subseteq & \SubA[B]\\
      \rotatebox[origin=c]{90}{$\subseteq$} && \rotatebox[origin=c]{90}{$\subseteq$} && \rotatebox[origin=c]{90}{$\subseteq$} && \rotatebox[origin=c]{90}{$=$}\\
      \SubA[G] & \subseteq & \SubA[BG]^{p_1} & \subseteq & \SubA[BG]^{p_2} & \subseteq & \SubA[B]
    \end{array}
  \end{align}
\end{lem}

Note that due to Lemma \ref{lem:krylovembedding} we have the analog
embeddings for the corresponding Krylov spaces.

As a closure of this section we consider the more generic case of
a linear operator $\bm{A}\in L(\bbR^{\dimA \times s}, \bbR^{\dimA \times s})$,
and introduce a classification for this type of operator.
This is in analogy of symmetry and definiteness in the scalar case
($s=1$).

\begin{defn}
  \label{def:bsa_bpd}
  Let $\bm{A} \in L(\bbR^{\dimA \times s}, \bbR^{\dimA \times s})$ be a
  linear operator on $\bbR^{\dimA \times s}$ and $\SubA$ a
  \mbox{*-subalgebra} with a block inner product
  $\sproduct[\SubA]{\cdot}{\cdot}$.
  We call $\bm{A}$
  \begin{itemize}
  \item block self-adjoint (BSA), if for all
    $X,Y\in\bbR^{\dimA\times s}$ holds
    \begin{align}
      \sproduct[\SubA]{\bm{A}X}{Y} = \sproduct[\SubA]{X}{\bm{A}Y}.
    \end{align}
    \begin{samepage}
  \item block positive definite (BPD), if
    \begin{enumerate}[label=\alph*)]
    \item $\bm{A}$ is BSA and for all
      $X\in \bbR^{\dimA\times s}$ with full rank,
      $\sproduct[\SubA]{X}{\bm{A}X}$ is self-adjoint and positive definite and
    \item for all rank-deficient $X\neq 0$, $\sproduct[\SubA]{X}{\bm{A}X}$
      is self-adjoint, positive semi-definite and non-zero.
    \end{enumerate}
  \end{samepage}
\end{itemize}
\end{defn}

\begin{rem}
  Consider the following representation of the operator $\bm{A}$ that operates
  on the vectorization of $\bbR^{\dimA\times s}$,
  \begin{align}
    \bm{\hat{A}} = \begin{bmatrix}
      A_{1,1} & \cdots & A_{1,s}\\
      \vdots & & \vdots\\
      A_{s,1} & \cdots & A_{s, s}
    \end{bmatrix}.
  \end{align}
  That means
  \begin{align}
    (\bm{A}X)_{i} = \sum_{j=1}^{s} A_{ij}X_{j}.
  \end{align}
  Then we distinguish the following cases
  \begin{itemize}
  \item global ($\SubA[G]$):
    \begin{align}
      \bm{A} \text{ is BSA} &&\Leftrightarrow&&&
                              \bm{\hat{A}}
                              \text{ is symmetric}\\
      \bm{A}\text{ is BPD} &&\Leftrightarrow&&&
                              \bm{\hat{A}}
                             \text{ is symmetric positive definite}
    \end{align}
  \item parallel ($\SubA[P]$):
    \begin{align}
      \bm{A}\text{ is BSA} &&\Leftrightarrow&&
      & A_{i,j}=0\quad\forall j\neq i\in \left\{1,\ldots,s\right\}\quad \text{
        and}\\
                           &&&&
      & A_{i,i}\text{ is symmetric}\quad \forall i\in\left\{1,\ldots,s\right\}\\
      \bm{A}\text{ is BPD} &&\Leftrightarrow&&& \bm{A}\text{ is BSA
                                                and}\\
      &&&&&A_{i,i} \text{ is positive definite } \forall i=1,\ldots,s
    \end{align}
  \item block  ($\SubA[B]$):
    \begin{align}
      \bm{A}\text{ is BSA} &&\Leftrightarrow &&
      & A_{i,j}=0 \quad\forall j\neq i\in\left\{1,\ldots,s\right\}\text{ and}\\
                           &&&&
      & A_{i,i}=A_{j,j} \quad\forall j, i\in \left\{1,\ldots,s\right\}\text{ and}\\
                           &&&&
      & A_{i,i}\text{ is symmetric}\quad\forall i\neq j \in\left\{i,\ldots,s\right\}\\
      \bm{A}\text{ is BPD} &&\Leftrightarrow&&& \bm{A}\text{ is BSA
                                                and}\\
      &&&&&A_{i,i} \text{ is positive definite } \forall i=1,\ldots,s
    \end{align}
  \end{itemize}
  As we only consider block linear systems as defined by Equation
  \eqref{eq:blocklinearsystem}, the operator $\bm{A}$ defined by
  \begin{align}
    \left(\bm{A}X\right)_{i} = AX_{i}
  \end{align} is BSA if $A$ is symmetric and BPD if $A$ is symmetric positive
definite for all the mentioned cases.
\end{rem}
Finally we define orthogonally for the block inner product.
\begin{defn}[block orthogonality]
  Let $X,Y\in\bbR^{\dimA\times s}$ be two block vectors and $\bm{A}$ a
  BPD operator.
  We call $X,Y$
  \begin{itemize}
    \item $\SubA$-orthogonal if
    \begin{align}
      \sproduct[\SubA]{X}{Y} &= 0
    \end{align}
    holds.
  \item $\SubA$-$\bm{A}$-orthogonal if
    \begin{align}
      \sproduct[\SubA]{\bm{A}X}{Y} = 0
    \end{align}
    holds.
  \end{itemize}
\end{defn}

After reviewing the theoretical aspects of the block Krylov framework,
we take a look at the practical parts in the next sections.

\section{Implementation}

Our implementation builds upon the SIMD interface in the \cc~software
framework \dune~\cite{bastian2020dune,dune24:16,dunepaperI:08,dunepaperII:08}.
This ensures that we make use of the SIMD capabilities of the hardware and
offers a way to implement horizontal parallelism easily.
SIMD data types behave like a numeric type (e.g.\ \cpp{double}) but process
multiple values at once.
The coefficients of a SIMD data type are called the \textit{lanes}.
The type of one lane is called the \textit{scalar type} of the SIMD data type.

Supported SIMD data types include VCL~\cite{fog2017vcl} and
Vc~\cite{kretz2012vc,kretz2015extending}. Furthermore, \dune provides a
simple fallback implementation, \cpp{Dune::LoopSIMD}, that is based on
static loops and relies on compiler optimization for the exploration
of SIMD instructions.
In the future it is planed to support the SIMD features of the
\cc~standard once the parallelism TSv2\footnote{see for example
  \url{https://en.cppreference.com/w/cpp/experimental/parallelism_2}}
is merged into the standard.

The SIMD interface of \dune unifies the usage of the SIMD specific
operations that differ for different implementations.
The essential components of the interface for a SIMD data type
\cpp{T} are
\begin{itemize}
\item \cpp{Simd::lane(size_t l, T x)}: provides access to a
  single lane
\item \cpp{Simd::lanes()}: provides access to the number of
  lanes (SIMD width)
\item \cpp{SIMD::Scalar<T>}: the scalar data type
  (e.g.\ \cpp{double})
\end{itemize}
Further features include the evaluation of conditional expressions and
the implementation of the math functions in \cc, like \cpp{min},
\cpp{max}, \cpp{sin} etc.

\cpp{Dune::LoopSIMD} is a fallback implementation in \dune.
It inherits from \cpp{std::array} and
implements the SIMD interface by overloading all arithmetic operators.
For example, the implementation of the \cpp{operator+} can be
found in Listing \ref{lst:loopsimd_plus}.
The performance gain of this data type depends on the optimizations of
the compiler.
In particular one problem, in which recent compilers fail, are the use
of FMA operations in expressions like \cpp{a += alpha*b}, if it
involves many lanes.

Another feature of \cpp{Dune::LoopSIMD} is that it can be used to
concatenate another SIMD type to a large one.
For example VCL only implements types with the
hardware SIMD width, i.e.\ $4$ or $8$.
If we want to use larger SIMD types we use \cpp{Dune::LoopSIMD} to
concatenate multiple \cpp{Vec8d} or \cpp{Vec4d}.
For example, \cpp{Dune::LoopSIMD<Vec8d, 4>} is a SIMD data type with $32$ lanes.

\begin{c++}{Implementation of the \cpp{operator+} of \cpp{Dune::LoopSIMD}.}{lst:loopsimd_plus}
template<class T, std::size_t S, std::size_t A>
auto operator +(const LoopSIMD<T,S,A> &v,
                const LoopSIMD<T,S,A> &w) {
  LoopSIMD<T,S,A> out;
  for(std::size_t i=0; i<S; i++){
    out[i] = v[i] + w[i];
  }
  return out;
}
\end{c++}

To use a SIMD data type in a solver a vector type must be specified,
which represents the block vector space $\bbR^{\dimA\times s}$.
This could be achieved by using the SIMD data type as the
\cpp{field_type} in the \cpp{Dune::BlockVector} template.
Practically, this represents a row major storage of the block vectors.
For this vector type a \cpp{Dune::LinearOperator} that
represents the operator $A:\bbR^{\dimA\times s}\to\bbR^{\dimA\times
  s}$ can be implemented for example by a \cpp{Dune::MatrixAdatper}
that takes a sparse matrix and turns it into a linear operator.
A setup of the linear operator type can be seen in Listing \ref{lst:linear_operator_setup}.
Preconditioners can be set up based on the same vector type.

\begin{c++}{Setup of a linear operator type that operates on block%
    vectors based on a sparse matrix.}{lst:linear_operator_setup}
typedef Dune::LoopSIMD<double, 16> field_type;
typedef Dune::BlockVector<field_type> vector_type;
typedef Dune::BCRSMatrix<double> matrix_type;
typedef Dune::MatrixAdapter<matrix_type, vector_type, vector_type> linear_operator_type;
\end{c++}

As we want to be flexible with the choice of the *-subalgebra, the block inner
product is implemented in a generic fashion.
The existing solvers in \dune[ISTL]\cite{dune-istl} are build upon an interface
called \cpp{Dune::ScalarProduct}.
We extend this concept of a scalar product, i.e.\ the return type of the scalar
product is not a scalar, but a \cpp{Block}, which is a type-erasure container
that provides the functionalities shown in Listing \ref{lst:block}.

\begin{c++}{The Block interface used to represent elements in %
    $\bbR^{s\times s}$.}{lst:block}
template<class X>
class Block{
  // the BAXPY operation x += alpha*Y*sigma
  void axpy(const scalar alpha, X& x, const X& y) const;
  Block& invert(); // inverts the block
  Block& transpose(); // transposes the block
  Block& add(const Block& other); // add another block
  Block& scale(const scalar& factor); // scale with a scalar
  // multiplication with other Block
  Block& leftmultiply(const Block& other);
  Block& rightmultiply(const Block& other);
};
\end{c++}
This design enables us to implement the kernels for different
*-subalgebras in a specialized way.
The fallback implementation is the parallel case $\SubA[P]$,
which was the default behavior in \dune before.

We extend the interface further by a function \cpp{inormalizer(X& x)} that
computes and returns the normalizer $\normalizer[X]{\SubA}$ and normalizes the
block vector \cpp{x} with respect to the computed normalizer.
The \texttt{i} as a prefix is inspired by the non-blocking MPI functions and
indicates that the function returns a \cpp{Future} (see Section
\ref{sec:communication_implementation}).
In the sequential case, the normalizer is computed using the
LAPACK~\cite{lapack99} function \cpp{xGEQRF}, which uses Householder
transformations to compute the QR decomposition.
An elaborate discussion about how to compute the normalizer in the parallel case
can be found in Section \ref{sec:tall-skinny-qr}.

\section{Performance Analysis}
Now we look at the performance characteristics of the building blocks that are
needed to build a block Krylov method.
These are
\begin{samepage}
  \begin{itemize}
  \item \texttt{BOP}: Applying the operator $A$
  \item \texttt{BDOT}: Compute the block inner product
  \item \texttt{BAXPY}: Block vector update
  \end{itemize}
\end{samepage}
As already mentioned, the implementation of the normalizer relies
on LAPACK in the sequential case.
Therefore, we do not discuss its performance here.
The performance of the preconditioner depends of course of its choice.
For simplicity, we assume that the preconditioner behaves similar to
\texttt{BOP}.

We assume in this section that the operator is an assembled sparse
matrix in CSR format with $z$ non-zeros. We further assume that the
column index in the CSR format needs as much space in memory as the
coefficient (e.g.\ 64-bit for \cpp{int} and \cpp{double}).
This is the unit in which we denote data size.
The row indices for the sparse matrix are neglected.
This leads to a total memory requirement for the matrix of $2z$.
Together with the input and output block vector $\beta_{\texttt{BOP}}
= 2z + 2s\dimA$ values must be transferred from the main memory to
the registers.
For the \texttt{BOP} operation $\omega_{\texttt{BOP}} = 2sz$
floating-point operations are necessary.
Hence, we get an operational intensity of $\frac{sz}{z+sn}$.
This means the operational intensity is higher (better) for more right-hand
sides $s$ or more non-zeros $z$.

The situation is a bit more sophisticated for the \texttt{BDOT} and
\texttt{BAXPY} kernels. However, both kernels behave quite similar.
Both operate on two block vectors that must be loaded from the main
memory, which yields $\beta_{\texttt{BDOT}} = 2ns$.
The difference between the kernels is that the \texttt{BAXPY} kernel
writes one block vector back to the main memory.
Therefore, we have $\beta_{\texttt{BAXPY}} = 3ns$ memory transfers.
We assume that the data of the *-subalgebra element can be cached and therefore
does not need to be communicated through the memory hierarchy.

The number of floating-point operations depend on the *-subalgebra.
In this analysis we consider the cases $\SubA[BP]^{p}$ and $\SubA[BG]^{p}$
from Definition \ref{defn:subalgebras}.
For both, \texttt{BDOT} and \texttt{BAXPY}, the number of
floating-point operations increases quadratic with $p$ and we have
$\omega_{\texttt{BDOT}} = \omega_{\texttt{BAXPY}} = 2np^{2}q$.
This yields an arithmetic intensity of $p$ for \texttt{BOP} and
$\frac23 p$ for \texttt{BAXPY}.
\begin{table}
  \centering
  \caption[Performance relevant characteristics for the \texttt{BOP},
\texttt{BDOT} and \texttt{BAXPY} kernels.]{Performance relevant characteristics
for the \texttt{BOP}, \texttt{BDOT} and \texttt{BAXPY} kernels.
    The columns denote the number of floating-point operations $\omega$, amount
of data loaded from main memory $\beta$ and the arithmetic intensity
$\frac{\omega}{\beta}$. The number of non-zeros in $A$ are denoted by $z$.
  }
  \renewcommand{\arraystretch}{1.5}
  \begin{tabular}{L{4em}ccc}\toprule
    & $\omega$ & $\beta$ & arith. intensity\\
    \midrule
    \texttt{BOP}
    & $2sz$ & $2z + 2sn$ & $\frac{sz}{z + sn}$\\
    \texttt{BDOT}
    & $2np^2q$ & $2ns$ & $p$\\
    \texttt{BAXPY}
    & $2np^2q$ & $3ns$ & $\frac23 p$\\
    \bottomrule
  \end{tabular}
  \label{tbl:kernels}
\end{table}
Table \ref{tbl:kernels} summarizes the numerical characteristics of the kernels.

This is the great advantage of the presented framework.
The parameter $p$ can be tuned such that the arithmetic intensity
matches the properties of the hardware.
For many right-hand sides $s$ and small $p$ all kernels would be
memory-bound and the costs are independent of $p$, as the amount of
data that must be loaded does not depend on $p$.
Therefore, the parameter $p$ can be chosen as large such that the
$p^{2}$ scaling of the kernels does not have an effect.
Up to that $p$ the faster convergence of the block method comes for
free and the better arithmetical intensity of the \texttt{BOP} kernel
for large $s$ can be preserved.

To achieve the best performance the kernels must be implemented very
carefully.
In particular, one must ensure a good data locality.
For our implementation, we iterate over the rows of the block vectors
in chunks of $4$ rows.
We found this number experimentally and suppose that the optimal
number depends on the number of registers of the CPU and for how many
cycles a FMA operation occupies the registers.
Within these chunks we iterate over the rows and compute
the corresponding matrix-matrix products.
The implementations of the matrix-matrix products are shown in
Listing \ref{lst:mm}.
This approach was already presented by
\citeauthor[]{stewart2008block}~\cite[]{stewart2008block}.

\begin{c++}{Implementation of the inner matrix-matrix products. Block %
    vector rows are iterated in chunks of size \lstinline{ChunkSize} %
    to increase the arithmetical intensity.}{lst:mm}
// computes c += a^t b
template<class SIMD, size_t ChunkSize>
void mtm(const std::array<SIMD, ChunkSize>& a,
         const std::array<SIMD, ChunkSize>& b,
         std::array<SIMD, lanes<SIMD>()>& c){
  for(size_t i=0; i<ChunkSize; ++i){
    for(size_t j=0; j<lanes<SIMD>(); ++j){
      c[j] += lane(j, a[i])*b[i];
    }
  }
}

// computes c += a b
template<class SIMD, size_t ChunkSize>
void mm(const std::array<SIMD, ChunkSize>& a,
        const std::array<SIMD, lanes<SIMD>()>& b,
        std::array<SIMD, ChunkSize>& c){
  for(size_t i=0; i<lanes<SIMD>(); ++i){
    for(size_t j=0; j<ChunkSize; ++j){
      c[j] += lane(i, a[j])*b[i];
    }
  }
}
\end{c++}

\section{Numerical Experiments}
To compare the different methods in practice, we executed several
tests.
In this chapter, all tests are carried out on our compute server, which is
an Intel Skylake-SP Xeon Gold 6148 with $\SI{377}{\giga\byte}$ main memory.
To make the results as reproducible as possible, we deactivate the
turbo mode.
In the described setting the system has a theoretical peak performance of
$\SI{76.8}{\giga\flop\per\second}$ $\left(
  =\SI{2.4}{\giga\hertz} * \SI{32}{\flop \per \cycle}\right)$ using
one core.
Measurements show a memory bandwidth of
$\SI{13.34}{\giga\byte\per\second}$, measured with the
\texttt{daxpy} benchmark of the likwid-bench suite~\cite{treibig2010likwid}.

Multi-core tests are executed on $20$ cores of the machine, which is
one NUMA node.
In this setting the frequency reduces to $\SI{2.2}{\giga\hertz}$ when
using AVX-512 instructions,
leading to a theoretical peak performance of
$\SI{1408}{\giga\flop\per\second}$
$\left( = 20 * \SI{2.2}{\giga\hertz} * \SI{32}{\flop\per\cycle}\right)$.
As the cores share the same memory connection, the memory bandwidth
does not scale with the number of cores.
The measured memory bandwidth with $20$ cores is
$\SI{98.47}{\giga\byte\per\second}$, also measured with the
\texttt{daxpy} benchmark.

In a first test series we compare the run-times of the building blocks
\texttt{BOP}, \texttt{BAXPY} and \texttt{BDOT} for the block-parallel and
block-global methods.
\begin{figure}[t]
  \begin{subfigure}[t]{0.5\textwidth}
    \centering
    \resizebox{\textwidth}{!}{
      \input{generated/block_benchmark/bop_2DFD_seq.pgf}
    }
    \caption{\texttt{BOP} with 2D Finite-Differences matrix.}
  \end{subfigure}
    \begin{subfigure}[t]{0.5\textwidth}
    \centering
    \resizebox{\textwidth}{!}{
      \input{generated/block_benchmark/bop_3DQ1_seq.pgf}
    }
    \caption{\texttt{BOP} with 3D Q1-Finite-Elements matrix.}
  \end{subfigure}\\[5ex]
    \begin{subfigure}[t]{0.5\textwidth}
    \centering
    \resizebox{\textwidth}{!}{
      \input{generated/block_benchmark/bdot_seq_fixK.pgf}
    }
    \caption{\texttt{BDOT} for the block-parallel and block-global case.}
  \end{subfigure}
    \begin{subfigure}[t]{0.5\textwidth}
    \centering
    \resizebox{\textwidth}{!}{
      \input{generated/block_benchmark/baxpy_seq_fixK.pgf}
    }
    \caption{\texttt{BAXPY} for the block-parallel and block-global case.}
  \end{subfigure}
  \caption[Microbenchmarks for kernels \texttt{BOP}, \texttt{BDOT} and
    \texttt{BAXPY} executed on one core.]{Microbenchmarks for kernels \texttt{BOP}, \texttt{BDOT} and
    \texttt{BAXPY} executed on one core. Crosses mark the measured
    data. Dotted lines mark the memory bound. Dashed lines mark the
    compute bound.}
  \label{fig:microbenchmarks}
\end{figure}
Figure \ref{fig:microbenchmarks} shows the runtimes of the kernels per
right-hand side in the single-core case.
\begin{figure}[t]
  \begin{subfigure}[t]{0.5\textwidth}
    \centering
    \resizebox{\textwidth}{!}{
      \input{generated/block_benchmark/bop_2DFD_para.pgf}
    }
    \caption{\texttt{BOP} with 2D Finite-Differences matrix.}
  \end{subfigure}
    \begin{subfigure}[t]{0.5\textwidth}
    \centering
    \resizebox{\textwidth}{!}{
      \input{generated/block_benchmark/bop_3DQ1_para.pgf}
    }
    \caption{\texttt{BOP} with 3D Q1-Finite-Elements matrix.}
  \end{subfigure}\\[5ex]
    \begin{subfigure}[t]{0.5\textwidth}
    \centering
    \resizebox{\textwidth}{!}{
      \input{generated/block_benchmark/bdot_para_fixK.pgf}
    }
    \caption{\texttt{BDOT} for the block-parallel and block-global case.}
  \end{subfigure}
    \begin{subfigure}[t]{0.5\textwidth}
    \centering
    \resizebox{\textwidth}{!}{
      \input{generated/block_benchmark/baxpy_para_fixK.pgf}
    }
    \caption{\texttt{BAXPY} for the block-parallel and block-global case.}
  \end{subfigure}
  \caption[Microbenchmarks for kernels \texttt{BOP}, \texttt{BDOT} and
    \texttt{BAXPY} executed on 20 cores.]{Microbenchmarks for kernels \texttt{BOP}, \texttt{BDOT} and
    \texttt{BAXPY} executed on 20 cores. Crosses mark the measured
    data. Dotted lines mark the memory bound. Dashed lines mark the
    compute bound.}
  \label{fig:microbenchmarks_parallel}
\end{figure}
Figure \ref{fig:microbenchmarks_parallel} shows the run-time of the
same kernels in the multi-core case.
We plotted the execution time per right-hand side ($t/s$).
For the \texttt{BOP} kernel we carried out the tests for different
values of $s$.
We tested two different matrix patterns.
One is a very sparse one, resulting from a 2D finite differences
discretization of a Poisson problem on a $1000\times 1000$ grid with
$z=5$ non-zeros per row.
The other one results from a 3D Q1 finite element discretization on a
$100\times 100\times 100$ grid with $z=27$ non-zero coefficients per
row.
For the SIMD interface we used the VC library and combine it with
\cpp{Dune::LoopSIMD} to assemble larger SIMD data types as described
in the previous section.

For the \texttt{BDOT} and \texttt{BAXPY} kernel we used $s=256$ and
carried out the tests for different $p$.
In the one-core test case we used $\dimA=\num{500000}$ and in the multi-core
test case we used $\dimA=\num{6000000}$.
Further numerical tests show that the run-time of these kernels scale
linearly with $s$.

We see that the measured behavior of the kernels matches our
theoretical expectation.
For the \texttt{BOP} kernel it turns out that for all $s$ the kernel
is memory bound and it performs more efficient with larger $s$.
We suppose that the slight increase of the runtime per $s$ for larger
$s$ in the 2D finite differences case is due to cache effects, as fewer
rows of the block vectors can be cached.

For the \texttt{BDOT} and \texttt{BAXPY} kernels we see that the
runtime is memory bound for $p \lesssim 16$ in the one-core case and
memory bound for $p \lesssim 64$ in the multi-core case.
In particular the runtime does not depend on $p$ in this
regime.
For larger $p$ the run-time increases quadratically, as expected.
The runtime of the block-global and block-parallel method does not
differ.
This was predicted by the theory as well.
We already published similar results in \cite[]{dreier2019strategies}.


\chapter{Block Conjugate Gradients Method}
\label{chap:blockcg}

For the solution of large sparse symmetric positive definite linear systems,
the Conjugate Gradients method combined with a proper preconditioner is the
method of choice.
In the following, we reformulate the block Conjugate Gradients (BCG) method
as proposed by \citeauthor{oleary1980block}~\cite[]{oleary1980block}
based on the general framework presented in the last chapter
and introduce a novel adaptive stabilization technique based on the paper by
\citeauthor{dubrulle2001retooling}~\cite{dubrulle2001retooling}.

We consider $A$ as a symmetric positive definite operator, as
this is a requirement of the BCG method.
Furthermore, in this chapter $\Prec\in L(\bbR^{\dimA},\bbR^{\dimA})$
denotes a symmetric positive definite preconditioner.

\section{Formulation of the Block Conjugate Gradients Method}

The objective of the BCG method in the $k$th iteration is to find an approximation
$X^{k} \in X^{0} + \krylov[\SubA]{\inverse{\Prec}A}{\inverse{\Prec}R^{0}}{k}$,
which minimizes the block energy error
\begin{align}
  \label{eq:energy-minimizer}
  X^{k} = \argmin_{Y\in
  X^{0} + \krylov[\SubA]{\inverse{\Prec}A}{\inverse{\Prec}R^{0}}{k}} \|Y-\blocksolution\|_{A,F}.
\end{align}

Before characterizing this minimization property in more detail, we look at a small
auxiliary lemma.
\begin{lem}
  \label{lem:trace_aux}
  Let $\eta \in \SubA$ with
  \begin{align}
    \trace{\eta \sigma} &= 0&\forall \sigma \in \SubA.
  \end{align}
  Then $\eta = 0$ must hold.
\end{lem}
\begin{proof}
  Assume that $\eta \neq 0$ and choose $\sigma = \transpose{\eta}$. It
  follows
  \begin{align}
    0 = \trace{\eta \sigma} = \trace{\eta \transpose{\eta}} > 0.
  \end{align}
  Because $\eta \transpose{\eta} \neq 0$ is a positive semi-definite matrix, for
  which the trace is the sum of its eigenvalues.
\end{proof}

With this lemma, we formulate the following theorem.
\begin{thm}
  \label{thm:minimizing-orthogonality}
  The minimization property \eqref{eq:energy-minimizer} is equivalent
  to the orthogonality condition
  \begin{align}
  \label{eq:minimizing-orthogonality}
    \sproduct[\SubA]{R^{k}}{X} &= 0&\forall
 X\in\krylov[\SubA]{\inverse{\Prec}A}{\inverse{\Prec}R^{0}}{k},
  \end{align}
  where $R^{k} = AX^{k} -B$ is the residual for the approximation $X^{k}$.
\end{thm}
\begin{proof}
  For any $X\in
  \krylov[\SubA]{\inverse{\Prec}A}{\inverse{\Prec}R^{0}}{k}$, we define
  the coercive functional
  \begin{align}
    \mcJ_{X}(\varepsilon) &= \frac12 \norm[A,F]{X^{k} +
                            \varepsilon X - \blocksolution}^{2}\\
                      &= \frac12 \trace{\sproduct[\SubA]{X^{k} +\varepsilon X
                        - \blocksolution}{R^{k}  +\varepsilon AX}}.
  \end{align}
  Note that the minimization of \eqref{eq:energy-minimizer} is
  equivalent to the minimization of $\mcJ_{X}$ for all
  $X\in\krylov[\SubA]{\inverse{\Prec}A}{\inverse{\Prec}R^{0}}{k}$.
  Due to the linearity of the trace and the block inner product, the
  differential of $\mcJ_{X}$ computes
  \begin{align}
    D_{\varepsilon}\mcJ_{X}(\varepsilon) &=
                                       \trace{\sproduct[\SubA]{R^{k}}{X}}
                                       + \varepsilon \trace{\sproduct[\SubA]{X}{AX}}.
  \end{align}
  Here we used that $\trace{\sproduct[\SubA]{X}{Y}} =
  \trace{\transpose{\sproduct[\SubA]{Y}{X}}} = \trace{\sproduct[\SubA]{Y}{X}}$.
  For the first implication, we assume that $X^{k}$ is a minimizer of
  \eqref{eq:energy-minimizer}. Hence, we have
  \begin{align}
    0 &= D_{\varepsilon}\mcJ_{X\sigma}(0) =
        \trace{\sproduct[\SubA]{R^{k}}{X\sigma}}
    &\forall X\in\krylov[\SubA]{\inverse{\Prec}A}{\inverse{\Prec}R^{0}}{k}, \sigma\in\SubA.
  \end{align}
  From Lemma \ref{lem:trace_aux}, we obtain
  \begin{align}
    \sproduct[\SubA]{R^{k}}{X} &= 0 & \forall
                                      X\in\krylov[\SubA]{\inverse{\Prec}A}{\inverse{\Prec}R^{0}}{k}.
  \end{align}
  For the other implication, we assume that
  $\sproduct[\SubA]{R^{k}}{X}=0$ for all
  $X\in\krylov[\SubA]{\inverse{\Prec}A}{\inverse{\Prec}R^{0}}{k}$.
  This yields
  \begin{align}
    D_{\varepsilon}\mcJ_{X}(\varepsilon)=\varepsilon\trace{\sproduct[\SubA]{X}{AX}}.
  \end{align}
  Hence, $D_{\varepsilon}\mcJ_{X}(0)$ vanishes for all $X \in \krylov[\SubA]{\inverse{\Prec}A}{\inverse{\Prec}R^{0}}{k}$.
  Thus, $X^{k}$ is a minimizer of \eqref{eq:energy-minimizer}.
\end{proof}

Now we deduce the formulas for the method.
The method computes an $A$-$\SubA$ orthogonal basis
$\left\{P^{j}\right\}_{j=0}^{k}$ of the block Krylov space
$\krylov[\SubA]{\inverse{\Prec}A}{\inverse{\Prec}R^{0}}{k+1}$.
This basis is used to update the initial guess and residual iteratively
\begin{align}
  \label{eq:cg_update_x}
  X^{k+1} &= X^{k} + P^{k}\lambda^{k}\\
  \label{eq:cg_update_r}
  R^{k+1} &= R^{k} -AP^{k}\lambda^{k}.
\end{align}

To compute the coefficient $\lambda^{k}\in\SubA$, let an
$A$-$\SubA$-block orthogonal basis $\{P^{j}\}_{j=0}^{k}$ be given and let
$X^{k}$ be a minimizer of \eqref{eq:energy-minimizer}.
By Theorem \ref{thm:minimizing-orthogonality}, we obtain the minimizer
of \eqref{eq:energy-minimizer} for the following block Krylov space by
\begin{align}
  0 &= \sproduct[\SubA]{R^{k+1}}{X} & \forall X \in
 \krylov[\SubA]{\inverse{\Prec}A}{\inverse{\Prec}R^{0}}{k+1}\\
    &= \sproduct[\SubA]{R^{k}}{X} -
      \sproduct[\SubA]{AP^{k}\lambda^{k}}{X}.\\
  \intertext{For $X$ we choose the basis $\{P^{j}\}_{j=0}^{k}$ and get}
  0 &= \sproduct[\SubA]{R^{k}}{P^{j}} -
      \sproduct[\SubA]{AP^{k}\lambda^{k}}{P^{j}}  & \forall j=0,\ldots,k.\
  \intertext{Due to the $A$-$\SubA$-orthogonality, this yields}
  \label{eq:rpzero}
  0 &= \sproduct[\SubA]{R^{k}}{P^{j}} & \forall j=0,\ldots,k-1\\
  \intertext{and}
  \label{eq:lambda_def}
  \lambda^{k} &= \inverse{\left(\sproduct[\SubA]{P^{k}}{AP^{k}}\right)}
  \sproduct[\SubA]{P^{k}}{R^{k}}.
\end{align}
By Theorem \ref{thm:minimizing-orthogonality}, Equation \eqref{eq:rpzero} holds,
as $X^{k}$ is a minimizer in the Krylov space $\krylov[\SubA]{A}{R^{0}}{k}$.
We use Equation \eqref{eq:lambda_def} as a definition for $\lambda^{k}$.

The next basis vector $P^{k+1}$ is then obtained by $A$-$\SubA$-orthogonalizing
the preconditioned residual $\inverse{\Prec}R^{k+1}$ against the previous basis
vectors, it reads
\begin{align}
  P^{k+1} &= \inverse{\Prec}R^{k+1}
            - \sum_{j=0}^{k} P^{j}\inverse{\left(\sproduct[\SubA]{P^{j}}{AP^{j}}\right)}
            \sproduct[\SubA]{P^{j}}{A\inverse{\Prec}R^{k+1}}.
\end{align}
Next, we show that for $j=0,\ldots,k-1$, the coefficient in the
orthogonalization vanishes.
As $A$ and
$\Prec$ are symmetric and $\inverse{\Prec}AP^{j} \in
\krylov[\SubA]{\inverse{\Prec}A}{\inverse{\Prec}R^{0}}{j+2}$, we have
\begin{align}
  \sproduct[\SubA]{P^{j}}{A\inverse{\Prec}R^{k+1}}
  &= \sproduct[\SubA]{\inverse{\Prec}AP^{j}}{R^{k+1}} = 0,
\end{align}
by using the orthogonality from Theorem \ref{thm:minimizing-orthogonality}.
Hence, we get the update formula
\begin{align}
  \label{eq:cg_update_p}
  P^{k+1} &= M^{-1}R^{k+1} + P^{k}\beta^{k}\\
  \intertext{with}
  \beta^{k} &= \inverse{\left(\sproduct[\SubA]{P^{k}}{AP^{k}}\right)}
              \sproduct[\SubA]{\inverse{\Prec}AP^{k}}{R^{k+1}}.
\end{align}

To reduce the number of block inner products, we reformulate
the coefficients $\lambda^{k}$ and $\beta^{k}$ as follows
\begin{align}
  \lambda^{k}
  &=\inverse{\left(\sproduct[\SubA]{P^{k}}{AP^{k}}\right)}
    \sproduct[\SubA]{P^{k}}{R^{k}}\\
  &=\inverse{\left(\sproduct[\SubA]{P^{k}}{AP^{k}}\right)}\sproduct[\SubA]{\inverse{\Prec}R^{k}
    - P^{k-1}\beta^{k-1}}{R^{k}}\\
  \label{eq:cg_lambda}
 &=\inverse{\left(\sproduct[\SubA]{P^{k}}{AP^{k}}\right)} \sproduct[\SubA]{\inverse{\Prec}R^{k}}{R^{k}},\\
  \beta^{k}
  &= \inverse{\left(\sproduct[\SubA]{P^{k}}{AP^{k}}\right)}
              \sproduct[\SubA]{\inverse{\Prec}AP^{k}}{R^{k+1}}\\
  &=
    \inverse{\left(\sproduct[\SubA]{P^{k}}{AP^{k}}\right)}\transinv{\lambda^{k}}
    \sproduct[\SubA]{R^{k}-R^{k+1}}{\inverse{\Prec}R^{k+1}}\\
  \label{eq:cg_beta}
  &= \inverse{\left(\sproduct[\SubA]{\inverse{\Prec}R^{k}}{R^{k}}\right)}
    \sproduct[\SubA]{\inverse{\Prec}R^{k+1}}{R^{k+1}}.
\end{align}
Here, we used the orthogonality relation
\eqref{eq:minimizing-orthogonality} and the update formula for the
residual \eqref{eq:cg_update_r}.
This reduces the necessary block inner products to
\begin{align}
  \label{eq:cg_alpha}
  \alpha^{k} &= \sproduct[\SubA]{P^{k}}{AP^{k}}\\
  \intertext{and}
  \label{eq:cg_rho}
  \rho^{k} &= \sproduct[\SubA]{\inverse{\Prec}R^{k}}{R^{k}}.
\end{align}

\begin{algorithm}[t]
  \caption{Block Conjugate Gradients Method}
  \label{alg:bcg}
  \begin{algorithmic}
    \State $R^0 = B-AX^0$
    \State $P^0 = \inverse{\Prec}R^0$
    \State $\rho^{0} = \sproduct[\SubA]{P^{0}}{R^{0}}$
    \For{$k = 0,\ldots$ until convergence}
    \State $Q^k = AP^k$
    \State $\alpha^k = \sproduct[\SubA]{P^k}{Q^k}$
    \State $\lambda^k = \inverse{\left(\alpha^k\right)}\rho^k$
    \State $X^{k+1} = X^{k} + P^k\lambda^k$
    \State $R^{k+1} = R^{k} - Q^k\lambda^k$ \label{algl:residual}
    \State break if $\|R^{k+1}\| < \tol$
    \State $Z^{k+1} = \inverse{\Prec}R^{k+1}$
    \State $\rho^{k+1} = \sproduct[\SubA]{Z^{k+1}}{R^{k+1}}$
    \State $\beta^k = \inverse{\left(\rho^{k}\right)}\rho^{k+1}$
    \State $P^{k+1} = Z^{k+1} + P^k\beta^k$
    \EndFor
  \end{algorithmic}
\end{algorithm}

Putting together Equations \eqref{eq:cg_update_x},
\eqref{eq:cg_update_r}, \eqref{eq:cg_update_p}, \eqref{eq:cg_lambda},
\eqref{eq:cg_beta}, \eqref{eq:cg_alpha} and \eqref{eq:cg_rho}, we
obtain Algorithm \ref{alg:bcg}.

In Algorithm \ref{alg:bcg}, we choose $\|R^{k}\| < \tol$ as a break criteria,
where we do not specify which norm is used.
A natural choice would be the Frobenius norm which is the Euclidean
norm on the block vector space.
However, another possibility would be to choose the maximum column
norm
\begin{align}
  \|R^{k}\|_{\infty} := \max\left\{\|R^{k}_{i}\|_{2}\,\big |\, i=1,\ldots,s \right\}.
\end{align}
That ensures that the residual norm of each column is smaller than
$\tol$.
In all our numerical tests, we used the latter, as this is the desired condition
in most applications.

\section{Convergence}
\label{sec:cg_convergence}
We start with a general result that holds for all choices of
*-subalgebras $\SubA$.
From that result, we derive statements about the convergence in the
elementary cases.
These statements can be combined to obtain statements about the
convergence of the combined *-subalgebras defined in Definition~\ref{defn:elementary-algebras}.

\begin{lem}[Generic convergence result]
  \label{lem:bcg_convergence_inf}
  For the error $E^{k}$ of the $k$th step of the BCG method, the
  estimation
  \begin{align}
    \label{eq:cg_convergence_inf}
    \|E^{k}\|_{A,F} \leq \inf_{\mathcal{Q}^{k}}
    \|\mathcal{Q}^{k}(\inverse{\Prec}A)\circ E^{0}\|_{A,F}
  \end{align}
  holds,
  where the infimum is taken over all $\SubA$-valued polynomials
  $\mathcal{Q}^{k}\in\bb{P}_{\SubA}^{k}$ of degree $k$ with absolute coefficient
  $\Identity$.
\end{lem}

\begin{proof}
  We use $\SubA$-valued polynomials to represent the
  energy error.
  By Lemma \ref{lem:blockkrylov_polynomial},
  we can represent the $k$th error of the BCG method as
  \begin{align}
    E^{k} &= \blocksolution - X^{k}\\
          &= \blocksolution - X^{0} - \mathcal{P}^{k}(\inverse{\Prec}A)\circ \inverse{\Prec}R^{0}\\
          &= \left(\Identity -
            \mathcal{P}^{k}(\inverse{\Prec}A)\inverse{\Prec}A\right)\circ E^{0},
  \end{align}
  for some $\SubA$-valued polynomial $\mathcal{P}^{k} \in
  \bb{P}_{\SubA}^{k-1}$.
  With $\mathcal{Q}^{k}(\polyX) = \Identity - \mathcal{P}^{k}(\polyX)\polyX$ and
  taking the $A$-Frobenius norm, we obtain
  \begin{align}
    \|E^{k}\|_{A,F} &= \|\mathcal{Q}^{k}(\inverse{\Prec}A)
                      \circ E^{0}\|_{A,F}.
  \end{align}
  As the BCG method minimizes the $A$-Frobenius norm in the Krylov
  space, we get
  \begin{align}
        \|E^{k}\|_{A,F} &= \inf_{\mathcal{Q}^{k}}\|\mathcal{Q}^{k}(\inverse{\Prec}A)
                          \circ E^{0}\|_{A,F}
  \end{align}
  by taking the infimum of all $\mathcal{Q}^{k}\in\bb{P}^{k}_{\SubA}$ of this shape.
\end{proof}

The next lemma gives a concrete error bound for
all *-subalgebras that we consider in this work.
It is the generalization of Theorem \ref{thm:cg_convergence}.
However, the given bound is not sharp in all cases.
\begin{lem}
  Let $\SubA$ be a *-subalgebra of $\bbR^{s\times s}$, that contains the
  identity of $\bbR^{s\times s}$.
  Then we have
  \begin{align}
    \|E^{k}\|_{A,F} \leq
    2\left(\frac{\sqrt{\kappa} -1 }{\sqrt{\kappa} +1}\right)^{k}\|E^{0}\|_{A,F},
  \end{align}
  where $\kappa$ denotes the condition number of the preconditioned
  operator ${\Prec^{-\frac12}A\Prec^{-\frac12}}$.
\end{lem}
\begin{proof}
  We are using Theorem \ref{lem:bcg_convergence_inf} and choose
  \begin{align}
    \mathcal{Q}^{k}(\polyX) = \tilde{T}^{k}(\polyX)\Identity
  \end{align}
  in Equation \eqref{eq:cg_convergence_inf}, where $\tilde{T}^{k}$ are the
  scaled Chebyshev polynomials defined by Equation \eqref{eq:chebyshev}, scaled with
  respect to the eigenvalues $\lambda_{\min}$ and $\lambda_{\max}$ of the
  symmetrically preconditioned operator $\Prec^{-\frac12}A\Prec^{-\frac12}$.
  As this operator is symmetric positive-definite, it is similar to the diagonal
  matrix of its eigenvalues $\Lambda$, denoted by
  \begin{align}
    \Prec^{-\frac12}A\Prec^{-\frac12} = V\Lambda \inverse{V},
  \end{align}
  where $V$ is the orthonormal matrix of the eigenvectors.
  Similar to the proof of Theorem \ref{thm:cg_convergence}, we compute
  \begin{align}
    \|\tilde{T}^{k}(\inverse{\Prec}A)\circ E^{0}\|_{A,F}^{2}
    &= \trace{\transpose{\left(\tilde{T}^{k}(\inverse{\Prec}A) E^{0}\right)}
      A\tilde{T}^{k}(\inverse{\Prec}A) E^{0}}\\
    &= \trace{\transpose{\left(\tilde{T}^{k}(\Lambda)\inverse{V}
      \Prec^{\frac12}E^{0}\right)}
      \Lambda\tilde{T}^{k}(\Lambda)\inverse{V}\Prec^{\frac12}E^{0}}\\
    &= \trace{\transpose{\left(\inverse{V}\Prec^{\frac12}E^{0}\right)}
      \Lambda^{\frac12}\tilde{T}^{k}(\Lambda)^{2}\Lambda^{\frac12}
      \inverse{V}\Prec^{\frac12}E^{0}}\\
    &\leq \max_{i=0}^{\dimA} |\tilde{T}^{k}\left(\lambda_{i}\right)|^{2}\,
      \trace{\transpose{E^{0}}A E^{0}}\\
    &\leq 4\left(\frac{\sqrt{\kappa}-1}{\sqrt{\kappa}+1}\right)^{2k}
      \|E^{0}\|_{A,F}^{2}.
  \end{align}
  Applying the square-root completes the proof.
\end{proof}

This theorem also applies for the block *-subalgebra $\SubA[B]$.
However, in this case the estimation can be improved.
O'Leary showed the following convergence result to estimate the error of the
classical BCG method.
\pagebreak
\begin{thm}[Convergence of the block Conjugate Gradients Method {\cite[Theorem~5]{oleary1980block}}]
  \label{thm:blockcg}
  For the energy-error of the $i$th column in the $k$th iteration
  $\|E^k_{i}\|_A$ of the BCG method, the following
  estimation holds:
  \begin{align*}
    \norm[A]{E^k_i} &\leq c_1 \mu^k\\
    \text{with }
    \mu = \frac{\sqrt{\cond_s}-1}{\sqrt{\cond_s}+1}, \cond_s
                    &= \frac{\lambda_\dimA}{\lambda_s} \text{ and some constant } c_1>0,
  \end{align*}
  where $\lambda_1\leq \ldots \leq \lambda_\dimA$ denote the eigenvalues of the
  preconditioned matrix ${\Prec^{-\frac12}A\Prec^{-\frac12}}$.
  The constant $c_1$ depends on $s$ and the initial error $E^0$ but
  not on $k$ or $i$.
\end{thm}

The proof also makes use of on Lemma \ref{lem:bcg_convergence_inf}
but the construction of the polynomials is much more sophisticated
and technical.
As we want to concentrate on the practical aspects in this work, we
refer the reader to \cite[]{oleary1980block} for the rigorous proof.

The theorem holds for the classical BCG method ($\SubA[B]$).
However, as the block-parallel method is only a data-parallel version
of the block method the same convergence rate holds with
$s=p$ for the $\SubA[BP]^p$ method, it reads
\begin{align}
  \mu = \frac{\sqrt{\kappa_{p}}-1}{\sqrt{\kappa_{p}}+1}.
\end{align}

The following lemma gives us a convergence rate for the block-global method.
\begin{lem}[Theoretical convergence rate of the block-global method]
  \label{lem:global_cg_convergence}
  The theoretical convergence rate of a block-global method using
  $\SubA[BG]^p$ is
  \begin{align}
    \hat{\mu} =  \frac{\sqrt{\hat{\kappa}_p}-1}{\sqrt{\hat{\kappa}_p}+1},
    \quad\text{with}\quad\hat{\kappa}_p = \frac{\lambda_\dimA}{\lambda_{\left\lceil
        \frac{p}{q} \right\rceil}}.
  \end{align}
\end{lem}
\begin{proof}
    A block-global method is equivalent to solve the $qn$-dimensional system
  \begin{align}
    \begin{pmatrix}
      A\\
      & A\\
      && \ddots\\
      &&& A\\
    \end{pmatrix}
    \begin{pmatrix}
      X_1 & \cdots & X_p\\
      X_{p+1} & \cdots & X_{2p}\\
      &\vdots\\
      X_{s-p+1} & \cdots & X_s
    \end{pmatrix}
                           =
    \begin{pmatrix}
      B_1 & \cdots & B_p\\
      B_{p+1} & \cdots & B_{2p}\\
      &\vdots\\
      B_{s-p+1} & \cdots & B_s
    \end{pmatrix}
  \end{align}
  with the classical block Krylov method with $p$ right-hand sides.
  The matrix of this system has the same eigenvalues as $A$ but with $q$ times
  the multiplicity.
  Thus, the $p$-smallest eigenvalue is $\lambda_{\left\lceil \frac{p}{q} \right\rceil}$.
  Therefore and by applying Theorem \ref{thm:blockcg}, we deduce the theoretical
  convergence rate.
\end{proof}
This result makes the block-global methods irrelevant for practical use.
In particular for $q>1$, the block-parallel method would perform
better while the building blocks are similarly expensive, as we have seen in
Chapter \ref{chap:blockkrylovframework}.

\section{Residual Re-Orthonormalization}
\label{ssec:blockcg_residual_ortho} Algorithm \ref{alg:bcg} requires that
$\rho^{k}=\sproduct[\SubA]{Z^{k}}{R^{k}}$ and
$\alpha^{k}=\sproduct[\SubA]{Q^{k}}{P^{k}}$ are invertible.
It is the case when the residual $R^{k}$ has full-rank.
In the scalar case, $s=1$, this is not a problem.
If the residual is rank-deficient, the linear system has been solved.

In the case, $s>1$, however, this leads to severe problems.
An interpretation of the rank deficiency of the residual is that a linear
combination is converged, because there exists a vector $y\in\bbR^{s}$, such
that
\begin{align}
  0 &= R^{k}y = AX^{k}y - By &\Leftrightarrow &&AX^{k}y &= By.
\end{align}
In exact arithmetic, this case appears in the iteration $\xi_{\SubA}(A,R^{0})$
defined by Equation \eqref{eq:blockkrylov_xi}.
Initially, \citeauthor[]{oleary1980block}~\cite[]{oleary1980block}
suggested to remove dependent vectors and continue the iteration with
a smaller block size.
This strategy is called \textit{deflation} and has two disadvantages
from our perspective.
Firstly, a numerical tolerance parameter must be introduced to check
for the \textit{numerical} rank-deficiency of the residual.
\citeauthor[]{langou2003iterative}~\cite[]{langou2003iterative}
showed in his PhD thesis that a badly chosen parameter could lead
to instabilities or slow down the convergence.
Secondly, we want to choose the block width $s$ as a multiple of the SIMD
width to facilitate SIMD-vectorization.
Deflating the system would change the block width, such that some
effort is needed to handle it in our SIMD setting and we would loose the
performance benefits from the exploration of the SIMD instructions.

\citeauthor{dubrulle2001retooling}~\cite[]{dubrulle2001retooling}
presented multiple approaches to mitigate the stabilization issues
without decreasing the block size.
The most promising approach is the orthonormalization of the residual in
every iteration by computing a QR decomposition.
The algorithm is very elegant without preconditioning, as $\rho$
simplifies to the identity.
As we use preconditioning, this does not hold anymore, as the
$\Prec$-product is used to compute $\rho$.
Therefore, we either need to make the orthonormalization with respect to
the $\Prec$-product, or $\rho$ must be computed explicitly after the
orthonormalization of the residual, which
needs an additional global communication.
We decided to use the normalizer for the
orthonormalization and compute $\rho$ explicitly thereafter,
which is the reason why we defined the normalizer also for
rank-deficient block vectors.

We store the transformation from the orthonormal residual $\bar{R}^{k}$ to the real
residual in the variable $\sigma^{k}\in\SubA$.
It can be updated as
\begin{align}
  \sigma^{k} = \gamma^{k}\sigma^{k-1},
\end{align} where $\gamma^{k}$ is the normalizer of the updated residual
$\tilde{R}^{k} = \bar{R}^{k-1}-Q^{k}\lambda^{k}$, i.e.
\begin{align}
  \tilde{R}^{k} = \bar{R}^{k}\gamma^{k}.
\end{align}
This transformation is then also used to update the solution
\begin{align}
  X^{k} = X^{k-1} + P^{k-1}\lambda^{k-1}\sigma^{k-1},
\end{align} because the search direction $P^{k-1}$ is obtained from the
transformed residual $\bar{R}^{k-1}$.
For the same reason, we must consider the normalizer in the orthogonalization
coefficient $\beta^{k}$, because $P^{k}$ and $P^{k-1}$ are transformed with
respect to $\sigma^{k}$ and $\sigma^{k-1}$, respectively.
Thus, we have
\begin{align}
  \beta^{k} = \left(\rho^{k-1}\right)^{-1}\transpose{\gamma^{k}}\rho^{k}.
\end{align}
The resulting algorithm is shown in Algorithm \ref{alg:bcgro}.
\begin{algorithm}[tbp!]
  \caption{BCG Method with Residual Re-Orthonormalization}
  \label{alg:bcgro}
  \begin{algorithmic}
    \State $R^0 = B-AX^0$
    \State $\bar{R}^{0}, \sigma^{0} = \normalizer[R^{0}]{\SubA}$
    \State $P^0 = \inverse{\Prec}\bar{R}^0$
    \State $\rho^{0} = \sproduct[\SubA]{P^{0}}{\bar{R}^{0}}$
    \For{$k = 0,\ldots$ until convergence}
    \State $Q^k = AP^k$
    \State $\alpha^k = \sproduct[\SubA]{P^k}{Q^k}$
    \State $\lambda^k = \inverse{\left(\alpha^k\right)}\rho^{k}$
    \State $X^{k+1} = X^{k} + P^k\lambda^k\sigma^{k}$
    \State $\tilde{R}^{k+1} = \bar{R}^{k} - Q^k\lambda^k$
    \State $\bar{R}^{k+1},\gamma^{k+1} = \normalizer[\tilde{R}^{k+1}]{\SubA}$
    \State $\sigma^{k+1} = \gamma^{k+1}\sigma^{k}$
    \State break if $\|\sigma^{k+1}\| < \tol$
    \State $Z^{k+1} = \inverse{\Prec}\bar{R}^{k+1}$
    \State $\rho^{k+1} = \sproduct[\SubA]{Z^{k+1}}{\bar{R}^{k+1}}$
    \State $\beta^{k} =
    \inverse{\left(\rho^{k}\right)}\transpose{\gamma^{k+1}} \rho^{k+1}$
    \State $P^{k+1} = Z^{k+1} + P^k\beta^{k}$
    \EndFor
  \end{algorithmic}
\end{algorithm}
Note that it is not necessary to compute the real residual for
checking the convergence criterion, as we have
\begin{align}
  \|R^{k}_{i}\| = \|Q^{k}\sigma^{k}_{i}\| = \|\sigma^{k}_{i}\|,
\end{align}
where $\sigma^{k}_{i}$ denotes the $i$th column of $\sigma^{k}$.

As the orthonormalization is expensive, it makes sense to skip it in
iterations in which it is not necessary.
To specify a criterion for the adaptive orthonormalization, we define
the diagonally scaled condition number.
\pagebreak[4]
\begin{defn}[Diagonally scaled condition number]
  For a symmetric matrix ${\alpha \in \bbR^{k\times k}}$, we define the
  diagonally scaled condition number $\kappa_{D}(\alpha)$ as
  \begin{align}
    \kappa_{D}(\alpha) = \kappa(\delta^{-\frac12}\alpha\delta^{-\frac12}),
  \end{align}
  where $\delta = \diag{\alpha}$ is the diagonal of $\alpha$ and
  $\kappa$ denotes the condition number.
\end{defn} In contrast to the condition number, the diagonally scaled condition
number equals $1$ for diagonal matrices.
This is desirable in particular in the parallel case.
As it is only a parallel version of the scalar CG, no re-orthonormalization is
necessary.
Using the usual condition number, in this case, could deliver high numbers if the
columns are scaled differently and would lead to superfluous
re-orthonormalizations.
We use the diagonally scaled condition number of $\alpha^{k}$ for an indicator
of the numerical rank deficiency of the residual.
To check this, we evaluate
\begin{align}
  \label{eq:reortho_criterion}
  \eta \kappa_{D}(\alpha^{k}) > \sqrt{\machineeps},
\end{align}
where $\eta$ is a tuning parameter and $\machineeps$ is the machine
precision of the used numerical type.
Algorithm \ref{alg:bcgaro} shows the resulting algorithm.
\begin{algorithm}[tbp]
    \caption{BCG Method with Adaptive Residual Re-Orthonormalization}
  \label{alg:bcgaro}
  \begin{algorithmic}
    \State $R^0 = B-AX^0$
    \If{$\eta > 0$}
    \State $\bar{R}^{0},\sigma^{0} = \normalizer[R^{0}]{\SubA}$
    \Else
    \State $\bar{R}^{0} = R^{0}$
    \State $\sigma^{0} = \Identity[\SubA]$
    \EndIf
    \State $P^0 = \inverse{\Prec}\bar{R}^0$
    \State $\rho^{0} = \sproduct[\SubA]{P^{0}}{\bar{R}^{0}}$
    \For{$k = 0,\ldots$ until convergence}
    \State $Q^k = AP^k$
    \State $\alpha^k = \sproduct[\SubA]{P^k}{Q^k}$
    \State $\lambda^k = \inverse{\left(\alpha^k\right)}\rho^{k}$
    \State $X^{k+1} = X^{k} + P^k\lambda^k\sigma^{k}$
    \State $\tilde{R}^{k+1} = \bar{R}^{k} - Q^k\lambda^k$
    \If{$\eta \kappa_{D}(\alpha^{k}) > \sqrt{\machineeps}$}
    \State $\bar{R}^{k+1},\gamma^{k+1} = \normalizer[\tilde{R}^{k+1}]{\SubA}$
    \State $\sigma^{k+1} = \gamma^{k+1}\sigma^{k}$
    \State break if $\|\sigma^{k+1}\| \leq \tol$
    \Else
    \State $\bar{R}^{k+1} = \tilde{R}^{k+1}$
    \State $\gamma^{k+1} = \Identity[\SubA]$
    \State $\sigma^{k+1} = \sigma^{k}$
    \State break if $\|\bar{R}^{k+1}\sigma^{k+1}\| \leq \tol$
    \EndIf
    \State $Z^{k+1} = \inverse{\Prec}\bar{R}^{k+1}$
    \State $\rho^{k+1} = \sproduct[\SubA]{Z^{k+1}}{\bar{R}^{k+1}}$
    \State $\beta^{k} =
    \inverse{\left(\rho^{k}\right)}\transpose{\gamma^{k+1}} \rho^{k+1}$
    \State $P^{k+1} = Z^{k+1} + P^k\beta^{k}$
    \EndFor
  \end{algorithmic}
\end{algorithm}

This approach assumes that the diagonally scaled condition number
increases continuously with the iterations.
This is not clear though.
However, numerical tests show that this approach works quite well.
Nevertheless, there is some mathematical background missing.
An alternative approach would be to roll-back one iteration if the
re-orthonormalization criterion \eqref{eq:reortho_criterion} is
satisfied.
This would increase the memory requirements by one block vector and leads to
some overhead as some computations must be redone.

Theoretically, the normalization could add artificial directions to the
orthogonal residual, if the residual is rank-deficient.
In the case where this direction is already contained in the Krylov space, this
has no effect, because it is orthogonal to the residual.
Otherwise, this would accelerate the convergence as it enhances the Krylov
subspace.
Whether these new directions could be chosen more cleverly and the comparison
with deflation strategies is an objective of
future work.

\section{Numerical Experiments}
As a first test series, we executed test runs to approve the
convergence theory developed in Section \ref{sec:cg_convergence}.
\begin{figure}[t]
  \centering
  \resizebox{\textwidth}{!}{
      \input{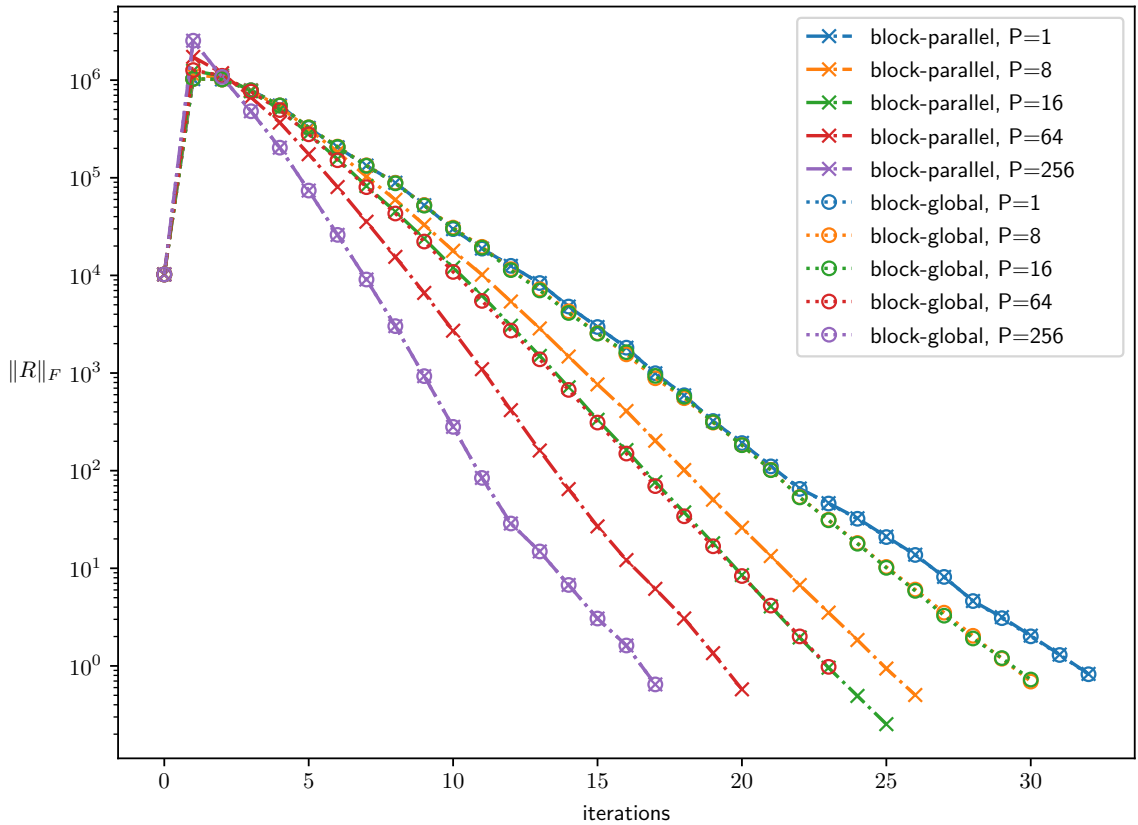}
    }
  \caption[Frobenius norm of the residual vs. iterations of BCG
    methods.]{Frobenius norm of the residual vs. iterations of BCG
    methods. Different choices for the *-subalgebra $\SubA$ are taken
    into account.
    Dashed lines with crosses decode the block-parallel method.
    Dotted lines with circles decode the block-global method.
    The colors decode the parameter $p$.
  }
  \label{fig:cg_convergence}
\end{figure}
In Figure \ref{fig:cg_convergence}, the convergence behavior for
different block sizes $p$ and the block-global and block-parallel case is shown.
As operator we used the \texttt{thermal2} matrix from the
SuiteSparse Matrix Collection \cite[]{lawry2002comb}.
This matrix provides a realistic size for tests on one machine and we
observed that the AMG preconditioner from the \dune framework worked
sufficiently good.
For the tests we used $s=256$ randomly generated right-hand sides.

The results confirm perfectly the theoretical expectations.
In particular, we obtain a faster convergence (per iteration) for larger block
size $p$.
Furthermore, we see that the convergence rates in the block-global method
$\SubA[BG]^{p}$ with $p$ up to $8$ are similar to the convergence rate of the
parallel method $\SubA[P]$.
This fits to the outcome of Lemma \ref{lem:global_cg_convergence}, as it
predicted that the convergence rate depends on the $\left\lceil \frac{p}{q}
\right\rceil$ smallest eigenvalue.
For $p\leq 8$, we have $\left\lceil \frac{p}{q} \right\rceil = 1$.
For the same reason, we see a connection of the convergence rate of the
block-global method with $p=64$, i.e.\ $\SubA[BG]^{64}$ and block-parallel method
with $p=16$, i.e.\ $\SubA[BP]^{16}$.

Note that the increase in the first iteration is due to the fact that
the Frobenius norm of the residual is plotted, instead of the
energy Frobenius norm of the error, which converges monotonically.

Further, we tested the methods behavior depending on the
re-orthonormalization parameter $\eta$.
\begin{table}
  \centering
    \caption[Iteration counts and numbers of residual re-orthonormalization for
single and double precision.]{Iteration counts and numbers of residual
re-orthonormalization for single and double precision. Blank cells indicate that
the method did not converged within 1000 iterations.
The parallel case is added for comparison. In all other rows the block-parallel
method with $p=64$ was used.}
  \begin{tabular}{rcccc}\toprule
    precision&\multicolumn{2}{c}{single}&\multicolumn{2}{c}{double}\\
    \cmidrule(lr){2-3} \cmidrule(lr){4-5}
    $\eta$ & iterations & ortho. & iterations & ortho.\\
    \midrule
    parallel & $868$ & $0$ & $514$ & $0$\\
    $0$ & & & & \\
    $0.1$ & $434$ & $10$ & $146$ & $3$\\
    $1$ & $344$ & $14$ & $137$ & $5$\\
    $10$ & $286$ & $30$ & $137$ & $5$\\
    $100$ & $229$ & $55$ & $104$ & $3$\\
    $1000$ & $231$ & $230$ & $88$ & $3$\\
    $\infty$ & $237$ & $238$ & $74$ & $75$\\
    \bottomrule
  \end{tabular}
  \label{tab:it_count_reortho}
\end{table} Table \ref{tab:it_count_reortho} shows iteration and
orthonormalization counts for different choices of $\eta$ and different
floating-point precision.
The \texttt{HB/1138\_bus} matrix from the SuiteSparse Matrix Collection
\cite[]{timothy2011university} with a SSOR preconditioner was used with $256$ randomly
generated right-hand sides.
The parallel method ($\SubA[P]$) is shown for comparison.
In all other rows the block-parallel with $p=64$ ($\SubA[BP]^{64}$) is used.
It can be seen that the number of iterations decreases if we use a higher
re-orthonormalization parameter $\eta$.
It can also be seen that a higher $\eta$ not necessarily leads to more
re-orthonormalizations.
Hence, the moment of the re-orthonormalization seems to be important, a fact
that can also be observed in the next experiment.

\begin{figure}[t]
  \centering
  \resizebox{\textwidth}{!}{
      \input{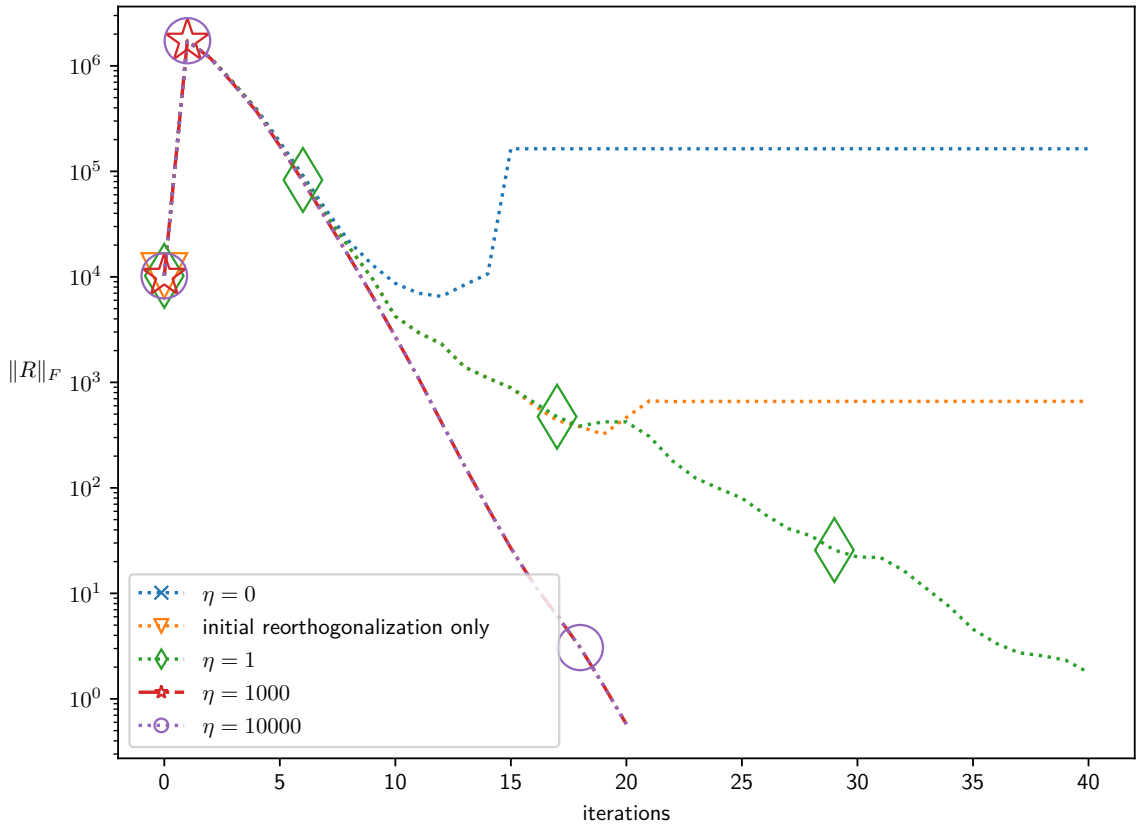}
    }
  \caption[Convergence of the BCG method for different re-orthonormalization
    parameters.]{Convergence of the BCG method for different re-orthonormalization
    parameters.
  The symbols mark the iteration in which a re-orthonormalization
  happened.
  Colors decode the different values of the re-orthonormalization parameter $\eta$.}
  \label{fig:cg_reortho}
\end{figure}
Figure \ref{fig:cg_reortho} shows the convergence history of the
BCG algorithm in the same setting used in Figure \ref{fig:cg_convergence},
for different re-orthonormalization parameters $\eta$.
We used the block-parallel setting with $p=64$.
We see that the re-orthonormalization is necessary to achieve
convergence and the re-orthonormalization must be reapplied during the
iteration - it is not sufficient to orthonormalize the initial
residual.
Furthermore, we see that choosing $\eta=1$ yields a converging method,
but does not ensure optimal convergence rates.
Only in cases where in iteration $0$ and $1$ a orthonormalization was
applied the convergence rate was optimal.
The additional re-orthonormalization in iteration $18$ for
$\eta=\num{10000}$ seems to be superfluous.

\begin{table}
  \centering
    \caption[Iteration counts, runtime and number of re-orthonormalizations for the
solution of several matrices from MatrixMarket with different $p$.]{Iteration
counts, runtime and number of re-orthonormalizations for the solution of several
matrices from MatrixMarket with different $p$.
The fastest time per row is marked with a green background.
Blank cells indicate that the method did not converge within $1000$ iterations.}
  \begin{tabular}{@{}p{16ex}*{9}{c}@{}}\toprule
    & \multicolumn{3}{c}{$p=1$} & \multicolumn{3}{c}{$p=32$} &
                                                               \multicolumn{3}{c}{$p=256$}\\
    \cmidrule(lr){2-4} \cmidrule(lr){5-7} \cmidrule(lr){8-10}
    Matrix & \#it & \#ro & $t$ & \#it & \#ro & $t$ & \#it &
                                                                  \#ro & $t$\\
    \midrule
     bcsstk14 & $251$ &    $0$ &                     $3.18$ &  $16$ &   $13$ & \cellcolor{green!30}$1.15$ &  $7$ &    $5$ &                       $1.79$ \\
 bcsstk15 & $573$ &    $0$ &                    $12.24$ &  $23$ &   $12$ & \cellcolor{green!30}$1.97$ & $10$ &    $5$ &                       $2.41$ \\
 bcsstk16 &  $26$ &    $0$ & \cellcolor{green!30}$1.02$ &   $9$ &    $1$ &                     $1.20$ &  $6$ &    $2$ &                       $1.64$ \\
 bcsstk17 &       &        &                            &  $73$ &   $71$ &                    $14.54$ & $16$ &   $15$ &   \cellcolor{green!30}$8.09$ \\
 bcsstk18 & $419$ &    $0$ &                    $21.89$ &  $49$ &    $5$ &                     $6.59$ & $15$ &    $4$ &   \cellcolor{green!30}$6.23$ \\
 s1rmq4m1 &  $85$ &    $0$ &                     $3.36$ &  $16$ &    $2$ & \cellcolor{green!30}$1.90$ &  $9$ &    $2$ &                       $2.81$ \\
 s1rmt3m1 & $166$ &    $0$ &                     $5.75$ &  $24$ &    $3$ & \cellcolor{green!30}$3.09$ & $12$ &    $2$ &                       $3.36$ \\
 s2rmq4m1 & $107$ &    $0$ &                     $4.24$ &  $17$ &    $3$ & \cellcolor{green!30}$2.02$ & $11$ &    $2$ &                       $3.33$ \\
 s2rmt3m1 & $226$ &    $0$ &                     $7.78$ &  $30$ &    $3$ & \cellcolor{green!30}$3.86$ & $14$ &    $2$ &                       $4.02$ \\
 s3dkq4m2 &       &        &                            & $217$ &    $8$ &                   $190.09$ & $65$ &    $9$ & \cellcolor{green!30}$138.10$ \\
 s3rmq4m1 & $197$ &    $0$ &                     $7.82$ &  $21$ &    $3$ & \cellcolor{green!30}$2.84$ & $12$ &    $3$ &                       $3.71$ \\
 s3rmt3m1 & $478$ &    $0$ &                    $16.80$ &  $33$ &    $5$ & \cellcolor{green!30}$3.03$ & $17$ &    $5$ &                       $5.13$ \\
 s3rmt3m3 & $442$ &    $0$ &                    $14.40$ &  $33$ &    $5$ & \cellcolor{green!30}$3.69$ & $15$ &    $4$ &                       $3.98$ \\
\bottomrule
  \end{tabular}
  \label{tab:smalltests}
\end{table}
Table \ref{tab:smalltests} shows the iteration counts, the number of
re-orthonormalizations and the runtime of the BCG method for the block-parallel
method and different value of $p$ for several symmetric positive definite
matrices of the SuiteSparse Matrix Collection \cite[]{lawry2002comb}.
We used an incomplete Cholesky preconditioner and the
re-orthonormalization parameter $\eta=\num{10000}$.
All systems are solved for $s=256$ randomly generated right-hand sides.
The fastest runtime per matrix is marked with a green background.
Missing numbers mark that no convergence was achieved within $\num{1000}$ iterations.
We aimed for a reduction of the residual in every column
by a factor of $\num{e-4}$.

The result shows that the number of iterations can be reduced
drastically by using a higher $p$.
We suppose that this is due to the weak preconditioner, that mainly smooths
the larger eigenvalues.
For example for the \texttt{bcsstk15} matrix, the number of iterations
was reduced by a factor of $\sim 25$ by using $p=32$ compared to the parallel case.
In some cases the block Krylov methods help to achieve convergence at
all, e.g.\ in the \texttt{bcsstk17} case, where the parallel case fails to
converge within $\num{1000}$ iterations, but the block methods converge within
$73$ iterations for the $p=32$ case.

For the \texttt{bcsstk16} matrix the parallel method is fastest,
although it need  the most iterations.
This is due to the re-orthonormalization costs, as the other building
blocks are similar expensive for the $p=32$ case.
An improvement of the re-orthonormalization criteria is one of the goals of
future work.


\chapter{Block GMRes Method}
\label{chap:blockgmres}
In the previous chapter we looked at the block Conjugate Gradients
method, which is the method of choice for symmetric positive definite problems.
We now discuss the block GMRes (BGMRes)
method~\cite[]{vital1990etude}, which is a block version of the GMRes
method by \citeauthor[]{saad1986gmres}~\cite[]{saad1986gmres}.
In contrast to the BCG method, it does not have any
requirements on the operator.
As a downside, the GMRes method can not make use of a short
recurrence, i.e.\ the memory and arithmetically costs per iteration
increase with every iteration.
Nevertheless, it is one of the most important Krylov methods in
practice.
Note that for symmetric indefinite problems there are also block
versions of the MinRes~\cite[]{soodhalter2015block} and conjugate
residual~\cite[]{zhang2013novel} method, which are not subject of this
work.
Good introductions into the classical block GMRes method can be found in the
monographs of \citeauthor[]{gutknecht2007block}~\cite[]{gutknecht2007block} and
\citeauthor[]{saad2003iterative}~\cite{saad2003iterative}.

We will formulate the BGMRes method based on the block Krylov
framework.
Recently, \citeauthor[]{kubinova2020admissible}~\cite[]{kubinova2020admissible}
presented a paper that also discusses the BGMRes method in this
framework and contains some results about the convergence of the
method. In addition, a generalization of the Givens rotations that are
used in the non-block case to triangulate the Hessenberg matrix is described.
We pick up this generalization in the first section and formulate the
BGMRes method.
We present some simple convergence results in the second section and
refer the reader to \cite[]{kubinova2020admissible} for a more
elaborate discussion.
Finally, we present some numerical experiments that give further insights into
the convergence behavior of the BGMRes method for different *-subalgebras.

\section{Formulation of the Block GMRes Method}
The BGMRes method is based on the block Arnoldi
process~\cite[]{arnoldi1951the,ruhe1979implementation}, which computes
an orthonormal basis of the Krylov space and was originally invented
to compute the eigenvalues of an operator.
It is based on the block Gram-Schmidt orthogonalization process,
see Algorithm \ref{alg:block-arnoldi}.
\begin{algorithm}[t]
  \caption{Block Arnoldi Method}
  \label{alg:block-arnoldi}
  \begin{algorithmic}
    \State Let $V^{0}\in\bbR^{\dimA\times s}$ with $\sproduct[\SubA]{V^{0}}{V^{0}}=0$ be given.
    \For{$k=1,\ldots k_{\max}$}
    \State $V^{k} = AV^{k-1}$
    \For{$j=0,\ldots,k-1$}
    \State $\eta_{j,k-1} = -\sproduct[\SubA]{V^{j}}{V^{k}}$
    \State $V^{k} \gets V^{k} + V^{j}\eta_{j,k-1}$
    \EndFor
    \State $V^{k},\eta_{k,k-1} \gets \normalizer[V^{k}]{\SubA}$
    \EndFor
\end{algorithmic}
\end{algorithm}

The coefficients $\eta$ build a block matrix
$\mathcal{H}\in\SubA^{k+1\times k}$, that has block Hessenberg form,
i.e.\ all blocks below the first off-diagonal under the diagonal are zero.
The resulting basis $\mathcal{V} = \left[V^{0},\ldots,V^{k_{\max}}\right]$
satisfies the so called \textit{block Arnoldi relation}
\begin{align}
  A\tilde{\mathcal{V}} = \mathcal{V}\mathcal{H},
\end{align}
where $\tilde{\mathcal{V}} = \left[V^{0},\ldots,V^{k_{\max}-1}\right]$.

As the normalizer is also defined for rank-deficient block vectors, we do not
get a breakdown in the case where $V^{k}$ is rank-deficient after the
Gram-Schmidt orthonormalization.
The normalization process adds additional directions to the Krylov space in this
case.
Theoretically the orthogonalization must be repeated to orthogonalize the
additional directions to the previous block vectors.
However, numerical experiments show, that this is not necessary, even if the
rigorous analysis of this effect is still missing.


Algorithm \ref{alg:block-arnoldi} uses the modified Gram-Schmidt procedure,
meaning the computation of the block inner products and the vector updates are
interleaved.
The modified Gram-Schmidt procedure is more stable than the classical
Gram-Schmidt procedure which computes all block inner products in advance.
However, as the block inner product computes multiple inner products
simultaneously the stability could be affected.
This could be mitigated by either using the ``real'' modified Gram-Schmidt that
considers the columns of the block vectors individually or by doing a
re-orthogonalization like presented by
\citeauthor[]{bjorck1994numerics}~\cite[]{bjorck1994numerics}.
\citeauthor[]{buhr2014numerically}~\cite[Algorithm 1]{buhr2014numerically}
presented an adaptive re-orthogonalization strategy for the Gram-Schmidt
procedure, that could be applied to decide adaptively whether a
re-orthonormalization is necessary.

The goal of the BGMRes method is to compute an update
$U^{k}\in\krylov[\SubA]{A}{R^{0}}{k}$
for the initial guess $X^{0}$ that solves the minimization problem
\begin{align}
  \label{eq:gmres_mini}
  U^{k} =
  \argmin_{Y\in\krylov[\SubA]{A}{R^{0}}{k}}
  \|B - AX^{0} - AY\|_{F}.
\end{align}
In other words, the Frobenius norm of the block residual is minimized.
With help of the Arnoldi basis the update reads
\begin{align}
  U^{k} = \tilde{\mathcal{V}}\bm{\zeta}^{k},
\end{align} where $\bm{\zeta}^{k}\in\SubA^{k}$ denote the coefficients of
$U^{k}$ in the basis $\tilde{\mathcal{V}}$.
Using the block Arnoldi relation, the orthonormality of $\mathcal{V}$ and a
block QR decomposition $\mathcal{QR}=\mathcal{H}$, with $\mathcal{Q}\in
\SubA^{k+1\times k}$ and $\mathcal{R}\in\SubA^{k\times k}$ we rewrite the
minimization problem \eqref{eq:gmres_mini} as
\begin{align}
  \|B - AX^{0} - AU^{k}\|_{F} &= \|R^{0} - A\tilde{\mathcal{V}}\bm{\zeta}^{k}\|_{F}\\
                              &= \|R^{0} -
                                \mathcal{V}\mathcal{H}\bm{\zeta}^{k}\|_{F}\\
                              &= \|\transpose{\mathcal{V}}R^{0} -
                                \mathcal{H}\bm{\zeta}^{k}\|_{F}\\
                              &=
                                \|\transpose{\mathcal{Q}}\transpose{\mathcal{V}}R^{0}
                                - \mathcal{R}\bm{\zeta}^{k}\|_{F}.
\end{align}
This minimization problem can then be solved for $\bm{\zeta}^{k}$ by
block-wise backward-substitution.

In the non-block case the QR decomposition of $\mathcal{H}$ is computed
using Given rotations to eliminate the lower off-diagonal entries.
This can be generalized for our block Krylov framework.
For that, we compute a full QR decomposition of the diagonal and lower
off-diagonal entry, starting with
\begin{align}
  Q_{0}\begin{pmatrix}
    \rho_{0}\\
    0
  \end{pmatrix}
  &= \begin{pmatrix}
    \eta_{0,0}\\
    \eta_{1,0}
  \end{pmatrix}
  & \transpose{Q_{0}}Q_{0} &=
                             \begin{pmatrix}
                               \Identity&0\\
                               0&\Identity
                             \end{pmatrix}
                           & Q_{0}\in\SubA^{2\times 2}, \rho_{0}\in\SubA.
\end{align}
The lower off-diagonal element can then be eliminated by
\begin{align}
  \begin{pmatrix}
    \multicolumn{2}{c}{\multirow{2}{*}{\scalebox{1.5}{$\transpose{Q_{0}}$}}}\\
    \\
    & & \Identity\\
    & & & \ddots\\
    & & & & \Identity
  \end{pmatrix}
          \begin{pmatrix}
            \eta_{0,0} & & \multicolumn{2}{c}{\multirow{2}{*}{*}}\\
            \eta_{1,0} & \eta_{1,1} & &\\
            & \eta_{2,1} & \ddots\\
            & & \ddots
          \end{pmatrix}
      &=
        \begin{pmatrix}
          \rho_{0} & &\multicolumn{2}{c}{\multirow{2}{*}{\scalebox{1.5}{*}}}\\
          0 & & &\\
          & \eta_{2,1} & \ddots\\
          & & \ddots
        \end{pmatrix}.
\end{align}
The star indicates non-zero entries in the upper triangle.
This procedure is repeated to eliminate the other lower
off-diagonal entries.
The $Q$-factor of the QR decomposition of $\mathcal{H}$ is then build
by concatenating all the $Q$-factors of the smaller QR decompositions
\begin{align}
  \mathcal{Q} = \begin{pmatrix}
    \multicolumn{2}{c}{\multirow{2}{*}{\scalebox{1.5}{$\transpose{Q_{0}}$}}}\\
    \\
    & & \Identity\\
    & & & \Identity\\
    & & & & \ddots\\
    & & & & & \Identity
  \end{pmatrix}
            \begin{pmatrix}
              \Identity\\
              &\multicolumn{2}{c}{\multirow{2}{*}{\scalebox{1.5}{$\transpose{Q_{1}}$}}}\\
    \\
    && & \Identity\\
    && & & \ddots\\
    && & & & \Identity
  \end{pmatrix}
             \cdots
\end{align}
In the algorithm the transformation of $R^{0}$ and the
QR decomposition of $\mathcal{H}$ is performed on-the-fly.
The vector
\begin{align}
  \bm{\sigma} = \begin{pmatrix}
    \sigma^{0}\\
    \vdots\\
    \sigma^{k_{\max}}
  \end{pmatrix}
  = \transpose{\mathcal{Q}}\transpose{\mathcal{V}}R^{0}
\end{align}
is updated during the iteration by
\begin{align}
  \sigma^{0} = \normalizer[M^{-1}R^{0}]{\SubA},\\
  \begin{pmatrix}
    \sigma^{k}\\
    \sigma^{k+1}
  \end{pmatrix}
  \gets \transpose{Q_{i}}\begin{pmatrix}
    \sigma^{k}\\
    0
  \end{pmatrix}.
\end{align}
The Frobenius norm of $\sigma^{k+1}$ can be used to determine the
residual in the $k$-iteration, as
\begin{align}
  \|R^{k}\|_{F}
  &= \|R^{0}-A\tilde{\mathcal{V}}\bm{\zeta}^{k}\|_{F}\\
  &= \|R^{0}-\mathcal{V}\mathcal{QR}\bm{\zeta}^{k}\|_{F}\\
  &= \|\transpose{\mathcal{Q}}\transpose{\mathcal{V}}R^{0} -
    \mathcal{R}\bm{\zeta}^{k}\|_{F}\\
  &= \|\bm{\sigma} - R\bm{\zeta}^{k}\|_{F} = \|\sigma^{k+1}\|_{F}.
\end{align}

Algorithm \ref{alg:bgmres} shows the BGMRes algorithm.
Preconditioning can be easily implemented by adapting lines
\ref{algline:gmres_v0} and \ref{algline:gmres_vk}.

\begin{algorithm}[htb!]
  \caption{Block GMRes Method}
  \label{alg:bgmres}
  \begin{algorithmic}[1]
    \State $R^{0} = B-AX^{0}$
    \State $V^{0} \sigma^{0} = \normalizer[R^{0}]{\SubA}$ \label{algline:gmres_v0}
    \For{$k=0,\ldots,k_{\max}-1$}
    \State $V^{k+1} = AV^{k}$  \label{algline:gmres_vk}
    \For{$j=0,\ldots,k$}
    \State $\eta_{j,k} = -\sproduct[\SubA]{V^{j}}{V^{k+1}}$
    \State $V^{k+1} \gets V^{k+1} + V^{j}\eta_{j,k}$
    \EndFor
    \State $V^{k+1},\gamma \gets \normalizer[V^{k+1}]{\SubA}$
    \For{$j=0,\ldots,k-1$}
    \State $\begin{pmatrix}
      \eta_{j,k}\\
      \eta_{j+1,k}
    \end{pmatrix} \gets \transpose{Q_{j}}\begin{pmatrix}
      \eta_{j,k}\\
      \eta_{j+1,k}
    \end{pmatrix}$
    \EndFor
    \State $Q_{k}\begin{pmatrix}\eta_{k,k}\\0\end{pmatrix}
    \gets \begin{pmatrix}\eta_{k,k}\\ \gamma\end{pmatrix}$
    \Comment{Compute QR decompostion}
    \State $\begin{pmatrix}\sigma^{k}\\ \sigma^{k+1}\end{pmatrix}
    \gets \transpose{Q_{k}}\begin{pmatrix}\sigma^{k}\\ 0\end{pmatrix}$
    \State break if $\|\sigma^{k+1}\| \leq \tol$
    \EndFor
    \For{$l=k_{\max}-1,\ldots,0$}
    \Comment{back-substitution}
    \State $\sigma^{l} \gets \sigma^{l} - \sum_{j=l+1}^{k_{\max}-1} \eta_{l,j}\sigma^{j}$
    \State $\sigma^{l} \gets \inverse{\left(\eta_{l,l}\right)}\sigma^{l}$
    \EndFor
    \State $X^{k_{\max}} = X^{0} + \sum_{j=0}^{k_{\max}-1} V^{j}\sigma^{j}$
    \If{$\|\sigma^{k_{\max}}\| > \tol$}
    \State restart with $X^{0} \gets X^{k_{\max}}$
    \EndIf
  \end{algorithmic}
\end{algorithm}

\section{Convergence}
\label{sec:gmres_convergence}
For the GMRes method we do not have a general theoretical statement
about the convergence rate, like for the CG case.
Rather \citeauthor[]{greenbaum1996any}~\cite[]{greenbaum1996any}
showed that any non-increasing convergence curve for the GMRes method
is possible.
This result was recently generalized by
\citeauthor[]{kubinova2020admissible}~\cite[]{kubinova2020admissible}
for the BGMRes method.
A convergence theory for the BGMRes method for special classes of
operators was presented by
\citeauthor[]{simoncini1996convergence}~\cite[]{simoncini1996convergence}.

As the BGMRes method minimizes the Frobenius norm of the residual
in the block Krylov space we have
\begin{align}
  \|R^{k+1}\|_{F} \leq \|R^{k}\|_{F}.
\end{align}
This ensures a monotonic convergence of the Frobenius norm of the residual,
but it cannot be ensured that the residual actually decreases, see
Example \ref{ex:nonconverging_gmres}.

To deduce better estimations a-priori knowledge about the operator is
necessary.
We use again the polynomial representation to formulate an abstract
statement about the convergence rate, similar to the convergence proof
of the BCG method.
This statement could be used to deduce concrete estimations if
further assumptions on the operator are made.

\pagebreak[4]
\begin{lem}[Abstract convergence of BGMRes method]
  Let $\SubA$ be a *-subalgebra of $\bbR^{s\times s}$.
  For the residual $R^{k}$ of the $k$th step in the BGMRes
  method the following estimation holds
  \begin{align}
    \|R^{k}\|_{F}
    \leq \inf_{\mathcal{Q}^{k}} \|\mathcal{Q}^{k}(A)\circ R^{0}\|_{F}
    \leq \inf_{\mathcal{Q}^{k}} \|\mathcal{Q}^{k}(A)\|\|R^{0}\|_{F}.
  \end{align}
  The infimum is taken over all $\SubA$-valued polynomials
  $\mathcal{Q}^{k}\in\bb{P}_{\SubA}^{k}$ with absolute coefficient
  $\Identity$.
  The norm $\|\mathcal{Q}^{k}(A)\|$ denotes the
  operator norm of $\mathcal{Q}^{k}(A)$ in the space
  $\mathcal{L}(\bbR^{s\dimA}, \bbR^{s\dimA})$.
\end{lem}
\begin{proof}
  The BGMRes method computes the best approximation in the space
  ${X^{0} + \krylov[\SubA]{A}{R^{0}}{k}}$.
  As we can represent the elements in the Krylov space with
  $\SubA$-valued polynomials we obtain
  \begin{align}
    \label{eq:abstract_bgmres_convergence}
    R^{k} &= B - AX^{k}\\
          &= B - AX^{0} - AU^{k}\\
          &= R^{0} - A\mathcal{P}^{k-1}(A)\circ R^{0}\\
          &= \left(\Identity -
            A\mathcal{P}^{k-1}(A)\right)\circ
    R^{0}\\
          &= \mathcal{Q}^{k}(A)\circ R^{0},
  \end{align}
  where $\mathcal{Q}^{k}(\polyX) = \Identity - \polyX\mathcal{P}^{k}(\polyX)$.
  Applying the Frobenius norm and taking the infimum completes
  the proof, as the Frobenius norm on $\bbR^{\dimA\times s}$ and the
  Euclidean norm on the space $\bbR^{\dimA s}$ coincide.
\end{proof}

The next lemma gives an example how this estimation could be used to create more
concrete estimations.
This is a generalization for the block Krylov framework of Theorem 3.1 in the
work of \citeauthor[]{simoncini1996convergence}~\cite[Theorem
3.1]{simoncini1996convergence}.

\begin{lem}
If the operator $A$ is diagonalizable
\begin{align}
  A &= V\Lambda\inverse{V},
  & \text{with} &&\Lambda &= \diag{\lambda_{1},\ldots,\lambda_{n}}
\end{align}
estimation \eqref{eq:abstract_bgmres_convergence} can be precised as
\begin{align}
  \|R^{k}\|_{F} &\leq \kappa(V) \inf_{c_{1}\,\ldots,c_{k}\in\SubA}
              \max_{i=1}^{\dimA}
              \left\|\sum_{j=0}^{k}\lambda_{i}^{j}c_{j}\right\| \|R^{0}\|_{F},
\end{align}
where $c_{0} = \Identity$.
\end{lem}
\begin{proof}
  Let $c_{1},\ldots,c_{k}\in\SubA$ denote the coefficients of the
  $\SubA$-valued polynomial $\mathcal{Q}^{k}$.
  Then we write
  \begin{align}
    \left\|\mathcal{Q}^{k}(A)\right\|
    &= \left\|\sum_{j=0}^{k} A^j \otimes c_{i}\right\|\\
    &= \left\|\sum_{j=0}^{k} V\Lambda^{j}\inverse{V} \otimes c_{i}\right\|\\
    &= \left\|\left(V\otimes \Identity\right)\left(\sum_{j=0}^{k}
      \Lambda^{j} \otimes c_{i}\right)\left(\inverse{V}\otimes
      \Identity\right)\right\|\\
    &\leq \kappa(V)  \max_{i=1}^{\dimA} \left\|\sum_{j=0}^{k}\lambda_{i}^{j}c_{j}\right\|
  \end{align}
\end{proof}

The challenge is to choose good coefficients
$c_{1},\ldots,c_{k}\in\SubA$.
If all $\lambda^{i}$ are positive, then probably the scaled Chebyshev
polynomials would yield an estimation similar to the CG case.
See the recent paper of
\citeauthor[]{kubinova2020admissible}\cite[]{kubinova2020admissible} for a
detailed discussion.
Further convergence results of the classical BGMRes method can be found in the
paper of \citeauthor[]{simoncini1996convergence}~\cite[]{simoncini1996convergence}.

\section{Numerical Experiments}
As the theoretical convergence results are still quite vague yet, we rely on
numerical test to get an impression of the convergence behavior of the method
with respect to the different *-subalgebras.
As the BGMRes method minimizes the Frobenius norm of the residual we know that
for the same $p$ the block-parallel method converges faster than the block-global
method.
In both cases the larger the $p$ the better the convergence rate
(per iteration), cf. Lemma \ref{lem:algebraembedding} and Lemma
\ref{lem:krylovembedding}.

\begin{figure}[tb]
  \centering
  \resizebox{\textwidth}{!}{
      \input{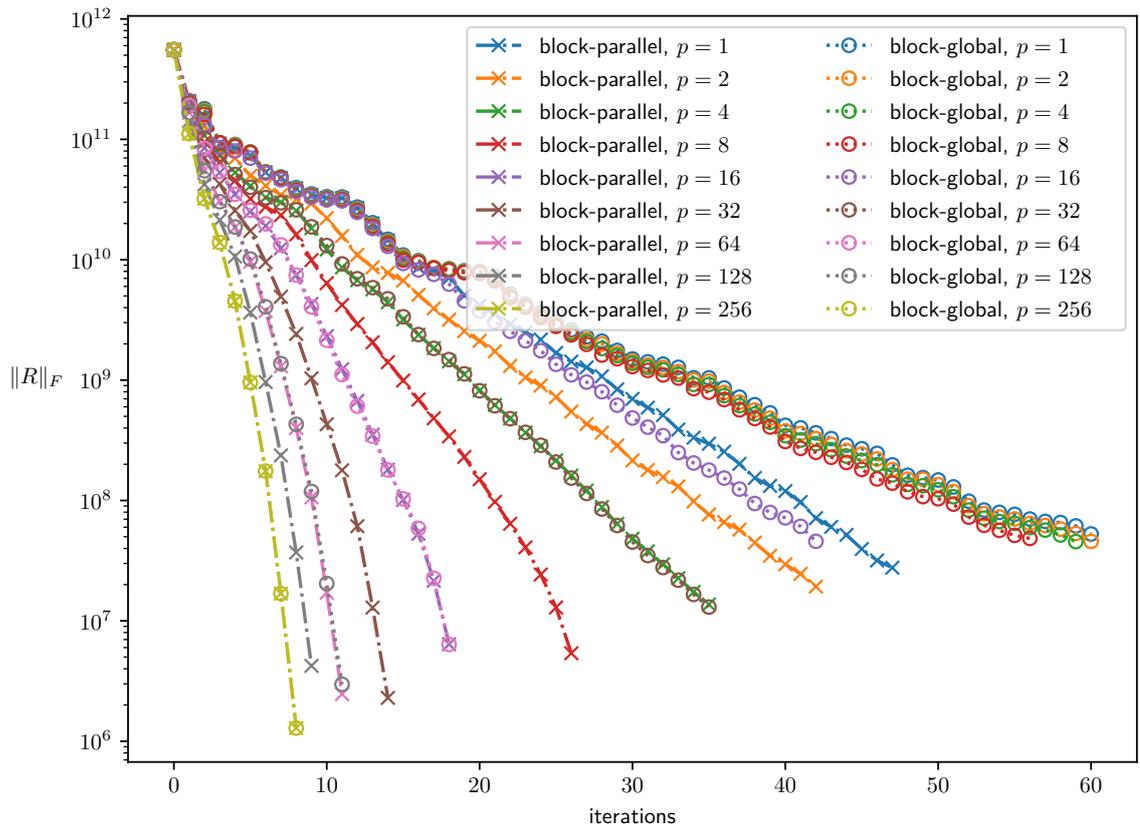}
    }
  \caption[Convergence of block GMRes method.]{Convergence of block GMRes method.
    Dashed lines with crosses decode the block-parallel methods.
    Dotted lines with circles denote the block-global method.
    Colors decode the blocking parameter $p$.}
  \label{fig:gmres_hybrid_vs_global}
\end{figure}
It is confirmed by the result presented in Figure
\ref{fig:gmres_hybrid_vs_global}.
It shows the convergence of the
BGMRes method for the
\texttt{Simon/raefsky3} matrix from the SuiteSparse Matrix
Collection~\cite[]{timothy2011university}.
The problem consists of $\num{21200}$ unknowns and originates from a computational
fluid dynamics problem.
We use an ILU(0) preconditioner and solve for $s=256$ randomly generated
right-hand sides until a reduction of the 2-norm of the residual for
every column by a factor of $10^{-4}$ is reached.

We see a relation of the convergence rates of the
$p$-block-global method and the $\frac{p}{q}$-block-parallel method,
like in the BCG case.
For example the convergence for $\SubA[BG]^{128}$ ($p=128$ block-global) and
$\SubA[BP]^{64}$ ($p=64$ block-parallel) is almost identical.
The same holds for $\SubA[BG]^{64}$ ($p=64$ block-global) and $\SubA[BP]^{16}$ ($p=16$ block-parallel).
That indicates that a similar result to Lemma
\ref{lem:global_cg_convergence} could also be possible for the BGMRes method.

Note that choosing a large restart parameter in the BGMRes method is crucial for
achieving good convergence.
Often the choice of that parameter is limited by the memory of the
machine.
This means the restart length directly competes with the number of
right-hand sides that can be used.
If this is an issue, probably the block BiCGStab method which is
considered in the next chapter is a better choice to solve the problem, as its
memory requirements are constant.


\chapter{Block BiCGStab Method}
\label{chap:blockbicgstab}
As a third block Krylov method we look at the Block BiCGStab (BBiCGStab)
method, a block version of the BiCGStab method presented by
\citeauthor[]{van1992bi}~\cite[]{van1992bi}.
This chapter is based on the paper by
\citeauthor[]{el2003block}~\cite[]{el2003block}.
They deduce the block BiCGStab method, i.e.\ the case $\SubA=\bbR^{s\times s}$.
We adapt this deduction for the block Krylov framework.

Like the GMRes method, the BiCGStab method was developed for
non-symmetric problems.
Unlike the GMRes method, it does not store a basis and is
therefore better suited for memory limited systems.
This advantage comes with the price, that no minimization property
is satisfied by the approximate solution.
Hence, it is difficult to develop theoretical convergence results.

As with the other methods, we reformulate the BBiCGStab method based on the
block Krylov framework in the first section.
In the second section we introduce a stabilization for the BBiCGStab method,
similar to that for the BCG method.
Up to the authors knowledge such a stabilization strategy was not presented before.
Finally, we present some numerical results in section
\ref{sec:blocknumericbicgstab}.

\section{Formulation of the Block BiCGStab Method}

Like the name suggests, the BBiCGStab method is based on the block BiCG (BBiCG)
method~\cite[]{oleary1980block}, but adds a stabilization step to
mitigate instabilities.
Another issue of the BBiCG method is that the transposed of the operator must be
applied to a block vector.
The BBiCG method actually solves an additional linear system with the transposed
operator and computes bases $(P^{i})_{i}$ and $(\tilde{P}^{i})_{i}$ of the
Krylov spaces $\krylov{\inverse{\Prec}A}{\inverse{\Prec}R^{0}}{k}$ and
$\krylov{\transinv{\Prec}\transpose{A}}{\transinv{\Prec}\tilde{R}^{0}}{k}$, for
some block vector $\tilde{R}^{0} \in \bbR^{\dimA\times s}$ and preconditioner $\Prec\in\bbR^{\dimA\times\dimA}$.
The residuals of the linear systems are projected onto the Krylov space of the
other system.
In the absence of rounding errors, this ensures that the method converges in at
least $\dimA$ iterations.
However, it is not the aim to proceed as many iterations, as even a direct solve
would be more efficient.
The BBiCG algorithm is shown in Algorithm \ref{alg:bbicg}.

\begin{algorithm}
  \caption{Block BiCG Method}
  \label{alg:bbicg}
  \begin{algorithmic}
    \State $R^{0} = B-AX^{0}$
    \State Choose $\tilde{R^{0}}$ arbitrarily with
    $\sproduct[\SubA]{\tilde{R}^{0}}{R^{0}}\neq 0$
    \State $P^{0} = \Prec^{-1}R^{0}, \tilde{P}^{0} =
    \transinv{\Prec}\tilde{R}^{0}$
    \State $\rho^{0} = \sproduct[\SubA]{\tilde{R}^{0}}{P^{0}}$
    \For{$k=0,\ldots$}
    \State $\alpha^{k} = \sproduct[\SubA]{\tilde{P}^{k}}{AP^{k}}$
    \State $X^{k+1} = X^{k} + P^{k}\inverse{(\alpha^{k})}\rho^{k}$
    \State $R^{k+1} = R^{k} -
    AP^{k}\inverse{(\alpha^{k})}\rho^{k}$
    \State break if $\|R^{k+1}\| \leq \tol$
    \State $\tilde{R}^{k+1} = \tilde{R}^{k} -
    \transpose{A}\tilde{P}^{k}\transinv{(\alpha^{k})}\transpose{\rho^{k}}$
    \State $Z^{k+1} = \Prec^{-1}R^{k+1}, \tilde{Z}^{k+1} =
    \transinv{\Prec}\tilde{R}^{k+1}$
    \State $\rho^{k+1} =
    \sproduct[\SubA]{\tilde{R}^{k+1}}{Z^{k+1}}$
    \State $P^{k+1} = Z^{k+1} +
    P^{k}\left(\rho^{k}\right)^{-1}\rho^{k+1}$
    \State $\tilde{P}^{k+1} = \tilde{Z}^{k+1} +
    \tilde{P}^{k}\transinv{\left(\rho^{k}\right)}\transpose{\rho^{k+1}}$
    \EndFor
  \end{algorithmic}
\end{algorithm}

\begin{lem}
  \label{lem:BCG_properties}
  The residual in the $k$th iteration is orthogonal to the Krylov space of the
  adjunct problem.
  \begin{align}
    \sproduct[\SubA]{V}{R^{k}} &= 0
    & \forall V\in\krylov[\SubA]{\transinv{\Prec}
      \transpose{A}}{\transinv{\Prec}\tilde{R}^{0}}{k}.
  \end{align}
\end{lem}
A proof can be found in the paper by \citeauthor[]{oleary1980block}~\cite[Lemma 1]{oleary1980block}.
If the operator and preconditioner are symmetric the method is
equivalent to the BCG method.

Next, we introduce some theory about orthogonality of $\SubA$-valued
polynomials leading to some further properties of the BBiCG method.
Based on these properties we define the BBiCGStab method, which uses
enhanced polynomials that satisfy the same properties.

\begin{lem}
  \label{lem:polynomial_recursion}
The variables in the BBiCG algorithm can be expressed in the form
\begin{align}
  R^{k} = \mathcal{R}^{k}(A\inverse{\Prec})\circ R^{0}\\
  \intertext{and}
  P^{k} = \mathcal{P}^{k}(\inverse{\Prec}A)\circ \inverse{\Prec}R^{0},
\end{align}
where $\mathcal{R,P}\in\mathbb{P}_{\SubA}^{k}$ are $\SubA$-valued
polynomials of degree $k$  defined by the recursion formulas
\begin{align}
  \mathcal{R}^{0}(\polyX) &= \Identity & \mathcal{R}^{k+1}(\polyX) &= \mathcal{R}^{k}(\polyX) -
  \polyX\mathcal{P}^{k}(\polyX)\lambda^{k} \label{eq:recursion_p}\\
  \intertext{and}
  \mathcal{P}^{0}(\polyX) &= \Identity & \mathcal{P}^{k+1}(\polyX) &= \mathcal{R}^{k+1}(\polyX)
                                                           + \mathcal{P}^{k}(\polyX)
                                                           \beta^{k},
\end{align}
with $\lambda^{k} = \inverse{\left(\alpha^{k}\right)}\rho^{k}$ and $\beta^{k}=
\inverse{\left(\rho^{k}\right)}\rho^{k+1}$.
\end{lem}
\begin{proof}
  We proof this by induction. The case $k=0$ is trivial.
  Assume that the relations hold for $k$.
  Then we have
  \begin{align}
    R^{k+1} &= R^{k} -
    AP^{k}\lambda^{k}\\
    &= \mathcal{R}^{k}(A\inverse{\Prec})\circ R^{0} -
      \left(A\mathcal{P}^{k}(\inverse{\Prec}A)\circ
      \inverse{\Prec}R^{0}\right)\lambda^{k}\\
            &= \mathcal{R}^{k}(A\inverse{\Prec})\circ R^{0} -
              \left(A\inverse{\Prec}\mathcal{P}^{k}(A\inverse{\Prec})
              \circ R^{0}\right)\lambda^{k}\\
    &= \left[\mathcal{R}^{k} -
      \polyX\mathcal{P}^{k}\lambda^{k}\right](A\inverse{\Prec})\circ R^{0}.
  \end{align}
  The second equation follows from the update formula of $P$
  \begin{align}
    P^{k+1} &= \inverse{\Prec}R^{k+1} +
              P^{k}\beta^{k}\\
            &= \inverse{\Prec}\mathcal{R}^{k+1}(A\inverse{\Prec})\circ R^{0} +
              \left(\mathcal{P}^{k}(\inverse{\Prec}A)\beta^{k}\right)\circ
              \inverse{\Prec}R^{0}\\
                &= \mathcal{R}^{k+1}(\inverse{\Prec}A)\circ \inverse{\Prec}R^{0} +
              \left(\mathcal{P}^{k}(\inverse{\Prec}A)\beta^{k}\right)\circ
                  \inverse{\Prec}R^{0}\\
            &=\left[\mathcal{R}^{k+1} +
              \mathcal{P}^{k}\beta^{k}\right] (\inverse{\Prec}A) \circ
              \inverse{\Prec}R^{0}.
  \end{align}
\end{proof}

Now we introduce some formality to describe orthogonal polynomials and show that
the polynomials of the recursion formulas of the BBiCG method satisfy these
orthogonality properties.
That formalism helps us to get rid of the transposed operator that must be
applied in the BBiCG method.
\begin{defn}[Formally orthogonal polynomials]
  We define the $\SubA$-valued linear functionals $\mathcal{C}$ and
  $\mathcal{C}^{(1)}$ on $\bb{P}_{\SubA}^{k}$ for any
  $\mathcal{P}\in\bb{P}_{\SubA}^{k}$ as
  \begin{align}
    \mathcal{C}(\mathcal{P}) &:=
    \sproduct[\SubA]{\transinv{\Prec}\tilde{R}^{0}}{\mathcal{P}(A\inverse{\Prec})\circ
                               R^{0}}\\
    \mathcal{C}^{(1)}(\mathcal{P}) &:= \mathcal{C}(\polyX\mathcal{P}).
  \end{align}
\end{defn}
\begin{lem}
  For the $\SubA$-valued polynomials $\mathcal{R}^{k}$ and
  $\mathcal{P}^{k}$ from Lemma \ref{lem:polynomial_recursion} we have
  \begin{align}
    \mathcal{C}(\mathcal{R}^{k}\mathcal{T}) &= 0 \label{eq:R_ortho}\\
    \intertext{and}
    \mathcal{C}^{(1)}(\mathcal{P}^{k}\mathcal{T}) &= 0 \label{eq:P_ortho}
  \end{align}
  for any $\mathcal{T}\in\bb{P}^{k-1}_{\SubA}$.
\end{lem}
\begin{proof}
  From Lemma \ref{lem:BCG_properties} we know that for all
  $i=0,\ldots,k-1$ it holds that
  \begin{align} \sproduct[\SubA]{(\transinv{\Prec}\transpose{A})^{i}\transinv{\Prec}\tilde{R}^{0}}{R^{k}}
    = 0.\\
    \intertext{It follows that}
    \sproduct[\SubA]{\transinv{\Prec}\tilde{R}^{0}}
    {(A\inverse{\Prec})^{i}\mathcal{R}^{k}(A\inverse{\Prec})\circ R^{0}}
    = 0,\\
    \intertext{hence}
    \mathcal{C}(\polyX^{i}\mathcal{R}^{k}) = 0.
  \end{align}
  Equation \eqref{eq:R_ortho} follows from the linearity of $\mathcal{C}$.

  Similarly, we proof Equation \eqref{eq:P_ortho}.
  By using
  \begin{align}
    AP^{k} &=
             \left(R^{k}-R^{k+1}\right)\inverse{\left(\lambda^{k}\right)}\\
    \intertext{we get}
    \mathcal{C}^{(1)}(\polyX^{i}\mathcal{P}) &= \mathcal{C}(\polyX^{i+1}\mathcal{P})\\
           &= \sproduct[\SubA]{\transinv{\Prec}\tilde{R}^{0}}
             {\left(A\inverse{\Prec}\right)^{i+1}\mathcal{P}^{k}(A\inverse{M})\circ R^{0})}\\
                                        &= \sproduct[\SubA]
    {\left(\transinv{M}\transpose{A}\right)^{i}\transinv{M}\tilde{R}^{0}}
    {AP^{k}}\\
           &=
             \sproduct[\SubA]
             {\left(\transinv{M}\transpose{A}\right)^{i}\transinv{M}\tilde{R}^{0}}
             {\left(R^{k}-R^{k+1}\right)\inverse{\left(\lambda^{k}\right)}}\\
           &= 0,
  \end{align}
  \sloppy
  where we used again Lemma \ref{lem:BCG_properties} and the fact that
  $\left(\transinv{M}\transpose{A}\right)^{i}\transinv{M}\tilde{R}^{0}
  \in
  \krylov[\SubA]{\transinv{M}\transpose{A}}{\transinv{M}\tilde{R}^{0}}{i}$.
  Equation \eqref{eq:P_ortho} follows by using the linearity of
  $\mathcal{C}^{(1)}$.
\end{proof}

For the stabilization we enhance the residual and
search direction by a polynomial $\mathcal{Q}^{k}$ as
\begin{align}
  R^{k} &=
          \left[\mathcal{R}^{k}\mathcal{Q}^{k}\right](A\inverse{M})\circ
          R^{0}\\
  P^{k} &=
          \left[\mathcal{P}^{k}\mathcal{Q}^{k}\right](\inverse{M}A)\circ
          \inverse{M}R^{0},
\end{align}
where $\mathcal{Q}^{k}\in\bb{P}_{\bbR}^{k}$ is a scalar
polynomial recursively defined by
\begin{align}
  \mathcal{Q}^0(\polyX) &= 1 & \text{and} &&\mathcal{Q}^{k+1}(\polyX) &= (1-\omega^{k}\polyX)\mathcal{Q}^{k}(\polyX).
\end{align}
The $\omega^{k}\in \bbR$ are chosen to minimize the residual norm.
By defining
\begin{align}
  S^{k} =
  \left[\mathcal{R}^{k+1}\mathcal{Q}^{k}\right](A\inverse{M})\circ R^{0},
\end{align}
we obtain the following update formulas
\begin{align}
  R^{k+1} &=
            \left[\mathcal{R}^{k+1}\mathcal{Q}^{k+1}\right](A\inverse{M})\circ
            R^{0}\\
          &= \left[\mathcal{R}^{k+1}\mathcal{Q}^{k} -
            \omega^{k}\polyX\mathcal{R}^{k+1}\mathcal{Q}^{k}\right](A\inverse{M})
            \circ R^{0}\\
          &= S^{k}-\omega^{k}A\inverse{M}S^{k},\\
P^{k+1} &=
            \left[\mathcal{P}^{k+1}\mathcal{Q}^{k+1}\right](\inverse{M}A)\circ
            \inverse{M}R^{0}\\
          &= \left[\mathcal{R}^{k+1}\mathcal{Q}^{k+1} +
            \mathcal{P}^{k}\mathcal{Q}^{k+1}\beta^{k}\right] (\inverse{M}A)\circ
            \inverse{M}R^{0}\\
          &= \inverse{M}R^{k+1} + \left[\mathcal{P}^{k}\mathcal{Q}^{k} -
            \omega^{k}\polyX\mathcal{P}^{k}\mathcal{Q}^{k}\right](\inverse{M}A)\circ
            \inverse{M}R^{0}\beta^{k}\\
        &= \inverse{M}R^{k+1} + \left(P^{k} - \omega^{k}
          \inverse{M}AP^{k}\right)\beta^{k},\\
  S^{k} &=
          \left[\mathcal{R}^{k+1}\mathcal{Q}^{k}\right](A\inverse{M})\circ
          R^{0}\\
          &= \left[\mathcal{R}^{k}\mathcal{Q}^{k} -
            \polyX\mathcal{P}^{k}\mathcal{Q}^{k}\lambda^{k}\right](A\inverse{M})\circ
            R^{0}\\
          &= R^{k}-AP^{k}\lambda^{k}.\\
\end{align}

Finally, we need to deduce the coefficients $\omega^{k},\lambda^{k}$
and $\beta^{k}$.
As mentioned before, the $\omega^{k}$ is determined by minimizing the
residual $R^{k+1} = S^{k}-\omega^{k}A\inverse{M}S^{k}$ in the
Frobenius norm, which yields
\begin{align}
  \omega^{k} = \frac{\sproduct[F]{A\inverse{\Prec}S^{k}}{S^{k}}}
  {\sproduct[F]{A\inverse{\Prec}S^{k}}{A\inverse{\Prec}S^{k}}}.
\end{align}
The other coefficients are determined by the formal orthogonality
condition for $\mathcal{R}^{k}$ and $\mathcal{P}^{k}$. We have
\begin{align}
  0 &= \mathcal{C}\left(\mathcal{R}^{k+1}\mathcal{Q}^{k}\right)\\
    &= \mathcal{C}\left(\mathcal{R}^{k}\mathcal{Q}^{k}\right) -
      \mathcal{C}\left(\polyX\mathcal{P}^{k}\mathcal{Q}^{k}\right)\lambda^{k}\\
    &= \sproduct[\SubA]{\transinv{\Prec}\tilde{R}^{0}}{R^{k}} -
      \sproduct[\SubA]{\transinv{\Prec}\tilde{R}^{0}}{AP^{k}}\lambda^{k}\\
  \Rightarrow \lambda^{k} &=
\inverse{\left(\sproduct[\SubA]{\transinv{\Prec}\tilde{R}^{0}}{AP^{k}}\right)}
                            \sproduct[\SubA]{\transinv{\Prec}\tilde{R}^{0}}{R^{k}}\\
  \intertext{and}
  0 &=
      \mathcal{C}^{(1)}(\mathcal{P}^{k+1}\mathcal{Q}^{k})\\
    &= \mathcal{C}^{(1)}(\mathcal{R}^{k+1}\mathcal{Q}^{k}) +
      \mathcal{C}^{(1)}(\mathcal{P}^{k}\mathcal{Q}^{k})\beta^{k}\\
    &=
      \sproduct[\SubA]{\transinv{\Prec}\tilde{R}^{0}}{A\inverse{\Prec}S^{k}}
      + \sproduct[\SubA]{\transinv{\Prec}\tilde{R}^{0}}{AP^{k}}\beta^{k}\\
  \Rightarrow \beta^{k} &=
 -\inverse{\left(\sproduct[\SubA]{\transinv{\Prec}\tilde{R}^{0}}{AP^{k}}\right)}
 \sproduct[\SubA]{\transinv{\Prec}\tilde{R}^{0}}{A\inverse{\Prec}S^{k}}.
\end{align}
Now we can formulate the BBiCGStab algorithm, see Algorithm \ref{alg:bbicgstab}.
\begin{algorithm}[tbp!]
  \caption{Block BiCGStab Method}
  \label{alg:bbicgstab}
  \begin{algorithmic}
    \State $R^{0} = B-AX^{0}$
    \State $P^{0} = \inverse{\Prec}R^{0}$
    \State Choose $\transinv{\Prec}\tilde{R}^{0}$ arbitrary
    (e.g.\ $\transinv{M}\tilde{R}^{0} = P^{0}$)
    \State $V^{0}  = P^{0}$
    \For{$k=0,\ldots$}
    \State $Q^{k} = AP^{k}$
    \State $\lambda^{k} =
    \inverse{\left(\sproduct[\SubA]{\transinv{\Prec}\tilde{R}^{0}}{Q^{k}}\right)}
    \sproduct[\SubA]{\transinv{\Prec}\tilde{R}^{0}}{R^{k}}$
    \State $S^{k} = R^{k} - Q^{k}\lambda^{k}$
    \State $Z^{k} = \inverse{\Prec}Q^{k}$
    \State $T^{k} = V^{k}-Z^{k}\lambda^{k}$
    \State $U^{k} = AT^{k}$
    \State $\omega^{k} =
    \frac{\sproduct[F]{U^{k}}{S^{k}}}{\sproduct[F]{U^{k}}{U^{k}}}$
    \State $X^{k+1} = X^{k} + P^{k}\lambda^{k} + \omega^{k}T^{k}$
    \State $R^{k+1} = S^{k} -\omega^{k}U^{k}$
    \State break if $\|R^{k+1}\| \leq \tol$
    \State $V^{k+1} = \inverse{\Prec}R^{k+1}$
    \State $\beta^{k} =
    -\inverse{\left(\sproduct[\SubA]{\transinv{\Prec}\tilde{R}^{0}}{Q^{k}}\right)}
    \sproduct[\SubA]{\transinv{\Prec}\tilde{R}^{0}}{U^{k}}$
    \State $P^{k+1} =
    V^{k+1}+\left(P^{k}-\omega^{k}Z^{k}\right)\beta^{k}$
    \EndFor
  \end{algorithmic}
\end{algorithm}
Note that $S$ and $R$ as well as $Q$ and $T$ and $V$ and $U$ can share
memory pairwise.
Furthermore, the algorithm does not apply the transposed of the
operator or preconditioner.

It is the subject of further investigation whether we could choose the stabilization
parameter $\omega^{k}$ in the space $\SubA$.
A naive implementation by the author did not work.
Furthermore, an adaption of the BiCGStab($\ell$) method would be interesting.
An approach was recently proposed by
\citeauthor[]{saito2014development}~\cite[]{saito2014development}.

\section{Residual Re-Orthonormalization}
Like in the BCG algorithm we can improve and stabilize the convergence by
adding a residual re-orthonormalization approach.
This helps to resolve rank-deficiencies in the residual.
Unfortunately, the BBiCGStab method is not as stable as the BCG
method, especially at lower precision.
We transform the intermediate residual using the normalizer, such that
\begin{align}
  \hat{S}^{k}\sigma^{k} &= S^{k},
\end{align}
where $\sigma^{k}$ is an element in the *-subalgebra $\sigma^{k}\in\SubA$.
All other variables are transformed similarly, using the same $\sigma^{k}$.
The transformed variables are denoted with a hat.

The recurrence coefficients $\hat{\lambda}$ and $\hat{\beta}$ are computed by
using the recurrence formulas of $P$ and $S$ and assuming that
$\sigma$ is invertible.
We have
\begin{align}
  \hat{P}^{k+1}\sigma^{k+1} &= P^{k+1}\\
                            &= V^{k+1} + \left(P^{k}-
                              \omega^{k} Z^{k}\right)\beta^{k}\\
                            &= \hat{V}^{k+1}\sigma^{k+1} +
                              \left(\hat{P}^{k}\sigma^{k} +
                              \omega^{k}\hat{Z}^{k}\sigma^{k}\right)\beta^{k}\\
                            &= \hat{V}^{k+1}\sigma^{k+1} +
                              \left(\hat{P}^{k} +
                              \omega^{k}\hat{Z}^{k}\right)\hat{\beta}^{k}\sigma^{k+1}
                              \intertext{where we define
                              $\hat{\beta}^{k}$ by}
                              \hat{\beta}^{k}\sigma^{k+1} &=
                              \sigma^{k}\beta^{k}.\\
  \intertext{This can be used to deduce a formula for $\hat{\beta}^k$}
                              \sigma^{k}\beta^{k} &=  \sigma^{k}
 \inverse{\left(\sproduct[\SubA]{\transinv{\Prec}\tilde{R}^{0}}{Q^{k}}\right)}
                                                    \sproduct[\SubA]{\transinv{\Prec}\tilde{R}^{0}}{U^{k}}\\
  &=
 \inverse{\left(\sproduct[\SubA]{\transinv{\Prec}\tilde{R}^{0}}{\hat{Q}^{k}}\right)}
    \sproduct[\SubA]{\transinv{\Prec}\tilde{R}^{0}}{\hat{U}^{k}}\sigma^{k+1}\\
    \Rightarrow \hat{\beta}^{k} &=  \inverse{\left(\sproduct[\SubA]{\transinv{\Prec}\tilde{R}^{0}}{\hat{Q}^{k}}\right)}
    \sproduct[\SubA]{\transinv{\Prec}\tilde{R}^{0}}{\hat{U}^{k}}.
\end{align}
Here we used, that $\hat{U}^{k}$ is transformed with respect to the new transformation
$\sigma^{k+1}$.
A similar computation can be made for $\hat{\lambda}^{k}$ which shows that we
have
\begin{align}
  \sigma^{k}\lambda^{k} &= \hat{\lambda}^{k}\sigma^{k}.
\end{align}
Thus
\begin{align}
  \hat{\lambda}^{k} = \inverse{\left(\sproduct[\SubA]{\transinv{\Prec}\tilde{R}^{0}}{\hat{Q}^{k}}\right)}
    \sproduct[\SubA]{\transinv{\Prec}\tilde{R}^{0}}{\hat{R}^{k}}.
\end{align}
The stabilization coefficient $\hat{\omega}^{k}$ is chosen to
minimize the transformed residual
\begin{align}
  \hat{\omega}^{k} &= \frac{\sproduct[F]{\hat{U}^{k}}{\hat{S}^{k}}}{\sproduct[F]{\hat{U}^{k}}{\hat{U}^{k}}}.
\end{align}

As with the BCG method, we apply the re-orthonormalization adaptively.
We decide on the basis of the diagonally scaled condition number
$\kappa_{D}(\sproduct[\SubA]{\hat{R}^{k}}{\hat{R}^{k}})$ whether we
re-orthogonalize the residual.
This quantity was a heuristically choice.
Additionally, we introduce a parameter $\eta$ to be able to tune the
re-orthonormalization behavior.
As shown in the following numeric section, this choice works properly.
The algorithm can be found in Algorithm \ref{alg:BlockBiCGStab_with_reortho}.
\begin{algorithm}
  \caption{BBiCGStab with Adaptive Residual Re-Orthonormalization}
  \label{alg:BlockBiCGStab_with_reortho}
  \begin{algorithmic}
    \State $R^{0} = B-AX^{0}$
    \If{$\eta \neq 0$}
    \State $\hat{R}^{0},\sigma^{0} = \normalizer[R^{0}]{\SubA}$
    \Else
    \State $\sigma^{0} = \Identity$
    \EndIf
    \State $\hat{P}^{0} = \inverse{\Prec}\hat{R}^{0}$
    \State Choose $\transinv{\Prec}\tilde{R}^{0}$
    (e.g.\ $\transinv{\Prec}\tilde{R}^{0}=\hat{P}^{0}$)
    \State $\hat{V}^{0} = \hat{P}^{0}$
    \For{$k=0,\ldots$}
    \State $\hat{Q}^{k} = A\hat{P}^{k}$
    \State $\hat{\lambda}^{k} =
    \inverse{\left(\sproduct[\SubA]{\transinv{\Prec}\tilde{R}^{0}}{\hat{Q}^{k}}\right)}
    \sproduct[\SubA]{\transinv{\Prec}\tilde{R}^{0}}{\hat{R}^{k}}$
    \State $\hat{Z}^{k} = \inverse{\Prec}\hat{Q}^{k}$
    \State $\hat{S}^{k} = \hat{R}^{k}
    -\hat{Q}^{k}\hat{\lambda}^{k}$
    \State $X^{k+\frac12} = X^{k} + \hat{P}\lambda^{k}\sigma^{k}$
    \If {$\eta\kappa_{D}\left(\sproduct[\SubA]{\hat{R}^{k}}{\hat{R}^{k}}\right)
      > \sqrt{\machineeps}$}
    \State $\hat{S}^k,\gamma^{k} \gets \normalizer[\hat{S}^{k}]{\SubA}$
    \State $\sigma^{k+1} = \gamma^{k}\sigma^{k}$
    \State $\hat{T}^{k} = \inverse{\Prec}\hat{S}^{k}$
    \Else
    \State $\hat{T}^{k} =
    \hat{V}^{k}-\hat{Z}^{k}\hat{\lambda^{k}}$
    \EndIf
    \State $\hat{U}^{k} = A\hat{T}^{k}$
    \State $\hat{\omega}^{k} =
    \frac{\sproduct[F]{\hat{U}^{k}}{\hat{S}^{k}}}{\sproduct[F]{\hat{U}^{k}}{\hat{U}^{k}}}$
    \State $X^{k+1} = X^{k+\frac12} + \hat{\omega}^{k}\hat{T}^{k}\sigma^{k}$
    \State $\hat{R}^{k+1} =
    \hat{S}^{k}-\hat{\omega}^{k}\hat{U}^{k}$
    \State break if $\|R^{k+1}\| \leq \tol$
    \State $\hat{V}^{k+1} = \inverse{\Prec}\hat{R}^{k+1}$
    \State $\hat{\beta}^{k} = -\inverse{\left(\sproduct[\SubA]{\transinv{\Prec}\tilde{R}^{0}}{\hat{Q}^{k}}\right)}
    \sproduct[\SubA]{\transinv{\Prec}\tilde{R}^{0}}{\hat{U}^{k}}$
    \State $\hat{P}^{k+1} = \hat{V}^{k+1} + \left(\hat{P}^{k} - \hat{\omega}^{k}\hat{Z^{k}}\right)\hat{\beta}^{k}$
    \EndFor
  \end{algorithmic}
\end{algorithm}

We avoid computing the inverse of $\gamma^{k}$, because it is badly conditioned
if $R^{k}$ is numerically rank deficient.
Therefore we apply the preconditioner directly to $\hat{S}^{k}$ to
compute $\hat{T}^{k}$ in the case where a re-orthonormalization is performed.

\section{Numerical Experiments}
\label{sec:blocknumericbicgstab}
Like for the other methods we perform an experiment to compare
the block-parallel and block-global methods for the BBiCGStab
method.
We used the same matrix and preconditioner as in the numerical
experiments for the BGMRes method and a re-orthonormalization
parameter of $\eta=100$.

The convergence behavior is not very smooth for the BiCGStab
method, therefore we decided to display the result in a table and not
in a plot.
The results are shown in Table \ref{tab:bicgstab_hybrid_vs_global}.
\begin{table}[tb!]
  \centering
    \caption[Convergence results of the BBiCGStab method.]{Convergence results of
the BBiCGStab method. Gray lines indicate that no convergence was achieved
within $500$ iterations.}
  \begin{tabular}{lSccc}
\toprule
               &     &  rate &  iterations &  re-orthogonalizations \\
Method & $p$ &       &             &                        \\
\midrule
block-global & 1   &\color{gray!70} 0.998 &\color{gray!70}         500 &\color{gray!70}                      0 \\
               & 2   &\color{gray!70} 0.998 &\color{gray!70}         500 &\color{gray!70}                      3 \\
               & 4   &\color{gray!70} 0.998 &\color{gray!70}         500 &\color{gray!70}                     18 \\
               & 8   &\color{gray!70} 1.023 &\color{gray!70}         500 &\color{gray!70}                     95 \\
               & 16  & 0.942 &         156 &                     59 \\
               & 32  & 0.815 &          59 &                     42 \\
               & 64  & 0.590 &          27 &                     23 \\
               & 128 & 0.428 &          13 &                     10 \\
               & 256 & 0.267 &          11 &                      8 \\
block-parallel & 1   & 0.929 &         131 &                      0 \\
               & 2   & 0.895 &          90 &                     58 \\
               & 4   & 0.829 &          54 &                     49 \\
               & 8   & 0.708 &          33 &                     30 \\
               & 16  & 0.649 &          24 &                     21 \\
               & 32  & 0.541 &          18 &                     15 \\
               & 64  & 0.437 &          14 &                     11 \\
               & 128 & 0.352 &          12 &                      9 \\
               & 256 & 0.267 &          11 &                      8 \\
\bottomrule
\end{tabular}
  \label{tab:bicgstab_hybrid_vs_global}
\end{table} We see a similar behavior as for the other methods - the larger the
$p$, the faster the convergence.
The convergence rate denotes the average factor by that the residual norm is
reduced per iteration.
However, we do not have a relationship between the convergence of the
$\frac{p}{q}$-block-parallel method and the corresponding
$p$-block-global method.
It seems that the block-global methods perform better with the BBiCGStab method.


Table \ref{tab:bbicgstab_iterations} shows the number of iterations and
re-orthonormalizations for the \texttt{Simon/raefsky3} problem from the
SuiteSparse Matrix collection \cite[]{timothy2011university}.
We used a block size of $s=32$ and the block method $\SubA[B]$.
The break criteria for the method was a reduction of the max-column norm by a
factor of $\num{e-7}$.
The results show that re-orthonormalization is necessary to achieve convergence,
but a larger re-orthonormalization parameter $\eta$ does not necessarily lead to
fewer iterations.
\begin{table}
  \centering
  \caption[Iterations of the BBiCGStab method with residual
re-orthonormalization.]{Iterations of the BBiCGStab method with residual
re-orthonormalization. Gray lines indicate that the method did not converged
within $2000$ iterations.}
  \begin{tabular}{Sccc}
\toprule
{} &  rate &  iterations &  re-orthogonalizations \\
$\eta$  &       &             &                        \\
\midrule
0.000   &\color{gray!70} 1.036 &\color{gray!70}        2000 &\color{gray!70}                      0 \\
0.001   & 0.967 &         299 &                    106 \\
0.010   & 0.967 &         294 &                    103 \\
0.100   & 0.967 &         286 &                    116 \\
1.000   & 0.968 &         295 &                    162 \\
10.000  & 0.964 &         258 &                    153 \\
100.000 & 0.963 &         261 &                    217 \\
\bottomrule
\end{tabular}
  \label{tab:bbicgstab_iterations}
\end{table}



%
\cleardoublepage


\epigraphhead[500]{\textit{We cannot solve the problems using the same kind of
    thinking we used when we created them.}\par\hfill\textsc{Albert Einstein}}

\part{Communication-aware Block Krylov Methods}
\label{part:communication}

\chapter{Challenges on Large Scale Parallel Systems}
\label{chap:testenvironment} In this part of the thesis we optimize the
algorithms presented in Part~\ref{part:blockkrylov} with respect to the global
communications, i.e.\ the inner products and orthonormalizations that must be
performed during the iteration.
This optimization is twofold.
On the one hand we reduce the number of synchronization points, i.e.\ compute
multiple block inner products simultaneously.
On the other hand we overlap the communication for the block inner products with
computation.

In the current research there are a couple of approaches to minimize the
communication overhead of Krylov methods.
The most popular are so called $s$-step Krylov methods
\cite[]{chronopoulos1989s,Chronopoulos1990sStep,hoemmen2010communication}.
These methods reduce the number of global communications by enlarging the Krylov
space by $s$ dimensions in every iteration.
The communication needed to find the minimizer in that space and
orthogonalize the basis could then be carried out chunk-wise.
While decreasing the number of messages by a factor of $s$,
these methods suffer from instabilities.
A lot of investigations have been made to mitigate these instabilities.
A rigorous analysis of the round-off errors in $s$-step methods can be found in
the PhD thesis of
\citeauthor[]{carson2015communication}~\cite[]{carson2015communication}.
The techniques to stabilize the methods are quite sophisticated and require
information about the spectrum of the operator.
Another problem is that we often use preconditioners that realize some kind of
coarse grid correction, which behave similar to a global communication.
Hence, the computation of the $s$-step basis would take as long as $s$ global
communications, which reduces the benefits of the method.
Our aim is to overlap the computation of the preconditioner with the block inner
product, such that every iteration effectively only needs one global
synchronization.
Recently, \citeauthor[]{eller2019scalable}~\cite{eller2019scalable} presented a
very elaborate performance study of pipelined Krylov solvers on modern HPC
hardware, showing that a clever organization of communication can yield a
significant performance improvement.

Another approach to solve linear systems while reducing the communication
overhead are asynchronous iterations,
cf. \citeauthor[]{chazan1969chaotic}~\cite[]{chazan1969chaotic} and
\citeauthor[]{frommer2000asynchronous}~\cite[]{frommer2000asynchronous}.
The idea behind these methods is that every process iterates in its own speed
and communicates the updates asynchronously.
In particular these methods are non-deterministic.
The strength of this method is that they can handle unreliable networks and
heterogeneous architectures well, as the computation proceeds even if one node
is delayed.

Finally, the block Krylov methods presented in Part \ref{part:blockkrylov}
already reduce the communication overhead as they reduce the number of
iterations.
This is used by the enlarged Krylov methods, like presented by
\citeauthor[]{grigori2017reducing}~\cite[]{grigori2017reducing} who make use of
this advantage for single right-hand side problems.

In the first section of this chapter we look at the implementation of the
non-blocking collective communication in our code.
Thereafter, we introduce a benchmark for quantifying the available time during a
collective communication for computations.
In the last section we look at the TSQR algorithm that is used to compute a QR
decomposition of a tall-skinny matrix with only one collective communication.

We use the terms ``global communication'', ``collective communication'' and
``reduction'' synonymous.
They all denote the communication procedure that is carried out to compute for
example a global sum.

All numerical tests in this part are carried out on $128$ nodes of the
supercomputer PALMAII at the University of Münster.
Each node is an Intel Xeon Gold 6140 18C (Skylake) with $\SI{2.3}{\giga\hertz}$
connected by a $\SI{100}{\giga\bit\per\s}$ Intel Omni Path network.
We use the optimized Intel MPI library Version 2018 Update 4 in all experiments.
Furthermore, we always use \texttt{I_MPI_ASYNC_PROGRESS=1}, if not stated
otherwise.
See Section \ref{sec:collectivecommunicationbenchmark} for a discussion about
that parameter.

\section{Implementation}
\label{sec:communication_implementation}
The technical foundation for asynchronous collective communication was
introduced with the MPI 3
standard\footnote{\url{https://www.mpi-forum.org/mpi-30/}}.
Within this standard definition, functions for non-blocking collective
communication were introduced.
These functions does not only allow to overlap the communication with
computation, but also to have multiple communications in progress
simultaneously.
In particular, the \cpp{MPI_Iallreduce} function which is needed for the
asynchronous computation of inner block products.
In the \dune framework, we introduced an abstraction layer for the
asynchronicity which makes use of the future-concept, shown in Listing
\ref{lst:future}.
It provides the \cpp{ready()} method for checking the status of the
communication, as well as the \cpp{wait()} method to block the execution until
the communication has finished.
In addition, the communicated data could be obtained with the \cpp{get()}
method.
The \cpp{Future} object encapsulate the communicated data and the
\cpp{MPI_STATUS} handle.
We made this available in the \dune framework as part of the \dune[common]
module.

\begin{c++}{The concept of a \cpp{Dune::Future}.}{lst:future}
template<class T>
class Future{
  bool valid();
  bool ready();
  void wait();
  T get();
};
\end{c++}

\begin{c++}{Example: How to use the Future interface of the %
    \cpp{ScalarProduct}.}{lst:scalarproduct}
  Future<Block<X>> gamma_future = sp->idot(x,y);
  op->apply(x,z); // This computation is overlapped with communication
  Block<X> gamma = gamma_future.get();
\end{c++}

Furthermore, we extended the \cpp{ScalarProduct} interface by a function
\cpp{idot}, that returns a \cpp{Future<Block<X>>} and a function
\cpp{inormalize}, computing the normalizer of a block vector that is
returned in a \cpp{Future<Block<X>>}.
This allows us to write concurrent code in a \cc~way, without thinking too much
about resource allocation and technical details.
An example how to use this interface is given in Listing \ref{lst:scalarproduct}.
For the sequential case, e.g.\ if no MPI is available we provide a fallback
implementation that implements a Future that does nothing but hold the
object.

As all the methods in this work, we plan to make the code publicly available in the
\dune[istl] module of the \dune framework.
A prototype implementation can be found in the GitLab repository of the
author\footnote{\url{https://gitlab.dune-project.org/nils.dreier/dune-common}}.
\pagebreak[0]
\section{Collective Communication Benchmark}
\label{sec:collectivecommunicationbenchmark} To quantify the costs of the
collective communication in our environment we implemented a benchmark.
The benchmark is inspired by the one presented by
\citeauthor[]{lawry2002comb}~\cite[]{lawry2002comb}.
As we want to overlap the communication with computation we measure how much of
the communication time is available for computations.
For that we proceed as follows.
We measure the time that is needed for initiating a global communication and
immediately waiting for it (\lstinline{MPI_Wait}).
We call this time the base time ($t_{\texttt{base}}$).
After that, we execute the same test again, but introduce a busy wait between the
initiation and the finalization of the global communication.
The time we measure is called the iteration time ($t_{\texttt{iter}}$).
The time we spend in the busy wait is called the work time
($t_{\texttt{work}}$).
To mitigate outliers, we execute all these measurements \num{1000} times and take
the average.
We start with a work time equal to one quarter of the base time.
After each repetition, the work time is doubled.
We repeat this until the iteration time is larger than two times the base time.
After that we determine the overhead as
\begin{align}
  t_{\texttt{ovhd}} = t_{\texttt{iter}}-t_{\texttt{work}},
\end{align}
and the available time for computation is
\begin{align}
  t_{\texttt{avail}} = t_{\texttt{base}} - t_{\texttt{ovhd}}.
\end{align}

The MPI standard only guaranties that progress in the communication is
made if calls to MPI are made.
However, some MPI implementations can be configured such that they
proceed even if no calls are made.
For Intel MPI this switch is called \texttt{I_MPI_ASYNC_PROGRESS}.
To distinguish this case, we introduce two tests.
In the first test \texttt{NB_sleep} we do not make any MPI calls in the
busy wait.
In the second case \texttt{NP_active} we call \cpp{MPI_Status} on the
\cpp{MPI_Request} during the busy wait.

\begin{figure}
  \centering
  \input{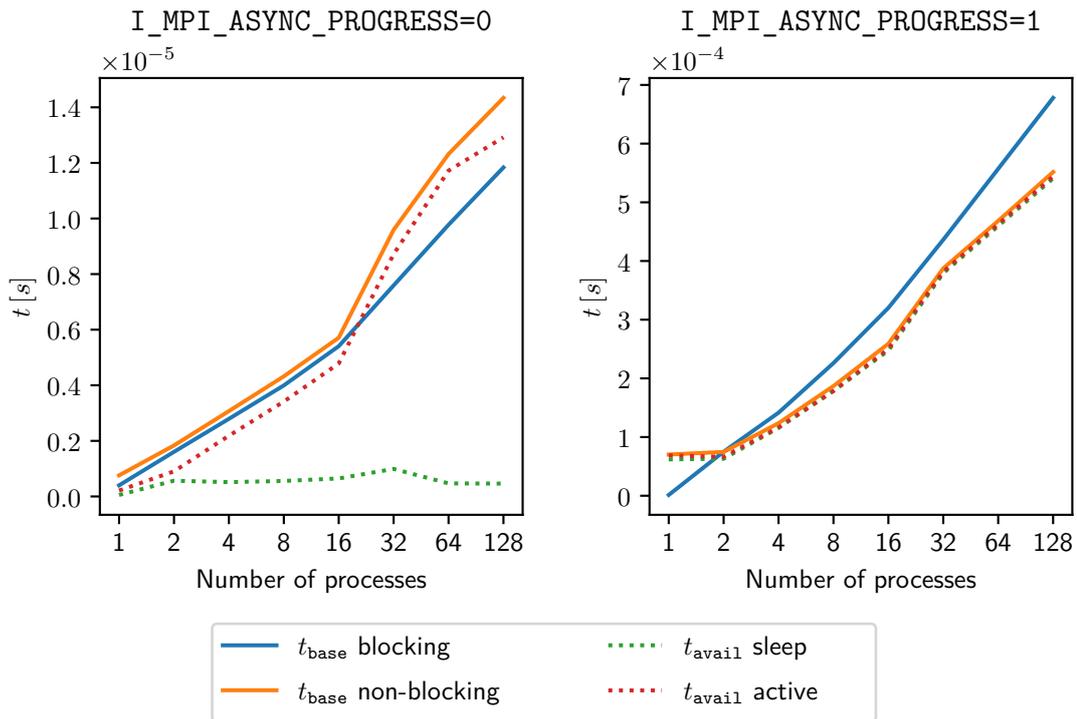}
  \caption[Collective communication benchmark.
  \texttt{I\_MPI\_ASYNC\_PROGRESS} on vs. off.]{Collective communication
benchmark.
    \texttt{I\_MPI\_ASYNC\_PROGRESS} on vs. off.
    Solid lines show the duration of the communication.
    Dotted lines show the portion of time that can be used for computations.
  }
  \label{fig:collective-communication-benchmark}
\end{figure}

Figure \ref{fig:collective-communication-benchmark} shows the results for our
test environment for the \cpp{MPI_Iallreduce} communication.
If \texttt{I_MPI_ASYNC_PROGRESS} is turned off, communication-computation
overlap is only feasible in the \texttt{NB_active} case.
In particular for a large number of processes most of the communication time can
be used for computation.
If \texttt{I_MPI_ASYNC_PROGRESS} is turned on, the \texttt{NB_active} and
\texttt{NB_sleep} case do not differ, i.e.\ communication computation overlap is
also feasible in the \texttt{NB_sleep} case, where no MPI calls are made during
the communication.
In this setting $99\%$ of the communication time can be used for computation.
However, the overall time in the case where \texttt{I_MPI_ASYNC_PROGRESS} is
turned on is much higher than in the other case.
Notice the different scaling of the y-axes.
From a practical point of view it does only make sense to activate
\texttt{I_MPI_ASYNC_PROGRESS} if the overlap is exploited aggressively, because
the overall communication time increase by a factor of approximately $50$.
An implementation of the benchmark can also be found in the GitLab repository of
the
author\footnote{\url{https://gitlab.dune-project.org/nils.dreier/dune-common/-/tree/master/dune/common/parallel/benchmark}}.

An alternative to turning \texttt{I_MPI_ASYNC_PROGRESS} on, was
presented by
\citeauthor[]{wittmann2013asynchronous}~\cite[]{wittmann2013asynchronous}.
In their approach, they spawn a dedicated thread that makes MPI calls
regularly to ensure progress for non-blocking communication, even if no calls to
MPI are made from the user code.


Let us now transfer these results to a hypothetical exascale machine, that
consists of the same nodes as our test environment.
One node in our test environment has a peak flop rate of
$\SI{2.65}{\tera\flop\per\second}$.
Hence, for an exascale machine $P\approx \num{380000}$ nodes are needed.
From Figure \ref{fig:collective-communication-benchmark} we can see that the
collective communication scale like $2\log_{2}(P)\si{\micro\second}$
(if \texttt{I_MPI_ASYNC_PROGRESS=0}).
Thus, a global reduction on this machine would take
\begin{align}
  2\log_{2}(380000)\si{\micro\second} \approx \SI{37}{\micro\second}.
\end{align}

\begin{figure}
  \centering
  \input{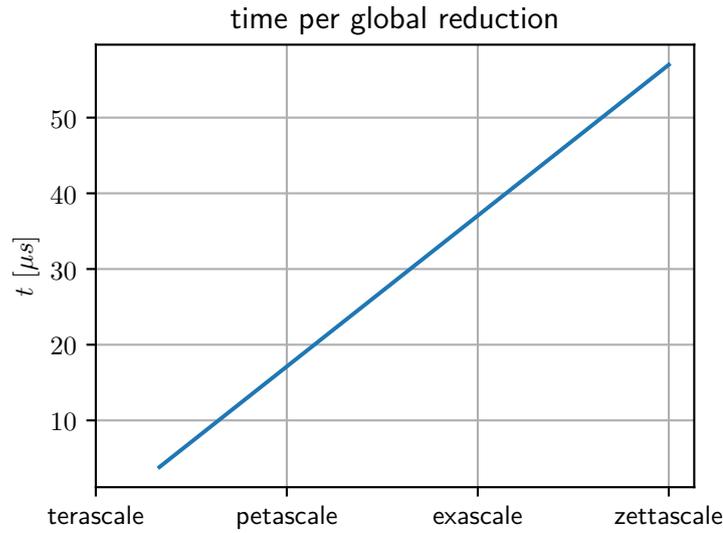}
  \caption[Time per global reduction vs. peak performance of the
    hypothetical machine.]{Time per global reduction vs. peak performance of the
    hypothetical machine.
  The $x$-axis is scaled logarithmically.}
  \label{fig:exascale_prediction}
\end{figure}

Figure \ref{fig:exascale_prediction} shows the time that is needed for a
global reduction on a hypothetical machine, as described above.
It shows that the costs of a global reduction will grow for larger
systems.
The costs increase from petascale to exascale by
$2\log_{2}(1000)\si{\micro\second} \approx \SI{20}{\micro\second}$,
equally from exascale to zettascale.
However, this assumes a quite optimal model.
It is not clear that the scaling is still true for machines of this
size, as small disturbance on one node would block the entire machine.

In particular, if the solution of a system is really time critical, e.g.\ weather
forecasts, the only possibility to solve it within the time constraints is to
use larger machines.
This would probably lead to relative small systems per node for which the
communication overhead plays a significant role.

\section{The TSQR Algorithm}
\label{sec:tall-skinny-qr}
Before we consider the block Krylov methods in the next chapters, we
look at the QR decomposition algorithm used in the distributed
setting to compute the normalizer.
Due to the shape of the block vectors for which the QR decomposition is
computed, this algorithm is called Tall-Skinny QR (TSQR).
It falls into the category of divide-and-conquer algorithms.

The algorithm was presented by
\citeauthor[]{demmel2008communication}~\cite[]{demmel2008communication,demmel2012communication}.
It reduces the communication for computing a QR decomposition
\begin{align}
  X = Q\sigma
\end{align}
of a block vector $X\in\bbR^{\dimA\times s}$, where $s \ll \dimA$,
$Q\in\bbR^{\dimA\times s}$ with $\transpose{Q}Q=\Identity$ and $\sigma
\in\bbR^{s\times s}$ is upper triangular.
We will review this algorithm and adapt it to our framework.
We start with the definition of the $\SubA$-QR decomposition.

\begin{defn}[$\SubA$-QR decomposition]
  Let $X\in\bbR^{\dimA\times s}$ be a block vector.  We call a
  decomposition of the type
  \begin{align}
    X &= Q\sigma &&& Q\in\bbR^{\dimA\times s},
  \end{align}
  where $\sigma\in\SubA$ is upper triangular and
  \begin{align}
    \sproduct[\SubA]{Q}{Q} &= \Identity
  \end{align}
  a $\SubA$-QR decomposition.
\end{defn}

Note that $\sigma$ is a $\SubA$-normalizer of $X$, with respect to the
normalized block vector $Q$.

The algorithm is based on the following computation.
Consider a block vector $X$ which is distributed onto $P$ processes
\begin{align}
  X = \begin{bmatrix}
    X_{1}\\
    \vdots\\
    X_{P}
  \end{bmatrix}
  \in\bbR^{\dimA\times s},
\end{align}
where $X_{p}$ is stored on process $p$.
This vector can be decomposed without communication into
\begin{align}
  \begin{bmatrix}
    X_{1}\\
    \vdots\\
    X_{P}
  \end{bmatrix} &= \begin{bmatrix}
    Q_{1} \sigma_{1}\\
    \vdots\\
    Q_{P}\sigma_{N}
  \end{bmatrix}
  = \begin{bmatrix}
    Q_{1}\\
    &\ddots\\
    &&Q_{P}
  \end{bmatrix}\begin{bmatrix}
    \sigma_{1}\\
    \vdots\\
    \sigma_{P}
  \end{bmatrix}
\end{align}
where
\begin{align}
\label{eq:local-qr}
  X_{p} &= Q_{p}\sigma_{p}
\end{align}
is a $\SubA$-QR decomposition for all $p=1,\ldots,P$.
The $R$-factors of these decompositions are gathered on one process and
decomposed into
\begin{align}
  \begin{bmatrix}
    \sigma_{1}\\
    \vdots\\
    \sigma_{P}
  \end{bmatrix} &= \begin{bmatrix}
    \tilde{\sigma}_{1}\\
    \vdots\\
    \tilde{\sigma}_{P}
  \end{bmatrix}\sigma &
                        \tilde{\sigma}_{1},\ldots,\tilde{\sigma}_{P},\sigma
                        \in \SubA,
\end{align}
with $\sigma$ is upper triangular and
\begin{align}
  \sum_{p=1}^{P} \transpose{\tilde{\sigma}_{p}}\tilde{\sigma}_{p}
  & =\Identity.
\end{align}
The $\tilde{\sigma}_{p}$ are scattered back to the respective processes.
Then the decomposition
\begin{align}
  X = \begin{bmatrix}
    Q_{1} \tilde{\sigma}_{1}\\
    \vdots\\
    Q_{P}\tilde{\sigma}_{N}
  \end{bmatrix}\sigma
\end{align}
forms a $\SubA$-QR decomposition, since
\begin{align}
  \sproduct[\SubA]{\begin{bmatrix}
    Q_{1} \tilde{\sigma}_{1}\\
    \vdots\\
    Q_{P}\tilde{\sigma}_{N}
  \end{bmatrix}}{\begin{bmatrix}
    Q_{1} \tilde{\sigma}_{1}\\
    \vdots\\
    Q_{P}\tilde{\sigma}_{N}
  \end{bmatrix}} = \sum_{p=1}^{P}
  \transpose{\tilde{\sigma}_{p}}\sproduct[\SubA]{Q_{p}}{Q_{p}}\tilde{\sigma}_{p}
  = \sum_{p=1}^{P}
  \transpose{\tilde{\sigma}_{p}}\tilde{\sigma}_{p} = \Identity
\end{align}
and $\sigma$ is upper triangular.

On large systems with many processes this approach can be carried out
recursively. In that case a tree hierarchy is created in the processes
and the local decomposition \eqref{eq:local-qr} is
computed with the same algorithm.
The communication pattern is the same as for a
\cpp{MPI_Allreduce}.
The MPI standard allows for implementing custom reduction methods.
This is sufficient for computing the $R$-factor, like presented in the
work of
\citeauthor[]{langou2010computing}~\cite[]{langou2010computing}.
Unfortunately, the MPI standard does not allow to define a custom
function for the back scattering of the $\tilde{\sigma}$, which would
be needed to compute the $Q$-factor in the rank deficient case.

We will see how we can use the idea of this algorithm to create a new
orthogonalization algorithm for BGMRes that is communication
optimal and stable.

An alternative to the TSQR algorithm is the CholeskyQR
algorithm~\cite[]{stathopoulos2002block}.
It also uses only one reduction to compute the tall-skinny QR factorization
of a block vector and relies on the Cholesky factorization of the block
inner product
\begin{align}
  \transpose{\lambda}\lambda &= \sproduct[\SubA]{X}{X}&\lambda\in\SubA,
                                                        \text{ upper
                                                        right triangular.}
\end{align}
The normalized vector could then be computed as
\begin{align}
  \label{eq:choleskyqr}
  Q &= X\transinv{\lambda}.
\end{align}
From \eqref{eq:choleskyqr} we see that it inverts the Cholesky
factor.
Therefore, it only works for full-rank $X$.
In particular, it is unsuitable for our stabilization strategies.


\chapter{Pipelined Block CG Method}
\label{chap:communicationcg} Our BCG algorithm with adaptive
re-orthonormalization, Algorithm \ref{alg:bcgaro}, uses three blocking global
communications per iteration -- minimization, orthogonalization and convergence
check/re-orthonormalization.
The present chapter aims at reducing this.
As a first step we fuse multiple inner products, such that the communication can
be carried out simultaneously, reducing the number of synchronization points per
iteration.
In the second section we make use of the non-blocking features described in the
last section to overlap the communication with computation.
Finally, we compare all the derived variants with respect to their communication
overhead.

\section{Fusing Inner Block Products}

We fuse the communication of the convergence check with the
orthogonalization, i.e.\ the communication for $\eta^{k}$ and $\rho^{k}$.
Furthermore, the update of $X$ could be delayed such that it can be done during
this communication.
These optimizations can be made without introducing any additional variables or
operations.
The resulting algorithm is shown in Algorithm \ref{alg:2R-CG}.
\begin{algorithm}
    \caption{BCG with Two Reductions (2R-BCG)}
  \label{alg:2R-CG}
  \begin{algorithmic}[1]
    \State $R^0 = B-AX^0$
    \If{$\eta > 0$}
    \State $\bar{R}^{0},\sigma^{0} = \normalizer[R^{0}]{\SubA}$
    \Else
    \State $\bar{R}^{0} = R^{0}$
    \State $\sigma^{0} = \Identity$
    \EndIf
    \State $P^0 = \inverse{\Prec}\bar{R}^0$
    \State $\rho^{0} = \sproduct[\SubA]{P^{0}}{\bar{R}^{0}}$
    \For{$k = 0,\ldots$ until convergence}
    \State $Q^k = AP^k$
    \State \initcommunication{alpha}$\alpha^k = \sproduct[\SubA]{P^k}{Q^k}$
    \State \finalizecommunication{alpha}$\lambda^k = \inverse{\left(\alpha^k\right)}\rho^{k}$
    \State $\tilde{R}^{k+1} = \bar{R}^{k} - Q^k\lambda^k$
    \If{$\eta \kappa_{D}(\alpha^{k}) > \sqrt{\machineeps}$}
    \State $\bar{R}^{k+1},\gamma^{k+1} = \normalizer[\tilde{R}^{k+1}]{\SubA}$
    \State $\sigma^{k+1}=\gamma^{k+1}\sigma^{k}$
    \State $\eta^{k+1} = \|\sigma^{k+1}\|$
    \Else
    \State $\bar{R}^{k+1} = \tilde{R}^{k+1}$
    \State $\gamma^{k+1} = \Identity$
    \State $\sigma^{k+1} = \sigma^{k}$
    \State \initcommunication{eta}$\eta^{k+1} = \|\bar{R}^{k+1}\sigma^{k+1}\|$
    \EndIf
    \State $X^{k+1} = X^{k} + P^k\lambda^k\sigma^{k}$ \label{algline:x_update}
    \State $Z^{k+1} = \inverse{\Prec}\bar{R}^{k+1}$
    \label{algline:z_update}
    \State \initcommunication{rho}$\rho^{k+1} = \sproduct[\SubA]{Z^{k+1}}{\bar{R}^{k+1}}$
    \State \finalizecommunication{eta}break if $\eta^{k+1} \leq \tol$
    \State \finalizecommunication{rho}$\beta^{k} = \inverse{\left(\rho^{k}\right)}\transpose{\gamma^{k+1}}\rho^{k+1}$
    \State $P^{k+1} = Z^{k+1} + P^k\beta^k$
    \EndFor
  \end{algorithmic}
  \drawcommunication{alpha/1.8em,eta/1.8em,rho/2.1em}
\end{algorithm} Arrows mark the initiation and finalization of the global
communications.
The algorithm is arithmetically equivalent to Algorithm \ref{alg:bcgaro}.
It could be slightly further improved by checking after line
\ref{algline:x_update} and \ref{algline:z_update} whether the communication of
$\eta^{k}$ is already finished and then break if the convergence criterion is
satisfied.
This would terminate the algorithm a bit earlier.

As a next step, we reduce the number of synchronization points to
one.
Unfortunately, that is not possible without adding memory and
computational overhead.
We introduce the auxiliary variable $U^{k} = AZ^{k}$, that can be computed
right after the application of the preconditioner.
Thus we can precompute $\alpha^{k+1}$ by
\begin{align}
  \alpha^{k+1} &= \sproduct[\SubA]{P^{k+1}}{Q^{k+1}}\\
             &= \sproduct[\SubA]{Z^{k+1} + P^{k}\beta^{k}}{U^{k+1} +
               Q^{k}\beta^{k}}\\
             &= \sproduct[\SubA]{Z^{k+1}}{U^{k+1}} +
               \sproduct[\SubA]{Z^{k+1}}{Q^{k}\beta^{k}} +
               \sproduct[\SubA]{P^{k}\beta^{k}}{U^{k+1}} +
               \sproduct[\SubA]{P^{k}\beta^{k}}{Q^{k}\beta^{k}}\\
             &= \delta^{k+1} +
               \sproduct[\SubA]{Z^{k+1}}{Q^{k}}\beta^{k} +
               \transpose{\left(\sproduct[\SubA]{Z^{k+1}}{Q^{k}}\beta^{k}\right)}
               +
               \transpose{\beta^{k}}\alpha^{k}\beta^{k} \label{eq:alpha_k+1}
\end{align}
with $\delta^{k+1} := \sproduct[\SubA]{Z^{k+1}}{U^{k+1}}$.
We simplify this further by computing
\begin{align}
  \sproduct[\SubA]{Z^{k+1}}{Q^{k}} &= \sproduct[\SubA]{Z^{k+1}}{\bar{R}^{k}
                                    -\tilde{R}^{k+1}}
                                     \inverse{\left(\lambda^{k}\right)}\\
                                   &=\sproduct[\SubA]{Z^{k+1}}{\bar{R}^{k}}
                                     \inverse{\left(\lambda^{k}\right)}
                                     -\sproduct[\SubA]{Z^{k+1}}{\tilde{R}^{k+1}}
                                     \inverse{\left(\lambda^{k}\right)}\\
                                   &=-\sproduct[\SubA]{Z^{k+1}}{\bar{R}^{k+1}}
                                     \gamma^{k+1}\inverse{\left(\rho^{k}\right)}
                                     \alpha^{k}\\
                                   &= -\transpose{\beta^{k}}\alpha^{k},
\end{align} where we used that $\sproduct[\SubA]{Z^{k+1}}{\bar{R}^{k}} =
\sproduct[\SubA]{\bar{R}^{k+1}}{Z^{k}} = 0$, as
$Z^{k}\in\krylov[\SubA]{\inverse{\Prec}A}{\inverse{\Prec}R^{0}}{k}$ and that
$\rho^{k+1}$ is symmetric.
Equation \eqref{eq:alpha_k+1} simplifies then to
\begin{align}
  \alpha^{k+1} &= \delta^{k+1} - \transpose{\beta^{k}}\alpha^{k}\beta^{k}.
\end{align}

We amend Algorithm \ref{alg:2R-CG} by computing $\delta^{k+1}$ together with
$\rho^{k+1}$ right after the preconditioner application.
The resulting algorithm can be found in Algorithm \ref{alg:1R-CG}.

\begin{algorithm}
  \caption{BCG with One Reduction (1R-BCG)}
  \label{alg:1R-CG}
  \begin{algorithmic}[1]
    \State $R^0 = B-AX^0$
    \If{$\eta > 0$}
    \State $\bar{R}^{0},\sigma^{0} = \normalizer[R^{0}]{\SubA}$
    \Else
    \State $\bar{R}^{0} = R^{0}$
    \State $\sigma^{0} = \Identity$
    \EndIf
    \State $P^0 = \inverse{\Prec}\bar{R}^0$
    \State $\rho^{0} = \sproduct[\SubA]{P^{0}}{\bar{R}^{0}}$
    \State $Q^0 = AP^0$
    \State $\alpha^0 = \sproduct[\SubA]{P^0}{Q^0}$
    \For{$k = 0,\ldots$ until convergence}
    \State $\lambda^k = \inverse{\left(\alpha^k\right)}\rho^{k}$
    \State $\tilde{R}^{k+1} = \bar{R}^{k} - Q^k\lambda^k$
    \If{$\eta \kappa_{D}(\alpha^{k}) > \sqrt{\machineeps}$}
    \State $\bar{R}^{k+1},\gamma^{k+1} = \normalizer[\tilde{R}^{k+1}]{\SubA}$
    \State $\sigma^{k+1}=\gamma^{k+1}\sigma^{k}$
    \State $\eta^{k+1} = \|\sigma^{k+1}\|$
    \Else
    \State $\bar{R}^{k+1} = \tilde{R}^{k+1}$
    \State $\gamma^{k+1} = \Identity$
    \State $\sigma^{k+1} = \sigma^{k}$
    \State \initcommunication{eta}$\eta^{k+1} = \|\bar{R}^{k+1}\sigma^{k+1}\|$
    \EndIf
    \State $X^{k+1} = X^{k} + P^k\lambda^k\sigma^{k}$
    \State $Z^{k+1} = \inverse{\Prec}\bar{R}^{k+1}$
    \State \initcommunication{rho}$\rho^{k+1} = \sproduct[\SubA]{Z^{k+1}}{\bar{R}^{k+1}}$
    \State $U^{k+1} = AZ^{k+1}$
    \State \initcommunication{delta}$\delta^{k+1} = \sproduct[\SubA]{Z^{k+1}}{U^{k+1}}$
    \State \finalizecommunication{eta}break if $\eta^{k+1} \leq \tol$
    \State \finalizecommunication{rho}$\beta^k =
    \inverse{\left(\rho^{k}\right)}\transpose{\gamma^{k+1}}\rho^{k+1}$
    \State $P^{k+1} = Z^{k+1} + P^k\beta^k$
    \State $Q^{k+1} = U^{k+1} + Q^k\beta^k$
    \State \finalizecommunication{delta}$\alpha^{k+1} = \delta^{k+1} - \transpose{\beta^{k}}\alpha^{k}\beta^{k}$
    \EndFor
  \end{algorithmic}
  \drawcommunication{delta/1.5em,rho/1.8em,eta/2.1em}
\end{algorithm}

Note that, compared to Algorithm \ref{alg:2R-CG}, Algorithm \ref{alg:1R-CG}
needs memory for $2$ more block vectors.
One for the additional variable $U^{k}$ and one because $Z^{k}$ and $Q^{k}$ can not share
their memory anymore, since $Q^{k}$ is updated recursively.

\section{Overlap Computation and Communication}

With the introduction of asynchronous communication techniques we can overlap
the communication procedure with computations, like the operator or
preconditioner application.
In this section, we amend the algorithms of the last section by integrating
these techniques.

All these optimizations introduce additional computational and memory
overhead.
The benefit of asynchronicity must compensate for this overhead,
otherwise Algorithm \ref{alg:bcgaro} would perform better.
Especially for very sparse matrices or many right-hand sides, the vector
updates that are introduced become easily similar expensive as the
operator application or the global communication.
In that case, Algorithm \ref{alg:2R-CG} is probably a good choice, as
it is fairly optimized and does not introduce additional computational
overhead.

As a first step, we develop an algorithm based on Algorithm
\ref{alg:2R-CG} that overlaps the two block inner products with the
operator and preconditioner application, respectively.
The approach was suggested by
\citeauthor[]{gropp2010update}~\cite[]{gropp2010update} for the non-block CG
method.

To do so, we introduce two additional variables $U^{k}=AZ^{k}$ and
$V^{k}=\inverse{\Prec}Q^{k}$.
These variables are computed during the global communication of the block inner
products and then used to update $Q^{k}$ and $Z^{k}$ recursively
\begin{align}
  Q^{k+1} &= U^{k+1} + Q^{k}\beta^{k},\\
  Z^{k+1} &= \inverse{\Prec}\bar{R}^{k+1}\\
          &= \inverse{\Prec}\tilde{R}^{k+1}\inverse{\left(\gamma^{k+1}\right)}\\
          &= \inverse{\Prec}\bar{R}^{k}\inverse{\left(\gamma^{k+1}\right)} -
            \inverse{\Prec}Q^{k}\lambda^{k}\inverse{\left(\gamma^{k+1}\right)}\\
          &= Z^{k}\inverse{\left(\gamma^{k+1}\right)}
            - V^{k}\lambda^{k}\inverse{\left(\gamma^{k+1}\right)}.
\end{align}

For the reasons mentioned in Section \ref{ssec:blockcg_residual_ortho}, it is
not favorable to invert $\gamma^{k+1}$.
In the iterations where a re-orthonormalization is made we rely on the direct computation of
$Z^{k+1} = \inverse{\Prec}\bar{R}^{k+1}$.
It means that the preconditioner is applied twice in that iteration.
However, this case occurs rarely.

In honor to \citeauthor[]{gropp2010update}~\cite[]{gropp2010update},
we call this algorithm \textit{Gropp's BCG}, it is shown in
Algorithm \ref{alg:gropps_block_cg}.

\begin{algorithm}
  \caption{Gropp's BCG}
  \label{alg:gropps_block_cg}
  \begin{algorithmic}
    \State $R^0 = B-AX^0$
    \If{$\eta > 0$}
    \State $\bar{R}^{0},\sigma^{0} = \normalizer[R^{0}]{\SubA}$
    \Else
    \State $\bar{R}^{0} = R^{0}$
    \State $\sigma^{0} = \Identity$
    \EndIf
    \State $P^0 = \inverse{\Prec}\bar{R}^0$
    \State $Q^0 = AP^0$
    \State $Z^{0} = P^{0}$
    \State $\rho^{0} = \sproduct[\SubA]{P^{0}}{\bar{R}^{0}}$
    \For{$k = 0,\ldots$ until convergence}
    \State \initcommunication{alpha}$\alpha^k = \sproduct[\SubA]{P^k}{Q^k}$
    \State $V^{k} = \inverse{\Prec}Q^{k}$
    \State \finalizecommunication{alpha}$\lambda^k = \inverse{\left(\alpha^k\right)}\rho^{k}$
    \State $\tilde{R}^{k+1} = \bar{R}^{k} - Q^k\lambda^k$
    \If{$\eta \kappa_{D}(\alpha^{k}) > \sqrt{\machineeps}$}
    \State $\bar{R}^{k+1},\gamma^{k+1} = \normalizer[\tilde{R}^{k+1}]{\SubA}$
    \State $\sigma^{k+1}=\gamma^{k+1}\sigma^{k}$
    \State $\eta^{k+1} = \|\sigma^{k+1}\|$
    \State $Z^{k+1} = \inverse{\Prec}\bar{R}^{k+1}$
    \Else
    \State $\bar{R}^{k+1} = \tilde{R}^{k+1}$
    \State $\gamma^{k+1} = \Identity$
    \State $\sigma^{k+1} = \sigma^{k}$
    \State \initcommunication{eta}initiate $\eta^{k+1} = \|\bar{R}^{k+1}\sigma^{k+1}\|$
    \State $Z^{k+1} = Z^{k} - V^{k}\lambda^{k}$
    \EndIf
    \State $X^{k+1} = X^{k} + P^k\lambda^k\sigma^{k}$
    \State \initcommunication{rho}$\rho^{k+1} = \sproduct[\SubA]{Z^{k+1}}{\bar{R}^{k+1}}$
    \State $U^{k+1} = AZ^{k+1}$
    \State \finalizecommunication{eta}break if $\eta^{k+1} \leq \tol$
    \State \finalizecommunication{rho}$\beta^k =
    \inverse{\left(\rho^{k}\right)}\transpose{\gamma^{k+1}}\rho^{k+1}$
    \State $P^{k+1} = Z^{k+1} + P^k\beta^k$
    \State $Q^{k+1} = U^{k+1} + Q^k\beta^k$
    \EndFor
  \end{algorithmic}
  \drawcommunication{alpha/0em,rho/0em,eta/1em}
\end{algorithm}

The next algorithm is based on the one reduction variant of our BCG method, Algorithm
\ref{alg:1R-CG}.
We use the same technique as before and introduce auxiliary variables to resolve
dependencies and precompute the operator or preconditioner application.
We reuse the variable $V^{k}$ and the update formula of $Z^{k}$ from the
previous algorithm, but update $V^{k}$ recursively, too.
For that we introduce the variable $W^{k+1} = \inverse{\Prec}U^{k+1}$, allowing
us update $V^{k}$ as
\begin{align}
  V^{k+1} &= \inverse{\Prec}Q^{k+1}\\
        &= \inverse{\Prec}U^{k+1} + \inverse{\Prec}Q^{k}\beta^{k}\\
        &= W^{k+1} + V^{k}\beta^{k}.
\end{align}
The variable $W^{k+1}$ can be computed during the communication of
the block inner products.
The resulting algorithm is shown in Algorithm \ref{alg:1RP-bcgaro}.

\begin{algorithm}
  \caption{Partially Pipelined BCG (PPBCG)}
  \label{alg:1RP-bcgaro}
  \begin{algorithmic}[1]
    \State $R^0 = B-AX^0$
    \If{$\eta > 0$}
    \State $\bar{R}^{0},\sigma^{0} = \normalizer[R^{0}]{\SubA}$
    \Else
    \State $\bar{R}^{0} = R^{0} \qquad \sigma^{0} = \Identity$
    \EndIf
    \State $Z^{0} = P^0 = \inverse{\Prec}\bar{R}^0$
    \State $\rho^{0} = \sproduct[\SubA]{P^{0}}{\bar{R}^{0}}$
    \State $Q^0 = AP^0$
    \State $V^{0} = \inverse{\Prec}Q^{0}$
    \State $\alpha^0 = \sproduct[\SubA]{P^0}{Q^0}$
    \For{$k = 0,\ldots$ until convergence}
    \State $\lambda^k = \inverse{\left(\alpha^k\right)}\rho^{k}$
    \State $\tilde{R}^{k+1} = \bar{R}^{k} - Q^k\lambda^k$
    \If{$\eta \kappa_{D}(\alpha^{k}) > \sqrt{\machineeps}$}
    \State $\bar{R}^{k+1},\gamma^{k+1} = \normalizer[\tilde{R}^{k+1}]{\SubA}$
    \State $\sigma^{k+1}=\gamma^{k+1}\sigma^{k}$
    \State $\eta^{k+1} = \|\sigma^{k+1}\|$
    \State $Z^{k+1} = \inverse{\Prec}\bar{R}^{k+1}$
    \Else
    \State $\bar{R}^{k+1} = \tilde{R}^{k+1}$
    \State $\gamma^{k+1} = \Identity$
    \State $\sigma^{k+1} = \sigma^{k}$
    \State \initcommunication{eta}initiate $\eta^{k+1} = \|\bar{R}^{k+1}\sigma^{k+1}\|$
    \State $Z^{k+1} = Z^{k} - V^{k}\lambda^{k}$
    \EndIf
    \State $X^{k+1} = X^{k} + P^k\lambda^k\sigma^{k}$
    \State \initcommunication{rho}$\rho^{k+1} = \sproduct[\SubA]{Z^{k+1}}{\bar{R}^{k+1}}$
    \State $U^{k+1} = AZ^{k+1}$
    \State \initcommunication{delta}$\delta^{k+1} = \sproduct[\SubA]{Z^{k+1}}{U^{k+1}}$
    \State $W^{k+1} = \inverse{\Prec}U^{k+1}$
    \State \finalizecommunication{eta}break if $\eta^{k+1} \leq \tol$
    \State \finalizecommunication{rho}$\beta^k = \inverse{\left(\rho^{k}\right)}\transpose{\gamma^{k+1}}\rho^{k+1}$
    \State $P^{k+1} = Z^{k+1} + P^k\beta^k$
    \State $Q^{k+1} = U^{k+1} + Q^k\beta^k$
    \State $V^{k+1} = W^{k+1} + V^{k}\beta^{k}$
    \State \finalizecommunication{delta}$\alpha^{k+1} = \delta^{k+1} - \transpose{\beta^{k}}\alpha^{k}\beta^{k}$
    \EndFor
  \end{algorithmic}
  \drawcommunication{delta/1.4em,rho/2em,eta/2.6em}
\end{algorithm}

If the application of the preconditioner does not suffice to hide the
global communication, it could be beneficial to apply the same
strategy to precompute the operator application.
For that we introduce the variables $S^{k} = AV^{k}$ and $T^{k} =
AW^{k}$.
The update of $S^{k}$ is similar to that of $V^{k}$, it reads
\begin{align}
  S^{k+1} &= AV^{k+1}\\
          &= AW^{k+1} + AV^{k}\beta^{k}\\
          &= T^{k+1} + S^{k}\beta^{k}.
\end{align} The resulting algorithm was suggested by
\citeauthor[]{ghysels2014hiding}~\cite[]{ghysels2014hiding,ashby2012impact} for
the non-block CG.
Therefore, we call the block variant Ghysels' BCG, It is shown in Algorithm
\ref{alg:1RFP-bcgaro}.

\begin{algorithm}
  \caption{Ghysels' BCG}
  \label{alg:1RFP-bcgaro}
  \begin{algorithmic}[1]
    \State $R^0 = B-AX^0$
    \If{$\eta > 0$}
    \State $\bar{R}^{0},\sigma^{0} = \normalizer[R^{0}]{\SubA}$
    \Else
    \State $\bar{R}^{0} = R^{0} \qquad \sigma^{0} = \Identity$
    \EndIf
    \State $Z^{0} = P^0 = \inverse{\Prec}\bar{R}^0$
    \State $\rho^{0} = \sproduct[\SubA]{P^0}{\bar{R}^{0}}$
    \State $U^{0} = Q^0 = AP^0$
    \State $V^{0} = \inverse{\Prec}Q^{0}$
    \State $S^{0} = AV^{0}$
    \State $\alpha^0 = \sproduct[\SubA]{P^0}{Q^0}$
    \For{$k = 0,\ldots$ until convergence}
    \State $\lambda^k = \inverse{\left(\alpha^k\right)}\rho^{k}$
    \State $\tilde{R}^{k+1} = \bar{R}^{k} - Q^k\lambda^k$
    \If{$\eta \kappa_{D}(\alpha^{k}) > \sqrt{\machineeps}$}
    \State $\bar{R}^{k+1},\gamma^{k+1} = \normalizer[\tilde{R}^{k+1}]{\SubA}$
    \State $\sigma^{k+1} = \gamma^{k+1}\sigma^{k}$
    \State $\eta^{k+1} = \|\sigma^{k+1}\|$
    \State $Z^{k+1} = \inverse{\Prec}\bar{R}^{k+1}$
    \State $U^{k+1} = AZ^{k+1}$
    \Else
    \State $\bar{R}^{k+1} = \tilde{R}^{k+1}$
    \State $\gamma^{k+1} = \Identity$
    \State $\sigma^{k+1} = \sigma^{k}$
    \State \initcommunication{eta}initiate $\eta^{k+1} = \|\bar{R}^{k+1}\sigma^{k+1}\|$
    \State $Z^{k+1} = Z^{k} - V^{k}\lambda^{k}$
    \State $U^{k+1} = U^{k} - S^{k}\lambda^{k}$
    \EndIf
    \State $X^{k+1} = X^{k} + P^k\lambda^k\sigma^{k}$
    \State \initcommunication{rho}$\rho^{k+1} = \sproduct[\SubA]{Z^{k+1}}{\bar{R}^{k+1}}$
    \State \initcommunication{delta}$\delta^{k+1} = \sproduct[\SubA]{Z^{k+1}}{U^{k+1}}$
    \State $W^{k+1} = \inverse{\Prec}U^{k+1}$
    \State $T^{k+1} = AW^{k+1}$
    \State \finalizecommunication{eta}break if $\eta^{k+1} \leq \tol$
    \State \finalizecommunication{rho}$\beta^k =
    \inverse{\left(\rho^{k}\right)}\transpose{\gamma^{k+1}}\rho^{k+1}$
    \State $P^{k+1} = Z^{k+1} + P^k\beta^k$
    \State $Q^{k+1} = U^{k+1} + Q^k\beta^k$
    \State $V^{k+1} = W^{k+1} + V^{k}\beta^{k}$
    \State $S^{k+1} = T^{k+1} + S^{k}\beta^{k}$
    \State \finalizecommunication{delta}$\alpha^{k+1} = \delta^{k+1} -
    \transpose{\beta^{k}}\alpha^{k}\beta^{k}$
    \EndFor
  \end{algorithmic}
  \drawcommunication{delta/1.4em,rho/2em,eta/2.6em}
\end{algorithm}

Figure \ref{fig:flowdiagramm} and Table \ref{tab:overview} give an overview
of the Algorithms deduced in this chapter.
In Figure \ref{fig:flowdiagramm} the schematic program flow of one iteration is
shown. Red boxes represent computations, green
ones represent communications.
The flow direction is from top to bottom.
Boxes that are placed next to each other horizontally mean, that these
operations are executed simultaneously.

Table \ref{tab:overview} gives an overview of the performance relevant
characteristics.
That is, how much memory is needed for the block vectors, how many
\texttt{BAXPY} operations are performed per iterations and the number
of synchronization points per iteration.
We see how much overhead is introduced for optimizing the communication properties.
\begin{figure}
  \centering
  \documentclass[tikz]{standalone}
\usepackage{pgfmath} 
\usepackage{amsfonts}
\input{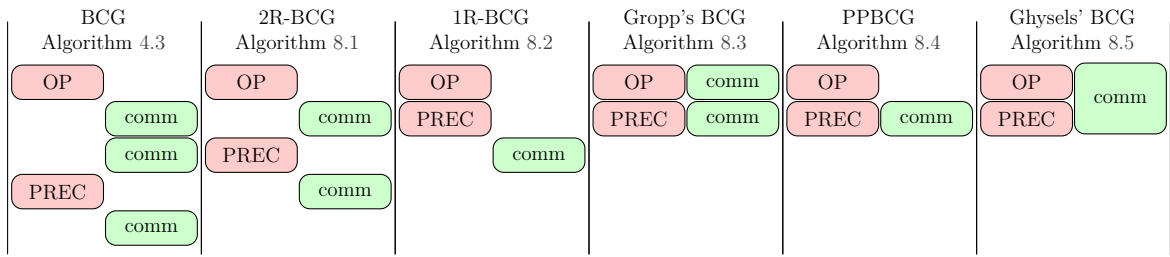}
\makeatletter
\newcommand*\ifcounter[1]{%
  \ifcsname c@#1\endcsname
    \expandafter\@firstoftwo
  \else
    \expandafter\@secondoftwo
  \fi
}
\makeatother
\newcounter{flowcounter}
\newcounter{step}
\begin{document}
\setcounter{flowcounter}{0}
\setcounter{step}{-1}

\newenvironment{flow}[1]{
  \def\bw{4.5em}
  \draw[black] ({(\theflowcounter+1)*(2*\bw+0.5em) -0.5*\bw -0.2em}, 3em) --
  ({(\theflowcounter+1)*(2*\bw+0.5em) - 0.5*\bw  -0.2em}, -8.5em);
  \tikzset{shiftblock/.style={xshift={\theflowcounter*(2*\bw+0.5em)},
      minimum width={\bw}}}
  \stepcounter{flowcounter}
  \draw[black] ({(\theflowcounter+1)*(2*\bw+0.5em) -0.5*\bw  -0.2em}, 3em) --
  ({(\theflowcounter+1)*(2*\bw+0.5em) - 0.5*\bw -0.2em}, -8.5em);
  \setcounter{step}{-1}
  \tikzset{block/.style={shiftblock,draw,rectangle,
      minimum height=1.7em,align=center,
      rounded corners,yshift=-1.8em*\thestep,
      /utils/exec={\stepcounter{step}}}}
  \tikzset{computation/.style={block,fill=red!20,}}
  \tikzset{communication/.style={block,fill=green!20,xshift=\bw+0.1em}}
  \tikzset{overlapped/.style={/utils/exec={\addtocounter{step}{-1}}}}
  \node[shiftblock, align=center] at (0.5*\bw, 2.5em) {#1};
}{}

\begin{tikzpicture}[scale = 0.65, transform shape]

  \begin{flow}{BCG\\Algorithm \ref{alg:bcgaro}}
      \node[computation] {OP};
      \node[communication] {comm};
      \node[communication] {comm};
      \node[computation] {PREC};
      \node[communication] {comm};
    \end{flow}

    \begin{flow}{2R-BCG\\Algorithm \ref{alg:2R-CG}}
      \node[computation] {OP};
      \node[communication] {comm};
      \node[computation] {PREC};
      \node[communication] {comm};
    \end{flow}

    \begin{flow}{1R-BCG\\Algorithm \ref{alg:1R-CG}}
      \node[computation] {OP};
      \node[computation] {PREC};
      \node[communication] {comm};
    \end{flow}

    \begin{flow}{Gropp's BCG\\Algorithm \ref{alg:gropps_block_cg}}
      \node[computation] {OP};
      \node[overlapped,communication] {comm};
      \node[computation] {PREC};
      \node[overlapped,communication] {comm};
    \end{flow}

    \begin{flow}{PPBCG\\Algorithm \ref{alg:1RP-bcgaro}}
      \node[computation] {OP};
      \node[computation] {PREC};
      \node[overlapped,communication] {comm};
    \end{flow}

    \begin{flow}{Ghysels' BCG\\Algorithm \ref{alg:1RFP-bcgaro}}
      \node[computation] {OP};
      \node[computation] {PREC};
      \node[overlapped,communication, minimum height=3.5em, yshift=1em] {comm};
    \end{flow}
  \end{tikzpicture}
\end{document}

  \caption[Schematic representation of the program flow in the
    different BCG variants.]{Schematic representation of the program flow in the
    different BCG variants.
    Red blocks are computational intensive procedures.
    Green blocks are communication.
  Horizontal alignment indicate concurrency.}
  \label{fig:flowdiagramm}
\end{figure}

\begin{table}
  \centering
  \caption{Arithmetical and communication properties of different BCG
algorithms.}
  \begin{tabular}{L{13em}ccccc}
    \toprule
    Method & vector storage & vector updates & synchronizations\\
    \midrule
    BCG (Alg.~\ref{alg:bcgaro}) & 4 & 3 & 3\\
    2R-BCG (Alg.~\ref{alg:2R-CG}) & 4 & 3 & 2\\
    1R-BCG (Alg.~\ref{alg:1R-CG}) & 6 & 4 & 1\\
    Gropp's BCG (Alg.~\ref{alg:gropps_block_cg}) & 6 & 5 & 2 (overlapped)\\
    PPBCG (Alg.~\ref{alg:1RP-bcgaro}) & 8 & 6 & 1 (overlapped)\\
    Ghysels' BCG (Alg.~\ref{alg:1RFP-bcgaro})& 10 & 8 & 1 (overlapped)\\
    \bottomrule
  \end{tabular}
  \label{tab:overview}
\end{table}

As the pipelining introduces additional \texttt{BAXPY} operations, it
also increases the sensitivity for round-off errors.
\citeauthor[]{cools2018analyzing}~\cite[]{cools2018analyzing} analyzed
the round-off errors for the pipelined CG method and proposed mitigation
strategies~\cite[]{cools2017numerically,cools2019numerically}.
These strategies include adaptive recomputation of the residual and
variables or introduction of shifts in the recurrence formulas.

For the pipelined block versions of the CG method we observe a similar
behavior.
In principle the strategies proposed by
\citeauthor[]{cools2019numerically} should be applicable to the
$\SubA$-BCG method as well.
However, the rigorous analysis and amendment of the algorithm is out
of the scope of this work.
Practical experiments showed that using a stronger preconditioner and using the
residual re-orthonormalization strategy mitigates the instabilities and
residual gap.

In extreme situations, it could be possible that the application of the operator
and preconditioner do not suffice to overlap the entire costs of a global
communication.
In these situations
\citeauthor[]{cornelis2018communication}~\cite[]{cornelis2018communication,cools2019improving}
proposed a deep-pipelined CG method, that overlays multiple iterations by a
global communication.
This method requires further additional \texttt{BAXPY} operations and memory
overhead.
As we consider systems with multiple right-hand sides, the operator and
preconditioner application is more expansive than in the non-block case.
Therefore, we do not consider deep-pipelining for our block methods.

\section{Numerical Experiments}

To quantify the advantages of the communication optimized BCG variants we run
strong scaling tests.
We already presented a similar study in \cite[]{bastian2020exa}.
The difficult part for running the tests is to choose a good problem size.
The communication costs and computation costs must be well-balanced to see the
differences of the algorithms.
We observed three different regimes.
Firstly, the computation costs are dominating the communication costs (e.g.\ in
the sequential case or for large problems).
Secondly, the communication and computation costs are similar.
This is the case where the overlap of communication and computation is
beneficial.
In the last case the communication costs are dominating (e.g.\ on a large
distributed system with a slow interconnect or for very small problems).
In the latter case we expect the method, that uses the least global reductions,
to be fastest.


With these considerations in mind we choose our test problem as a
plain Poisson problem, discretized with a $5$-point Finite Difference
stencil on a $500 \times 500$ grid.
We found this problem size experimentally by balancing the computational and
communicational costs for $P=16$ nodes.
We use a thread parallel block SSOR preconditioner.

\begin{figure}[t!]
  \centering
  \input{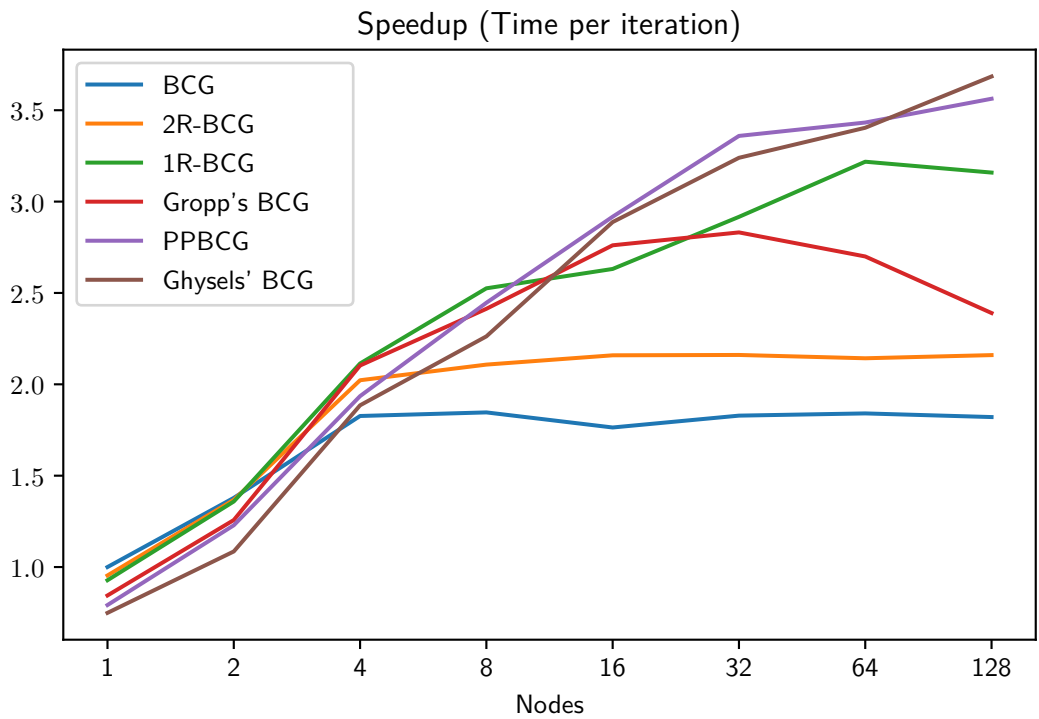}
  \caption[Strong scaling of the time per iteration for different CG
    variants.]{Strong scaling of the time per iteration for different CG
    variants.
  Speedup over sequential BCG.}
  \label{fig:parallel_cg}
\end{figure}

The result is shown in Figure \ref{fig:parallel_cg}.
It shows the speedup of the iterations with respect of the sequential case
($P=1$).
We can observe the three regimes.
For a smaller number of processes the BCG method (Algorithm \ref{alg:bcg}) and
2R-BCG (Algorithm \ref{alg:2R-CG}) are the fastest.
This matches our expectations as these method do not add computational overhead.
For a medium number of processes ($P\approx 16$) the methods that overlap
communication and computation are advantageous, which are Gropp's BCG, PPBCG and
Ghysels' BCG.
In the scaling limit the methods that use only one global communication per
iteration are fastest, which are 1R-BCG, PPBCG and Ghysels' BCG.

Note that the actual achieved speedup of $3.5\times$ for $P=128$ nodes is pretty
low.
This is due to the small problem size.
For $P=128$ nodes each node has approximately $\num{2000}$ DoFs.
As every node consists of $36$ cores, it results in approximately $55$ DoFs per
core.
Thus, the non-optimal scaling results are due to the bad performance of the
preconditioner for this tiny problem size.
However, we see that the optimization of global communication has an impact,
even if the scaling is hindered by other effects.


\chapter{Communication-optimized Block GMRes Method}
\label{chap:communicationgmres}

In this chapter, we analyze and optimize the communication effort of the BGMRes
method, introduced in Chapter \ref{chap:blockgmres}.
The orthogonalization in the Arnoldi process makes up the most communication
costs.
In particular, the communication overhead grows with every iteration.
In the literature, the classical Gram-Schmidt method is often used in
communication dominated settings.
This is fairly optimal with respect to communication but suffers from
instabilities for bad conditioned problems, cf. \citeauthor[]{giraud2005rounding}~\cite[]{giraud2005rounding}.

\citeauthor[]{hoemmen2010communication} presented in his
thesis~\cite[]{hoemmen2010communication} a Communication-Avoiding GMRes method
(CA-GMRes) that falls into the class of $s$-step methods.
That is, the Arnoldi procedure applies the operator multiple times and the
orthogonalization is carried out block-wise.
In principle, this approach could be combined with all strategies presented
in this chapter.
However, as with all $s$-steps methods, the use of CA-GMRes would require the
use of a stable $s$-step basis, which is not trivially constructed.

A pipelined method, like the one we have seen in the last chapter, was
presented by \citeauthor[]{ghysels2013hiding}~\cite[]{ghysels2013hiding}.
It overlaps the communication with the application of the operator and
preconditioner.
Nevertheless, as we consider block systems here, we concentrate on the Arnoldi
procedure, as it is more expensive than in the non-block case.
It should be costly enough to overlap the communication therein.

We develop and compare four variants of the BGMRes method in this chapter.
All methods alter the orthogonalization algorithm used in the Arnoldi
procedure.
We start with the classical Gram-Schmidt orthogonalization with
re-orthogonalization.
Thereafter, we see how we can overlap the communication with the computation of
the block inner products and block vector updates, leading to a pipelined
version of the Gram-Schmidt algorithm.
In addition, we develop a novel orthogonalization method based on the TSQR algorithm
that is optimal with respect to the communication overhead and preserves the
stability of the modified Gram-Schmidt algorithm.
We call this algorithm \textit{localized Arnoldi}, as it proceeds a local
orthogonalization procedure and then communicates the result with one reduction
communication to obtain the global result.
Finally, we present the result of a numerical test that compares all these methods.

\section{Classical Gram-Schmidt with Re-Orthogonalization}

\begin{algorithm}[tbp!]
  \caption{Modified Gram-Schmidt Orthogonalization}
  \label{alg:modified-gram-schmidt}
  \begin{algorithmic}
    \Input $X^{k+1}$ and $\SubA$-orthonormal basis $\mathcal{Q} = \left(Q^{i}\right)_{i=0}^{k}$
    \State $Q^{k+1} = X^{k+1}$
    \For{$i=0,\ldots,k$}
    \State $\rho_{i}^{k+1} = \sproduct[\SubA]{Q^{k+1}}{Q^{i}}$
    \State $Q^{k+1} \gets Q^{k+1} - Q^{i}\rho_{i}^{k+1}$
    \EndFor
    \State $Q^{k+1}, \rho_{k+1}^{k+1} \gets \normalizer[Q^{k+1}]{\SubA}$
    \State \Return $Q^{k+1}, \rho_{0}^{k+1},\ldots,\rho_{k+1}^{k+1}$
  \end{algorithmic}
\end{algorithm}
\begin{algorithm}[tbp!]
  \caption{Classical($\#it$) Gram-Schmidt Orthogonalization}
  \label{alg:classical-gram-schmidt}
  \begin{algorithmic}
    \Input $X^{k+1}$ and $\SubA$-orthonormal basis $\mathcal{Q} =
    \left(Q^{i}\right)_{i=0}^{k}$ and a parameter $\#it \in {1,2}$
    \State $Q^{k+1} = X^{k+1}$
    \For{$j=1,\ldots,\#it$}
    \For{$i=0,\ldots,k$}
    \State $\rho_{i}^{k+1} = \sproduct[\SubA]{Q^{k+1}}{Q^{i}}$
    \EndFor
    \For{$i=0,\ldots,k$}
    \State $Q^{k+1} \gets Q^{k+1} - Q^{i}\rho_{i}^{k+1}$
    \EndFor
    \EndFor
    \State $Q^{k+1}, \rho_{k+1}^{k+1} \gets \normalizer[Q^{k+1}]{\SubA}$
    \State \Return $Q^{k+1}, \rho_{0}^{k+1},\ldots,\rho_{k+1}^{k+1}$
  \end{algorithmic}
\end{algorithm}

The modified Gram-Schmidt method we used in Algorithm \ref{alg:bgmres} needs
$k+1$ global reductions to compute the next orthogonal basis vector.
Algorithm~\ref{alg:modified-gram-schmidt} shows the
modified Gram-Schmidt method.
It successively projects the new direction $X^{k+1}$ onto the orthogonal
compliments of the basis vectors $Q^{i}$.
The classical Gram-Schmidt method, shown in Algorithm
\ref{alg:classical-gram-schmidt}, does the same, but it computes the
projection factors $\rho_{i}^{k+1}$ in advance, such that the
communication for that could be carried out together.
Hence, it only needs two global communications -- one for the
orthogonalization and one for the normalization.

The numerical instabilities of the classical Gram-Schmidt method can be
mitigated by repeating the
orthogonalization~\cite{hoffmann1989iterative,bjorck1994numerics}.
In Algorithm \ref{alg:classical-gram-schmidt} this is implemented with the
parameter $\#it$.
\citeauthor[]{buhr2014numerically}~\cite[Algorithm 1]{buhr2014numerically}
proposed an algorithm that determines the number of re-orthogonalization
adaptively.
However, numerical experiments show that in most situations two iterations are
sufficient.

\section{Pipelined Gram-Schmidt}
\begin{wrapfigure}[21]{r}{0.4\textwidth}
  \centering
  \input{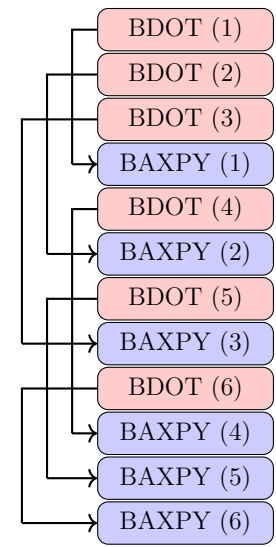}
  \caption[Schematic flow diagram of the pipelined Gram-Schmidt
orthogonalization Algorithm \ref{alg:pipelined-hybrid-gram-schmidt} for
$r=3$.]{Schematic flow diagram of the pipelined Gram-Schmidt orthogonalization
Algorithm \ref{alg:pipelined-hybrid-gram-schmidt} for $r=2$.
    Arrows indicate overlap of communication and computation.}
  \label{fig:pipelined_gramschmidt}
\end{wrapfigure}
Especially for many right-hand sides $s$, the \texttt{BAXPY} and
\texttt{BDOT} operations become equally expensive to or even more
expensive than the application of the operator or preconditioner.
In this case, it makes sense to overlap the reduction process with the
computation of the \texttt{BAXPY} and \texttt{BDOT} kernels.
\begin{algorithm}[tbp!]
  \caption{Pipelined Gram-Schmidt Orthogonalization}
  \label{alg:pipelined-hybrid-gram-schmidt}
  \begin{algorithmic}
    \Input $X^{k+1}$ and $\SubA$-orthonormal basis $\mathcal{Q} =
    \left(Q^{i}\right)_{i=0}^{k}$ and a parameter $1\leq r\leq k$
    \State $Q^{k+1} = X^{k+1}$
    \For{$i=0,\ldots,r-1$}
    \State $\rho_{i}^{k+1} = \sproduct[\SubA]{Q^{k+1}}{Q^{i}}$
    \EndFor
    \For{$i=r,\ldots,k$}
    \State $\rho_{i}^{k+1} = \sproduct[\SubA]{Q^{k+1}}{Q^{i}}$
    \State $Q^{k+1} \gets Q^{k+1} - Q^{i-r}\rho^{k+1}_{i-r}$
    \EndFor
    \For{$i=k-r,\ldots,k$}
    \State $Q^{k+1} \gets Q^{k+1} - Q^{i}\rho^{k+1}_{i}$
    \EndFor
    \State $Q^{k+1},\rho^{k+1}_{i+1} \gets \normalizer[Q^{k+1}]{\SubA}$
  \end{algorithmic}
\end{algorithm}
Algorithm \ref{alg:pipelined-hybrid-gram-schmidt} shows a pipelined
Gram-Schmidt method, that precomputes $r$ inner block products
before it starts to proceed the \texttt{BAXPY}s.
Every reduction process of the block inner product is overlapped with
$r$ \texttt{BAXPY} operations.
Figure \ref{fig:pipelined_gramschmidt} shows a schematic flow diagram
of the methods.
Arrows indicate the overlapping of communication and computation.
Another advantage is that not only computation is overlapped, but also $r+1$
global communications are computed in parallel. Hence, even in a communication
dominated environment, we could expect a speedup.

The case $r=0$ corresponds to the modified Gram-Schmidt process and the
case $r=k$ corresponds to the classical Gram-Schmidt process.
The algorithm is less stable for larger $r$.
To find a good parameter $r$, one must estimate a trade-of between
communication overlap and stability.

On unreliable networks, where the duration of the global communication
is not predictable, one could adapt Algorithm
\ref{alg:pipelined-hybrid-gram-schmidt} such that it checks in every
iteration whether a global communication is complete.
In that way, the parameter $r$ adapts automatically to the network performance.
The drawback is that this would lead to a non-deterministic behavior,
therefore we do not pursuit this strategy further.
This approach could be combined with the re-orthogonalization idea of
the previous section to achieve stability.

Up to the author's knowledge, this algorithm was not presented before, albeit the
basic idea is straight forward.
The reason might be that this method only makes sense, if the
\texttt{BAXPY} and \texttt{BDOT} operations consume a significant
amount of time which is the case for BGMRes with a relatively
sparse matrix or a lot of right-hand sides.

\section{Localized Arnoldi}
We now develop a orthogonalization method that is based on the
recursive structure of the TSQR algorithm.
Up to the authors knowledge, such algorithm was not presented before.
In principle, it is also useful in the non-block case.
For simplicity, we use the following notations:
\begin{itemize}
\item Vectors with elements in the *-subalgebra
  $\boldsymbol{\eta}\in\SubA^{k}$ are denoted by a small bold Greek
  letter.
\item square brackets $\left[\ldots\right]$ are used to concatenate
  matrices or vectors, both horizontally and vertically.
\end{itemize}

Our new method relies on an extended form of the block
QR decomposition. It assumes a distributed vector, i.e.\ 
\begin{align}
  X = \begin{bmatrix}
    X_{1}\\
    \vdots\\
    X_{P}
  \end{bmatrix},
\end{align} where $X_{p}$ is stored on processor $p$ for $p=1,\ldots,P$.
Let us start with the definition of the $\SubA$-QR decomposition for multiple block
vectors in the context of our block framework.
\begin{defn}[$\SubA$-QR decomposition for multiple block vectors]
  \label{defn:qr-decomposition}
  Let $\left(X^{i}\right)_{i=0}^{k} \in\bbR^{\dimA \times s}$,
  $\left(Q^{i}\right)_{i=0}^{k} \in\bbR^{\dimA \times s}$ and
  $\boldsymbol{\rho}^{0}, \ldots, \boldsymbol{\rho}^{k} \in \SubA^{k+1}$,
  then a decomposition of the form
  \begin{align}
    \left[X^{0}\ldots X^{k}\right] = \left[Q^{0}\ldots
    Q^{k}\right]\left[\boldsymbol{\rho}^{0} \ldots \boldsymbol{\rho}^{k}\right]
  \end{align} is called a $\SubA$-QR decomposition, if
$\left[\boldsymbol{\rho}^{0}\ldots \boldsymbol{\rho}^{k}\right]$ is upper
triangular and $Q$ satisfies
  \begin{align}
    \sproduct[\SubA]{Q^{i}}{Q^{j}} &= \delta_{ij}\Identity & \forall i,j=0,\ldots,k,
  \end{align}
  where $\delta_{ij}$ denotes the Kronecker symbol.
\end{defn}

The next definition generalizes the $\SubA$-QR decomposition for the distributed
setting.
\begin{defn}[Distributed $\SubA$-QR decomposition]
  \label{defn:distributed-qr-decomposition}
  A decomposition of the form
  \begin{align}
    \begin{bmatrix}
      X^{0}_{1} & \ldots & X^{k}_{1}\\
      \vdots & & \vdots\\
      X^{0}_{P} & \ldots & X^{k}_{P}
    \end{bmatrix} &= \begin{bmatrix}
      \left(\left[P^{0}_{1}\ldots P^{k}_{1}\right]\left[\boldsymbol{\zeta}_{1}^{0}
        \ldots\boldsymbol{\zeta}_{1}^{k}\right]\right)\\
      \vdots\\
      \left(\left[P^{0}_{P}\ldots
          P^{k}_{P}\right]\left[\boldsymbol{\zeta}_{P}^{0}\ldots
          \boldsymbol{\zeta}_{P}^{k}\right]\right)
    \end{bmatrix}\left[\boldsymbol{\rho}^{0}\ldots\boldsymbol{\rho}^{k}\right]
  \end{align}
  is called a distributed $\SubA$-QR decomposition, if for all
  $l=1,\ldots,P$ holds
  \begin{align}
    \sproduct[\SubA]{P^{i}_{l}}{P^{j}_{l}} &= \delta_{ij}\Identity
    \label{eq:local-ortho}
    \intertext{and}
    \label{eq:global-qr}
    \begin{bmatrix}
      \widetilde{\boldsymbol{\zeta}}^{0}_{1} & \ldots &
      \widetilde{\boldsymbol{\zeta}}^{k}_{i}\\
      \vdots & & \vdots\\
      \widetilde{\boldsymbol{\zeta}}^{0}_{P} & \ldots &
      \widetilde{\boldsymbol{\zeta}}^{k}_{P}
    \end{bmatrix} &= \begin{bmatrix}
      \boldsymbol{\zeta}^{0}_{1} & \ldots & \boldsymbol{\zeta}^{k}_{i}\\
      \vdots & & \vdots\\
      \boldsymbol{\zeta}^{0}_{P} & \ldots & \boldsymbol{\zeta}^{k}_{P}
    \end{bmatrix}\left[\boldsymbol{\rho}^{0}\ldots\boldsymbol{\rho}^{k}\right]
  \end{align} is a $\SubA$-QR decompostion with respect to the inner
  block product
  \begin{align}
    \llangle \begin{bmatrix}\boldsymbol{\zeta}_{1}\\
      \vdots\\
      \boldsymbol{\zeta}_{P}
    \end{bmatrix},
    \begin{bmatrix}\boldsymbol{\tau}_{1}\\
      \vdots\\
      \boldsymbol{\tau}_{P}
    \end{bmatrix} \rrangle_{\SubA} := \sum_{i=1}^{P}
    \transpose{\boldsymbol{\zeta}_{i}}\boldsymbol{\tau}_{i} \in \SubA.
  \end{align}
  The vectors
  $\widetilde{\boldsymbol{\zeta}}^{0}_{p},\ldots,\widetilde{\boldsymbol{\zeta}}_{p}^{k}
  \in \SubA^{k+1}$ are called the local $R$-factors.
\end{defn}

\begin{rem}\leavevmode
  \begin{itemize}
  \item The local parts might be rank-deficient, even if the global
    system $\left[X^{0}\ldots X^{k}\right]$ has full rank.
    In particular, it might occur that the local $R$-factors
    $\left[\widetilde{\boldsymbol{\zeta}}_{p}^{0}\ldots
      \widetilde{\boldsymbol{\zeta}}_{p}^{k}\right]$
    are not invertible.
  \item Equation \eqref{eq:local-ortho} needs a definition of the
    block inner product on the local part of the vector space.
    We use the definition
    \begin{align}
      \sproduct[\SubA]{P_{l}}{Q_{l}} := \sproduct[\SubA]{
      \begin{bmatrix}
        0\\
        \vdots\\
        0\\
        P_{l}\\
        0\\
        \vdots\\
        0
      \end{bmatrix}
      }{
      \begin{bmatrix}
        0\\
        \vdots\\
        0\\
        Q_{l}\\
        0\\
        \vdots\\
        0
      \end{bmatrix}
      }.
      \label{eq:local-block-inner-product}
    \end{align}
    However, as we will see, the only property that the local block
    inner product needs to satisfy for our theory is
    \begin{align}
      \sproduct[\SubA]{
      \begin{bmatrix}
        P_{1}\\
        \vdots\\
        P_{P}
      \end{bmatrix}
      }{
      \begin{bmatrix}
        Q_{1}\\
        \vdots\\
        Q_{P}
      \end{bmatrix}
      } = \sum_{l=1}^{P} \sproduct[\SubA]{P_{l}}{Q_{l}},
    \end{align}
    which is true for definition \eqref{eq:local-block-inner-product}
    and the block inner products of Definition \ref{defn:subalgebras}.
\end{itemize}
\end{rem}

The next lemma shows that a distributed $\SubA$-QR decompostion is
just a special case of the $\SubA$-QR decomposition.
\begin{lem}
  A distributed $\SubA$-QR decomposition (Definition \ref{defn:distributed-qr-decomposition}) is a $\SubA$-QR decomposition
  (Definition \ref{defn:qr-decomposition}).
\end{lem}
\begin{proof}
  We have to show that the columns of the matrix
  \begin{align}
    \begin{bmatrix}
      \left(\left[P^{0}_{1}\ldots
          P^{k}_{1}\right]\left[\boldsymbol{\zeta}_{1}^{0}
          \ldots\boldsymbol{\zeta}_{1}^{k}\right]\right)\\
      \vdots\\
      \left(\left[P^{0}_{P}\ldots P^{k}_{P}\right]
        \left[\boldsymbol{\zeta}_{P}^{0}\ldots\boldsymbol{\zeta}_{P}^{k}\right]\right)
    \end{bmatrix}
  \end{align}
  are $\SubA$-orthonormal and that $\left[\boldsymbol{\rho}^{0}\ldots\boldsymbol{\rho}^{k}\right]$ is upper triangular.
  The latter follows from the fact that \eqref{eq:global-qr} is a $\SubA$-QR
  decomposition.
  For the $i$th and $j$th column of the $Q$-factor it holds
  \begin{align}
    \sproduct[\SubA]{\begin{bmatrix}\left[P_{1}^{0}\ldots
          P_{1}^{k}\right] \boldsymbol{\zeta}_{1}^{i}\\
        \vdots\\
        \left[P_{P}^{0}\ldots
          P_{P}^{k}\right] \boldsymbol{\zeta}_{P}^{i}
      \end{bmatrix}
    }
      {\begin{bmatrix}\left[P_{1}^{0}\ldots
          P_{1}^{k}\right] \boldsymbol{\zeta}_{1}^{j}\\
        \vdots\\
        \left[P_{P}^{0}\ldots
          P_{P}^{k}\right] \boldsymbol{\zeta}_{P}^{j}
      \end{bmatrix}} &=
                       \sum_{p=1}^{P}
                       \sproduct[\SubA]{\left[P_{p}^{0}\ldots
                       P_{p}^{k}\right]\boldsymbol{\zeta}_{p}^{i}}
                       {\left[P_{p}^{0}\ldots P_{p}^{k}\right]\boldsymbol{\zeta}_{p}^{j}}\\
    &= \sum_{p=1}^{P}
      \transpose{\boldsymbol{\zeta}^{i}_{p}}\boldsymbol{\zeta}^{j}_{p}\\
    &= \llangle \boldsymbol{\zeta}^{i}, \boldsymbol{\zeta}^{j}\rrangle_{\SubA}= \delta_{ij}\Identity,
  \end{align}
  where we used that $\left(P^{i}_{p}\right)_{i=0}^{k}$ as well as
  $\left(\begin{bmatrix}
      \boldsymbol{\zeta}_{1}^{i}\\
      \vdots\\
      \boldsymbol{\zeta}_{P}^{i}
    \end{bmatrix}
  \right)_{i=0}^{k}$ are $\SubA$-orthonormal systems.
\end{proof}

As we want to apply this method in the Arnoldi iteration, we assume that
we have given a distributed $\SubA$-QR decomposition and want to
extend it by another vector $X^{k+1}$.
Inspired by the TSQR algorithm presented in Section
\ref{sec:tall-skinny-qr}, the idea is that a distributed $\SubA$-QR
decomposition can be extended by the direction $X^{k+1}$ by computing
a local $\SubA$-QR decomposition, followed by a global reduction of
the local $R$-factors.

Let the $Q$-factor
\begin{align}
  \begin{bmatrix}
    \left(\left[P_{1}^{0} \ldots
      P_{1}^{k}\right]\left[\boldsymbol{\zeta}_{1}^{0}\ldots
      \boldsymbol{\zeta}_{1}^{k}\right]\right)\\
    \vdots\\
    \left(\left[P_{P}^{0} \ldots
      P_{P}^{k}\right]\left[\boldsymbol{\zeta}_{P}^{0}\ldots
      \boldsymbol{\zeta}_{P}^{k}\right]\right)
  \end{bmatrix}
\end{align}
of a distributed $\SubA$-QR decomposition and a distributed block
vector $X$ be given.

The aim of the \textit{localized Arnoldi method} is to compute
$P_{1}^{k+1},\ldots,P_{P}^{k+1}$,
$\boldsymbol{\zeta}_{1}^{k+1},\ldots,\boldsymbol{\zeta}_{P}^{k+1}$
and $\boldsymbol{\rho}^{k+1}$ such that
\begin{align}
  \begin{bmatrix}
    X_{1}\\
    \vdots\\
    X_{P}
  \end{bmatrix}
  &=
    \begin{bmatrix}
      \left(\left[P_{1}^{0}\ldots
        P_{1}^{k+1}\right]
      \begin{bmatrix}
        \boldsymbol{\zeta}_{1}^{0} & \ldots &
        \boldsymbol{\zeta}_{1}^{k} &
        \multirow{2}{*}{$\boldsymbol{\zeta}_{1}^{k+1}$}\\
        \transpose{0} & \ldots &\transpose{0}
      \end{bmatrix}\right)\\
      \vdots\\
      \left(\left[P_{P}^{0}\ldots
        P_{P}^{k+1}\right]
      \begin{bmatrix}
        \boldsymbol{\zeta}_{P}^{0} & \ldots &
        \boldsymbol{\zeta}_{P}^{k} &
        \multirow{2}{*}{$\boldsymbol{\zeta}_{P}^{k+1}$}\\
        \transpose{0} & \ldots &\transpose{0}
      \end{bmatrix}\right)\\
    \end{bmatrix}\boldsymbol{\rho}^{k+1}
  \intertext{and}
  \sproduct[\SubA]{P_{l}^{i}}{P_{l}^{k+1}}
  &= \delta_{i,k+1}\Identity \qquad \forall l=1,\ldots,P,\quad i=0,\ldots,k+1\\
  \llangle \boldsymbol{\zeta}^{i}, \boldsymbol{\zeta}^{k+1}\rrangle_{\SubA}
  &= \delta_{i,k+1}\Identity  \qquad \forall i=0,\ldots,k+1.
\end{align} The local $Q$-factors $P_{l}^{k+1}$ are computed by applying any
stable orthogonalization method locally.
We use the modified Gram-Schmidt method in our examples.
After that, the global $\SubA$-QR decomposition of the local $R$-factors must be
computed.
This can be achieved by either gather it on one master process, or perform it
recursively in a reduction procedure.
The latter is preferable on large scale machines.
We describe the recursive procedure in the following.

The reduction procedure is performed on a tree on which every node stores its
local orthogonal basis
\begin{align}
  \left(\begin{bmatrix}
      \boldsymbol{\zeta}^{j}_{0}\\
      \boldsymbol{\zeta}^{j}_{1}
    \end{bmatrix}
  \right)_{j=0}^{k} \subset \SubA^{2k}.
\end{align} In particular, every node hold a state that is extended from
iteration to iteration.
Hence, in every iteration of the Arnoldi method, the same tree must be used.
At the beginning of every iteration, every process orthonormalizes its local
part of the new block vector $P^{k+1}_{p}$ to its local orthonormal basis.
The local R-factor $\tilde{\boldsymbol{\zeta}}^{k+1}_{p}\in\SubA^{k}$ is then send to its
parent.
During the reduction operation, every node in the reduction tree receives to
vectors $\tilde{\boldsymbol{\zeta}}^{k+1}_{0},
\tilde{\boldsymbol{\zeta}}^{k+1}_{1}$, that are stacked on each other and then
orthonormalized to its local basis.
The resulting R-factor $\boldsymbol{\rho}^{k+1}$ is send to its parent.
The root node does the same, but send its local orthonormal Q-factors
$\boldsymbol{\zeta}^{k+1}_{0}, \boldsymbol{\zeta}^{k+1}_{1}$, back to its
children.
This initiates the back-propagation.
The resulting $R$-factor of the root is the global one.
Every node in the tree receives the $Q$-factor from its parent and multiplies it
with its local $Q$-factors.
The result is then back-propagates to the children.
On the leafs the globally orthogonal basis can be obtained by multiplying the
received $Q$-factor to the locally orthogonal part of the block vector
$P^{k+1}_{p}$.

\begin{algorithm}[t]
  \caption{Reduction Process for the Localized Arnoldi Method}
  \label{alg:reduction-localized}
  \begin{algorithmic}
    \State Receive $\tilde{\boldsymbol{\zeta}}^{k+1}_{0}$ and
    $\tilde{\boldsymbol{\zeta}}^{k+1}_{1}$
    from children
    \For{$i=0,\ldots,k$}
    \State $\left(\boldsymbol{\rho}^{k+1}\right)_{i} = \sum_{j=0}^{i}
    \transpose{\left(\boldsymbol{\zeta}^{i}_{0}\right)_{j}}
    \left(\tilde{\boldsymbol{\zeta}}^{k+1}_{0}\right)_{j}
    +
    \transpose{\left(\boldsymbol{\zeta}^{i}_{1}\right)_{j}}
    \left(\tilde{\boldsymbol{\zeta}}^{k+1}_{1}\right)_{j}$
    \State $\tilde{\boldsymbol{\zeta}}^{k+1}_{0} \gets
    \tilde{\boldsymbol{\zeta}}^{k+1}_{0} -
    \begin{bmatrix}
      \boldsymbol{\zeta}^{i}_{0}\\
      0
    \end{bmatrix}\left(\boldsymbol{\rho}^{k+1}\right)_{i}$
    \State $\tilde{\boldsymbol{\zeta}}^{k+1}_{1} \gets
    \tilde{\boldsymbol{\zeta}}^{k+1}_{1} -
    \begin{bmatrix}
      \boldsymbol{\zeta}^{i}_{1}\\
      0
    \end{bmatrix}\left(\boldsymbol{\rho}^{k+1}\right)_{i}$
    \EndFor
    \State $\begin{bmatrix}
      \boldsymbol{\zeta}^{k+1}_{0}\\
      \boldsymbol{\zeta}^{k+1}_{1}
    \end{bmatrix}
    \left(\boldsymbol{\rho}^{k+1}\right)_{k+1} =
    \begin{bmatrix}
      \tilde{\boldsymbol{\zeta}}^{k+1}_{0}\\
      \tilde{\boldsymbol{\zeta}}^{k+1}_{1}
    \end{bmatrix}$
    \Comment {Normalization}
    \If{not root}
    \State Send $\boldsymbol{\rho}^{k+1}$ to parent
    \Else
    \State Send $\boldsymbol{\zeta}^{k+1}_{0}$ and
    $\boldsymbol{\zeta}^{k+1}_{1}$ back to children to start
    back-propagation
    \State broadcast $\boldsymbol{\rho}^{k+1}$ to all processes
    \EndIf
  \end{algorithmic}
\end{algorithm}
\begin{algorithm}[t]
  \caption{Back-Propagation for the Localized Arnoldi Method}
  \label{alg:back-propagation-localized}
  \begin{algorithmic}
    \State Receive $\widetilde{\boldsymbol{\rho}}^{k+1}$ from parent
    \State $\bar{\boldsymbol{\zeta}}^{k+1}_{0} =
    \sum_{i=0}^{k+1} \begin{bmatrix}
      \boldsymbol{\zeta}^{i}_{0}\\
      \boldsymbol{0}
    \end{bmatrix}\left(\widetilde{\boldsymbol{\rho}}^{k+1}\right)_{i}$
    \State $\bar{\boldsymbol{\zeta}}^{k+1}_{1} =
    \sum_{i=0}^{k+1} \begin{bmatrix}
      \boldsymbol{\zeta}^{i}_{1}\\
      \boldsymbol{0}
    \end{bmatrix}\left(\widetilde{\boldsymbol{\rho}}^{k+1}\right)_{i}$
    \If{not leaf}
    \State Send $\bar{\boldsymbol{\zeta}}^{k+1}_{0}$ and
    $\bar{\boldsymbol{\zeta}}^{k+1}_{1}$ back to children
    \EndIf
  \end{algorithmic}
\end{algorithm}
\begin{algorithm}
  \caption{Localized Arnoldi Method}
  \label{alg:localized-arnoldi}
  \begin{algorithmic}
    \State $\widetilde{\boldsymbol{\zeta}}^{0}_{p}, P_{p}^{0} =
    \normalizer[X^{0}_{p}]{\SubA}$ \Comment{local}
    \State Compute local and global $R$-factor $\boldsymbol{\zeta}^{0}_{p},
    \boldsymbol{\rho}^0 \in \SubA^{1}$ from $\widetilde{\boldsymbol{\zeta}}^{0}$
    using Alg. \ref{alg:reduction-localized}
    and \ref{alg:back-propagation-localized}
    \For{$k=0,\ldots,k_{\max}$}
    \State $X_{p}^{k} = \sum_{j=0}^{k}
    P_{p}^{j}\left(\boldsymbol{\zeta}^{k}_{p}\right)_{j}$
    \Comment assemble global orth. vector
    \State $P^{k+1} = AX^{k}$ \Comment requires neighbor communication
    \For{$j=0,\ldots,i$} \Comment local modified Gram-Schmidt
    \State $\left(\widetilde{\boldsymbol{\zeta}}^{k+1}_{p}\right)_{j} =
    \sproduct[\SubA]{P^{j}_{p}}{P^{k+1}_{p}}$
    \State $P^{k+1}_{p} \gets P^{k+1}_{p} - P^{j}_{p}\left(\widetilde{\boldsymbol{\zeta}}^{k+1}_{p}\right)_{j}$
    \EndFor
    \State $\widetilde{\boldsymbol{\zeta}}^{k+1}_{p}, P_{p}^{k+1} \gets
    \normalizer[P_{p}^{k+1}]{\SubA}$
    \State Compute global $R$-factors $\boldsymbol{\zeta}^{k+1}_{j},
    \boldsymbol{\rho}^{k+1} \in \SubA^{k+2}$ from
    $\widetilde{\boldsymbol{\zeta}}^{k+1}$ using Alg. \ref{alg:reduction-localized}
    and \ref{alg:back-propagation-localized}
    \EndFor
  \end{algorithmic}
\end{algorithm}

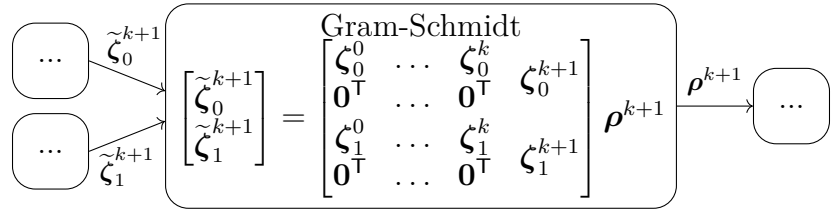
\begin{figure}[tp!]
  \centering
  \begin{tikzpicture}
      \tikzset{node/.style={node distance=1.2cm,draw,rectangle,rounded
          corners=0.3cm,minimum size=1cm}}
    \node [node] (N11) {...};
    \node [below of=N11,node] (N12) {...};

    \node [right=of N12,yshift=0.6cm,node,align=center] (N21) {
      Gram-Schmidt\\
      $\begin{bmatrix}
        \tilde{\boldsymbol{\zeta}}_{0}^{k+1}\\
        \tilde{\boldsymbol{\zeta}}_{1}^{k+1}
      \end{bmatrix}
      =
      \begin{bmatrix}
        \boldsymbol{\zeta}_{0}^{0} & \ldots&\boldsymbol{\zeta}_{0}^{k}
        &\multirow{2}{*}{$\boldsymbol{\zeta}_{0}^{k+1}$}\\
        \transpose{\boldsymbol{0}} & \ldots & \transpose{\boldsymbol{0}}\\
        \boldsymbol{\zeta}_{1}^{0} & \ldots&\boldsymbol{\zeta}_{1}^{k}
        &\multirow{2}{*}{$\boldsymbol{\zeta}_{1}^{k+1}$}\\
        \transpose{\boldsymbol{0}} & \ldots & \transpose{\boldsymbol{0}}\\
      \end{bmatrix}
      \boldsymbol{\rho}^{k+1}
      $};

    \node [right= of N21,node] (root) {...};

    \draw[->] (N11.east) -- ([yshift=0.2cm]N21.west)
    node[midway,above,xshift=0.1cm]
    {\footnotesize{$\tilde{\boldsymbol{\zeta}}^{k+1}_{0}$}};
    \draw[->] (N12.east) -- ([yshift=-0.2cm]N21.west)
    node[midway,below]{\footnotesize{$\tilde{\boldsymbol{\zeta}}^{k+1}_{1}$}};
    \draw[->] (N21.east) -- (root.west)
    node[midway,above]{\footnotesize{$\boldsymbol{\rho}^{k+1}$}};
\end{tikzpicture}
\caption{Reduction operation of the localized Arnoldi method.}
\label{fig:localized_arnoldi_reduction}
\end{figure}

\begin{figure}[tp!]
  \centering
  \begin{tikzpicture}
      \tikzset{node/.style={node distance=1.2cm,draw,rectangle,rounded
          corners=0.3cm,minimum size=1cm}}
    \node [node] (N11) {...};
    \node [below of=N11,node] (N12) {...};

    \node [right=of N12,yshift=0.6cm,node,align=center] (N21) {
      multiplication\\
      $\begin{bmatrix}
        \bar{\boldsymbol{\zeta}}^{k+1}_{1}\\
        \bar{\boldsymbol{\zeta}}^{k+1}_{2}\\
      \end{bmatrix}
      =
      \begin{bmatrix}
        \boldsymbol{\zeta}^{k+1}_{1}\\
        \boldsymbol{\zeta}^{k+1}_{2}\\
      \end{bmatrix}
      \tilde{\boldsymbol{\rho}}^{k+1}
      $};

    \node [right= of N21,node] (root) {...};

    \draw[->] ([yshift=0.2cm]N21.west) -- (N11.east)
    node[midway,above,xshift=0.2cm]
    {\footnotesize{$\bar{\boldsymbol{\zeta}}^{k+1}_{0}$}};
    \draw[->] ([yshift=-0.2cm]N21.west) -- (N12.east)
    node[midway,below]{\footnotesize{$\bar{\boldsymbol{\zeta}}^{k+1}_{1}$}};
    \draw[->] (root.west) -- (N21.east)
    node[midway,above]{\footnotesize{$\tilde{\boldsymbol{\rho}}^{k+1}$}};
\end{tikzpicture}
\caption{Back-propagation operation of the localized Arnoldi method.}
\label{fig:localized_arnoldi_back_propagation}
\end{figure}
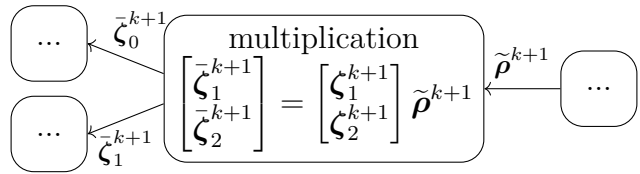

Algorithm~\ref{alg:reduction-localized} and
Algorithm~\ref{alg:back-propagation-localized} show the algorithms for the
reduction and back-propagation procedure.
Figure~\ref{fig:localized_arnoldi_reduction} and
Figure~\ref{fig:localized_arnoldi_back_propagation}
show the reduction and back-propagation operations, respectively.
Unfortunately, the current MPI standard does not allow to define a
custom back-propagation method.
Therefore, we implemented the reduction pattern using Point-to-Point
communications.
The used tree can be either a binary tree using the MPI rank number, or can be
deduced from a \cpp{MPI_Allreduce} pattern, where we use a custom reduction
function to deduce the tree.
A prototype implementation can be found in the gitlab repository of
the author\footnote{\url{https://zivgitlab.uni-muenster.de/n_drei02/tsqr_communication_pattern}}.

Algorithm \ref{alg:localized-arnoldi} shows the localized Arnoldi
method.
It uses the modified Gram-Schmidt method to compute the local
orthogonalization.

One small drawback of the methods can be found in the reduction
process.
As it can be seen in Algorithms \ref{alg:reduction-localized}
and \ref{alg:back-propagation-localized}, the arithmetical complexity
increases quadratically with the size $k$ of the basis.
This leads to an effort of $\mathcal{O}\left(k^{2}\right)$
for the reduction operation.
This term could become a bottleneck for large $k$.
Furthermore, the method computes $2(k+1)$ \texttt{BAXPY} and $k+1$
\texttt{BDOT} operations in the $k$th iteration.
This leads to an overall effort of $6kp^{2}q + \mathcal{O}\left(\log(P)k^{2}p^{2}q\right)$.

\begin{table}
  \centering
  \small
  \caption{Comparison of the arithmetical complexity, number of
    messages and stability for different orthogonalization methods in
    the BGMRes method.}
  \begin{tabular}{lccc}\toprule
    orthogonalization & arith.\ complexity & messages & stable\\
    \midrule
    classical & $4kp^{2}q\frac{\dimA}{P}$ & $\mathcal{O}(\log(P))$ & no\\
    classical(2) & $8kp^{2}q\frac{\dimA}{P}$ & $\mathcal{O}(\log(P))$ & yes\\
    modified & $4kp^{2}q\frac{\dimA}{P}$ & $\mathcal{O}(k\log(P))$ & yes\\
    pipelined($r$) hybrid & $4kp^{2}q\frac{\dimA}{P}$ & $\mathcal{O}(k\log(P))$ (overl.) & depend on $r$\\
    localized & $6kp^{2}q\frac{\dimA}{P} + \mathcal{O}\left(\log(P)k^{2}p^{2}q\right)$ & $\mathcal{O}(\log(P))$ & yes\\
    \bottomrule
  \end{tabular}
  \label{tab:bgmres_comparison}
\end{table} Table \ref{tab:bgmres_comparison} compares the localized Arnoldi
method with the other methods developed in this chapter.
Note that the localized Arnoldi strategy only proceeds one reduction
pattern and the back-propagation, while the classical orthogonalization method
proceeds two reductions and the broadcast of the result.
From the theoretical values in Table \ref{tab:bgmres_comparison}, the localized
method is competitive to the classical(2) method.
However, the methods differ in the arithmetical complexity and the number of
messages, as it only needs one reduction and back-propagation instead of four
reductions and four broadcasts.
For a BGMRes method with fewer iterations, or a low restart parameter, the
localized method would perform better, while for large $k$ the classical($2$)
orthogonalization method or the pipelined version would perform better.

\section{Numerical Experiments}

In this chapter, we introduced several orthogonalization methods for the BGMRes
method.
All resulting BGMRes methods are mathematically equivalent but differ in the
numerical stability and communication efforts.

\begin{figure}[tbp!]
  \centering
  \input{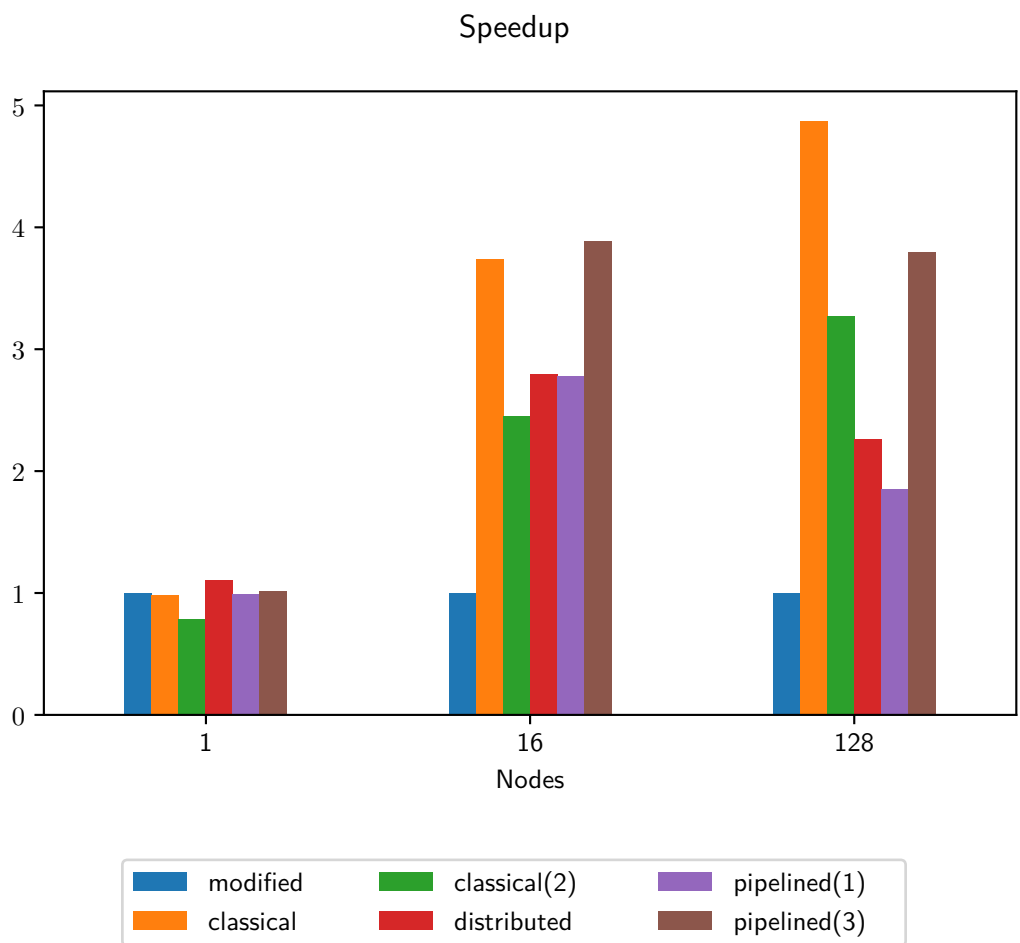}
  \caption[Speedup of runtime for different BGMRes variants compared to the
modified method.]{Speedup of runtime for different BGMRes variants compared to
the modified method.
    Colors decode the orthogonalization method.}
  \label{fig:parallel_gmres}
\end{figure}
To compare all these methods, we run the same test problem from Figure
\ref{fig:gmres_hybrid_vs_global}, which is the \texttt{Simon/raefsky3}
matrix.
The result is shown in Figure~\ref{fig:parallel_gmres}.
We used a block size of $p=4$ and restart the BGMRes method after $100$ iterations.
For a fixed number of processors $P$ all methods need almost the same
number of iterations.
Hence, stability is not an issue here.
As in the previous benchmark, we can observe the same three regimes.
In the sequential case, all methods take almost the same time.
Just the classical($2$) method is slightly slower as it
proceeds the orthogonalization procedure twice, which introduces overhead.
Running the test on a medium number of processors shows a significant speedup
of all the communication-optimized methods.
On $P=16$ nodes, the pipelined($3$) and classical methods are
fastest.
They reach a speedup of almost $4\times$ over the modified Gram-Schmidt method.
The other methods, classical($2$), localized Arnoldi and pipelined($1$) only
reach a speedup of almost $3\times$.

On the scaling limit ($P=128$) the classical method is clearly the fastest,
followed by the pipelined($3$) method.
The classical($2$) method performs significantly better than the localized
methods in this regime.
This could be due to the suboptimal scaling of the self-implemented reduction
operation or the $k^{2}$ term in the reduction operation.
Furthermore, the reduction pattern is not optimized for network hardware yet.
We assume that the performance of the localized method could be significantly
improved if this communication pattern would be supported by the MPI implementation.
A better implementation and quantification of the performance is an objective
of future work.

Overall, we see that the optimization of the orthogonalization
strategy can have significant impact on the performance of the BGMRes
method.


\chapter{Pipelined Block BiCGStab Method}
\label{chap:communicationbicgstab}

In this chapter, we apply the same techniques to the BBiCGStab method, as we used
for the BCG method in Chapter~\ref{chap:communicationcg}.
As the introduction of further recursions introduces sources of numerical
instability, we do not pursuit it very aggressively.
Numerical experiments showed that the mitigation of these instabilities is not
trivial and is therefore left for further work.


\begin{algorithm}[p!]
  \caption{Pipelined BBiCGStab with Adaptive Residual Re-Orthonormalization}
  \label{alg:BlockBiCGStab_with_reortho_and_comm_opti}
  \begin{algorithmic}
    \State $R^{0} = B-AX^{0}$
    \If{$\eta > 0$}
    \State $\hat{R}^{0},\sigma^{0} = \normalizer[R^{0}]{\SubA}$
    \Else
    \State $\sigma^{0} = \Identity$
    \EndIf
    \State $\hat{P}^{0} = \inverse{\Prec}\hat{R}^{0}$
    \State $\chi^{0} = \sproduct[\SubA]{\hat{R}^{0}}{\hat{R}^{0}}$
    \State Choose $\transinv{\Prec}\tilde{R}^{0}$
    (e.g.\ $\transinv{\Prec}\tilde{R}^{0} = \hat{P}^{0}$)
    \State $\hat{V}^{0} = \hat{P}^{0}$
    \For{$k=0,\ldots$}
    \State $\hat{Q}^{k} = A\hat{P}^{k}$
    \State \initcommunication{lambda}$\hat{\lambda}^{k} =
    \inverse{\left(\sproduct[\SubA]{\transinv{\Prec}\tilde{R}^{0}}{\hat{Q}^{k}}\right)}
    \sproduct[\SubA]{\transinv{\Prec}\tilde{R}^{0}}{\hat{R}^{k}}$
    \State $\hat{Z}^{k} = \inverse{\Prec}\hat{Q}^{k}$
    \State \finalizecommunication{lambda}$\hat{S}^{k} = \hat{R}^{k}
    -\hat{Q}^{k}\hat{\lambda}^{k}$
    \State $X^{k+\frac12} = X^{k} + \hat{P}\lambda^{k}\sigma^{k}$
    \If {$\eta\kappa_{D}\left(\chi^{k}\right)
      > \sqrt{\machineeps}$}
    \State $\hat{S}^k,\gamma^{k} \gets \normalizer[\hat{S}^{k}]{\SubA}$
    \State $\sigma^{k+1} = \gamma^{k}\sigma^{k}$
    \State $\hat{T} = \inverse{\Prec}\hat{S}$
    \Else
    \State $\hat{T}^{k} =
    \hat{V}^{k}-\hat{Z}^{k}\hat{\lambda^{k}}$
    \EndIf
    \State \initcommunication{chi}$\chi^{k+1} = \sproduct[\SubA]{\hat{S}^{k}}{\hat{S}^{k}}$
    \State $\hat{U}^{k} = A\hat{T}^{k}$
    \State \initcommunication{omega}$\hat{\omega}^{k} =
    \frac{\sproduct[F]{\hat{U}^{k}}{\hat{S}^{k}}}{\sproduct[F]{\hat{U}^{k}}{\hat{U}^{k}}}$
    \State \initcommunication{beta}$\hat{\beta}^{k} = -\inverse{\left(\sproduct[\SubA]{\transinv{\Prec}\tilde{R}^{0}}{\hat{Q}^{k}}\right)}
    \sproduct[\SubA]{\transinv{\Prec}\tilde{R}^{0}}{\hat{U}^{k}}$
    \State $\hat{W}^{k}  = \inverse{\Prec}U^{k}$
    \State \finalizecommunication{chi}break if $\|\chi^{k+1}\|$
    \State \finalizecommunication{omega}$X^{k+1} = X^{k+\frac12} + \hat{\omega}^{k}\hat{T}\sigma^{k}$
    \State $\hat{R}^{k+1} =
    \hat{S}^{k}-\hat{\omega}^{k}\hat{U}^{k}$
    \State $\hat{V}^{k+1} = \hat{T}^{k} - \hat{\omega}\hat{W}^{k}$
    \State \finalizecommunication{beta}$\hat{P}^{k+1} = \hat{V}^{k+1} + \left(\hat{P}^{k} - \hat{\omega}^{k}\hat{Z^{k}}\right)\hat{\beta}^{k}$
    \EndFor
  \end{algorithmic}
  \drawcommunication{lambda/0.6em,omega/1.2em,beta/0.6em,chi/1.8em}
\end{algorithm}

A study of the BBiCGStab method (Algorithm
\ref{alg:BlockBiCGStab_with_reortho}) shows that it is already in a good state
to overlap the communication with the preconditioner application without
changing a lot.
In every iteration two applications of the operator and two applications of the
preconditioner are performed.
Global communication is needed to compute the coefficients $\lambda^{k}, \beta^{k}$ and
$\omega^{k}$.
The global communication for $\lambda^{k}$ can already be overlapped with the
computation of $\hat{Z}^{k} = \inverse{\Prec}\hat{Q}^{k}$.
Furthermore, the computation of $\omega^{k}$ and $\beta^{k}$ can be fused.
To overlap that communication, we introduce a new variable $\hat{W}^{k} =
\inverse{\Prec}\hat{U}^{k}$ that can be computed during the communication, and
then be used to update $\hat{V}^{k}$ recursively, i.e.
\begin{align}
  \hat{V}^{k} = \inverse{\Prec}\hat{R}^{k} = \hat{T}^{k} - \hat{\omega}^{k}\hat{W}^{k}.
\end{align}
Together with this communication, we can also compute the product
$\sproduct[\SubA]{\hat{S}^{k}}{\hat{S}^{k}}$, which is then used to determine the
residual norm for the break criterion and for the decision, whether we need to
re-orthogonalize the residual in the next iteration.
In that way, we overlap all communication with the two preconditioner
applications.
The communication optimized variant of the BBiCGStab method can be found in
Algorithm \ref{alg:BlockBiCGStab_with_reortho_and_comm_opti}.
The arrows indicate which operations can overlap.

A further optimization would be to precompute the operator application, such
that it can also be computed during the communication.
An attempt by the author failed due to instabilities.
The analysis and mitigation of these instabilities is an objective of further
work.
It could be based on the work of
\citeauthor{cools2017communication}~\cite{cools2017communication,cools2019analyzing},
who analyzed the round-off errors for the pipelined non-block BiCGStab method.

\section{Numerical Experiments}

Although we do not spend a lot of effort in the optimization of the BBiCGStab
method, we want to quantify the performance in the parallel case.
We compare the BBiCGStab method (Algorithm \ref{alg:BlockBiCGStab_with_reortho})
with its optimized version (Algorithm
\ref{alg:BlockBiCGStab_with_reortho_and_comm_opti}).
For that, we use $s=4$ right-hand sides and compare the results for the parallel
($\SubA[P]$) and block ($\SubA[B]$) method.
The test problem is the same as in Section \ref{sec:blocknumericbicgstab}, i.e.\
the \texttt{raefsky3} matrix.
We use a ILU($0$) preconditioner with an additive Schwarz decomposition for the
parallelization.

\begin{figure}
  \centering
  \input{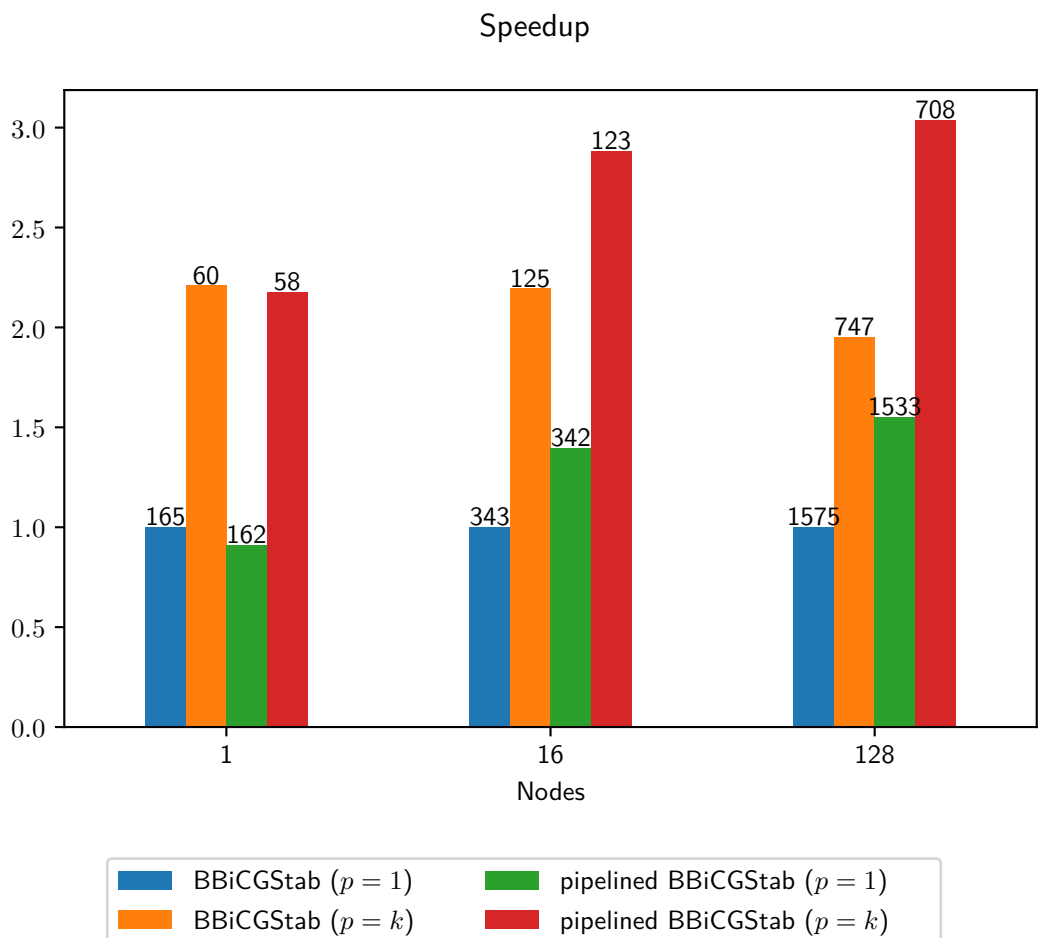}
  \caption[Speedup of runtime for different BiCGStab variants.]{Speedup of runtime for different BiCGStab variants.
    Numbers indicate the number of iterations that are needed to reduce the
residual by a factor of \num{e-7}.}
  \label{fig:parallel_bicgstab}
\end{figure}
Figure \ref{fig:parallel_bicgstab} shows the speedup of the different versions
compared to the parallel BBiCGStab method.
The numbers at the top of the bars show the number of iterations the
method needs to reduce the residual by a factor of \num{e-7}.
As expected in the sequential case the block method is faster than the
parallel method, but the optimization does not yield any remarkable benefits.
We see that the pipelined version is slightly slower in this regime
as it introduces arithmetical overhead.
However, the optimization of the algorithm does not affect the stability.
The pipelined BBiCGStab method actually needs slightly fewer iterations than
its not optimized counterparts.

In the parallel cases, the optimized versions are superior to the
non-optimized ones.
They are approximately $50\%$ faster.
This is due to the fewer global synchronizations they perform.
The increase of iterations on larger nodes is due to the weaker preconditioning,
as the domain is decomposed into more Schwarz domains.
Also in this case, the pipelined BBiCGStab method needs fewer iterations than
its BBiCGStab counterpart.



%
\bookmarksetup{startatroot}
\cleardoublepage
\chapter{Summary and Outlook}
\label{chap:summary}
\epigraph{\textit{Nothing in life is to be feared, it is only to be
    understood.}}{\scshape Marie Curie}

Finally, we summarize the achievements of this thesis and give an
outlook for future work and ideas how to transfer some of our approaches
to other contexts.

In Chapter \ref{chap:blockkrylovframework} we reviewed the block Krylov
framework by
\citeauthor[]{frommer2017block}~\cite[]{frommer2017block,frommer2019block}
and presented the novel block-global method, which turned out to be
irrelevant for practice as we have seen later.
We analyzed the framework and its building blocks concerning the
performance on modern hardware.
In particular, we considered the applicability of SIMD instructions
and their arithmetical intensity.
The advantage of the block Krylov framework over classical block
Krylov methods is that the blocking overhead could be balanced even
for a fixed large number of right-hand sides.
We saw that the vector update and inner product kernels perform with
constant time per right-hand side up to a certain blocking parameter.
This means that the faster convergence rate comes for free in this setting.

In the following chapters, we formulated the block variants of the CG, GMRes and
BiCGStab methods based on that block Krylov framework.
For the block CG method, we provided a convergence analysis which gave insights
into the behavior of the different block Krylov variants.
In particular, we saw why the novel block-global method is inferior to the
block-parallel method.
Furthermore, we introduced a novel stabilization strategy which stabilizes the
method and avoids the process of deflation.
Deflation, which is applied in a lot of other works in the literature, is
improper in our context, as we stick to the number of lanes given by the SIMD
interface.
For the stabilization strategy, we decided to orthonormalize the residual with
respect to the euclidean block inner product.
\citeauthor[]{dubrulle2001retooling}~\cite[]{dubrulle2001retooling} suggested in
his work to orthonormalize with respect to the inner product that is given by
the preconditioner.
This would simplify the algorithm but disqualifies the Householder algorithm for
orthonormalization.
The analysis and implementation of this approach is left for future work.
All the theoretical findings about the convergence rate and the benefits of the
stabilization strategy are supported by numerical tests.

In Chapter \ref{chap:blockgmres}, we addressed the block GMRes method.
Like for the block CG method, we formulated the method in the context
of the block Krylov framework.
To do so, we used a generalization of the Givens rotations to
triangulate the Hessenberg matrix.
In a numerical experiment, we observed that the convergence rates of
the different block Krylov variants are similarly connected as for the
CG method.

As a last block Krylov method, we considered the block BiCGStab method in
Chapter~\ref{chap:blockbicgstab}.
We provided a formulation based on orthogonal polynomials in the context of the
block Krylov framework.
Furthermore, we applied a similar residual re-orthonormalization strategy as for the
block CG method.
Numerical tests showed that also for the block BiCGStab method, the block-global
method performs inferior.
In future work, it would be interesting to investigate whether the stabilization
coefficient, which we have chosen as a scalar in our method, could be chosen as
an element of the *-subalgebra $\SubA$.
Besides that, higher level stabilization BiCGStab methods like presented by
\citeauthor*[]{saito2014development}~\cite[]{saito2014development} could be
applied in our context.

The second part of the thesis is about the optimization of the block
Krylov methods with respect to communication in distributed memory
systems.
For that, we discussed the conditions and challenges in distributed
systems and created a benchmark to measure the amount of time that can
be used for computation while collective communication is active.
Furthermore, we reviewed the TSQR algorithm which is a key building
block for block Krylov methods on distributed systems.

In the following chapter, we presented five variants of the block CG method with
different properties regarding communication by applying the approaches of
\citeauthor[]{ghysels2014hiding}~\cite[]{ghysels2014hiding} and
\citeauthor[]{gropp2010update}~\cite[]{gropp2010update}.
In the numerical tests on a medium size supercomputer, we observed three regimes
in that the methods behave differently.
This coincides with our theoretical expectations.
In the communication dominated regime, we achieved good speedups compared to
the default block CG method.

To optimize the block GMRes method with respect to communication, we
considered, in contrast to many approaches in the literature, the
orthogonalization method in the Arnoldi process.
This makes sense in our context, as the block vector update and the
inner block product are rather expensive so that overlapping it with
computation already suffice to hide the whole communication
costs.
Moreover, we presented a novel method for the orthogonalization
process that is based on the TSQR algorithm and allows an
orthogonalization with only one global reduction communication.
It would be interesting to examine how this method could be combined with other
communication-avoiding approaches, like $s$-step GMRes.
Another interesting subject would be to use this method for example in
the MINRes method.
In principle, the applicability of the method is not restricted to block
Krylov methods.
In the numerical experiments, we compared the different
orthogonalization methods and observed significant speed up compared
to the standard block GMRes method.

We reviewed shortly the amenability of pipelining
techniques to the block BiCGStab method in the last chapter of the thesis.
We found that some approaches can easily be applied but others
introduce too much numerical instability such that we did not further
pursuit them.
Thus, this would be an interesting subject for further work.
However, the approaches and optimization that we have applied already
yield a pretty good speedup in our numerical experiments.

In summary, we presented tailored methods to solve large sparse systems
on modern super computing hardware.
We proved their advantages in numerical tests on both, the node level
and on a large distributed system.
The author aims for the integration of the methods into the \dune[ISTL]
module and intends to present a merge request soon after submitting
this thesis to make the results of the thesis available to the community.

In the future, we want to provide a better comparison of the new stabilization
methods with deflation strategies, both experimentally and analytically.
We hope that this could give better insights for choosing the
re-orthonormalization parameter.
Moreover, we want to improve the implementation of the reduction and
back-propagation communication pattern and make it available for the community,
as it could be a generic building block for more methods like the localized
Arnoldi method.

In addition, we want to investigate how we could make the methods usable for a
broader range of problems.
One idea is to apply them in ODE solvers and compute multiple time steps
simultaneously.
This would lead formally to an approach like in parallel-in-time methods, which are
currently very popular in the scientific computing research community.
However, the gain of parallelism could also be used to apply block Krylov
methods instead of distributing it over the nodes.
This would not only decrease the inter-node communication, but also improve the
other aspects discussed in this work.

\appendix%

%

\cleardoublepage%
\phantomsection
\addcontentsline{toc}{chapter}{Bibliography}%
\nocite{*}
\printbibliography[category=cited]
%

%
\end{document}